\documentclass[11pt,reqno]{amsart}

\usepackage{amsthm,amssymb,amsmath,mathtools,bm}
\usepackage[colorlinks,linktocpage,linkcolor=blue]{hyperref}
\usepackage{todonotes,csquotes,lineno,graphicx}

\theoremstyle{definition}
\newtheorem{definition}{Definition}[section]
\theoremstyle{plain}
\newtheorem{theorem}[definition]{Theorem}
\newtheorem{lemma}{Lemma}[section]
\newtheorem{proposition}[definition]{Proposition}

\newtheorem{remark}[definition]{Remark}
\newtheorem{assumption}[definition]{Assumption}
\newtheorem{example}[definition]{Example}
\usepackage{geometry}
\numberwithin{equation}{section}

\allowdisplaybreaks[4]

\newcommand\del[1]{}

\newcommand{\V}{V}
\newcommand{\U}{U}
\newcommand{\E}{\mathbb E}
\newcommand{\R}{\mathbb R}
\newcommand{\1}{\mathbf 1}
\newcommand{\HS}{\mathcal L_2}
\newcommand{\cL}{\mathcal L}
\newcommand{\HH}{\mathcal H}
\newcommand{\deltaH}{\boldsymbol\delta}

\begin{document}

\title[Asymptotic error for SPDE driven by fBm]{Asymptotic error of the accelerated exponential Euler scheme for semilinear parabolic SPDEs driven by fractional Brownian motion}

\author{Diancong Jin}
\address{School of Mathematics and Statistics, Huazhong University of Science and Technology, Wuhan 430074, China; Hubei Key Laboratory of Engineering Modeling and Scientific Computing, Huazhong University of Science and Technology, Wuhan 430074, China}
\email{jindc@hust.edu.cn}

\author{Minoo Kamrani}
\address{School of Mathematics, Statistics and Computer Science, College of Science, University of Tehran, Tehran, Iran}
\email{m.kamrani@ut.ac.ir}

\author{Xu Wang}
\address{SKLMS, ICMSEC, Academy of Mathematics and Systems Science, Chinese Academy of Sciences, Beijing 100190, China, and School of Mathematical Sciences, University of Chinese Academy of Sciences, Beijing 100049, China}
\email{wangxu@lsec.cc.ac.cn}

\thanks{This work is supported by the National key R\&D Program of China (No. 2024YFA1015900), the National Natural Science Foundation of China (No. 12288201 and No. 12671472),  the CAS Project for Young Scientists in Basic Research (YSBR-087),  the Iran National Science Foundation (INSF) (No. 40507840), and the Fundamental Research Funds for the Central Universities (No. 2026BRSXB006).}


\begin{abstract}
We study the strong convergence and asymptotic error of the accelerated exponential Euler scheme for semilinear parabolic stochastic partial differential equations driven by additive cylindrical fractional Brownian motion with Hurst parameter $H\in(\frac12,1)$. The drift is globally Lipschitz and G\^ateaux differentiable, with a second directional derivative taking values in a negative-order space. Under a weighted Hilbert--Schmidt condition on the noise and suitable derivative bounds, we establish strong error estimates that reflect the interaction between noise regularity and drift regularity. In a regime where the strong order is one, the normalized error $m(X^m-X)$ converges in probability in $C([0,T];\V)$ to the unique continuous solution of a random linear integral equation. Its stochastic forcing is a divergence integral plus an explicit covariance correction. The proof combines a G\^ateaux chain rule for Malliavin derivatives, uniform spectral approximations, and time increment estimates. Localization gives the result for independent random initial values with finite moments, together with lower-moment convergence at each fixed time.
\end{abstract}

\maketitle


\section{Introduction}
Strong convergence rates quantify the magnitude of numerical errors, while asymptotic error analysis identifies the random process governing their leading behavior as the time step tends to zero. For stochastic partial differential equations (SPDEs), this requires understanding how temporal discretization interacts with an unbounded evolution operator and infinitely many noise modes. In this paper, we study the accelerated exponential Euler scheme for semilinear parabolic SPDEs driven by additive cylindrical fractional Brownian motion with Hurst parameter $H\in(\frac12,1)$. Our aim is to identify the leading error as a process on the entire time interval and to explain how the temporal correlation of fractional noise enters its limiting equation.

Limit theorems for numerical errors have a substantial history in finite-dimensional stochastic analysis. Jacod and Protter \cite{JP98} studied normalized Euler errors for equations driven by semimartingales. For fractional Brownian motion, Neuenkirch and Nourdin \cite{NN07} obtained limits for normalized errors of scalar Euler approximations, and Hu, Liu, and Nualart \cite{HLN16} showed how the normalization and limiting behavior of modified Euler schemes depend on the Hurst parameter. For additive fractional noise, Zhou, Hu, and Liu \cite{ZHL23} established first-order convergence and an asymptotic error limit for the backward Euler method, while Chen, Shen, and Yu \cite{CSY24} proved a first-order mean-square estimate for the exponential Euler method for stiff stochastic differential equations. Extending such an analysis to SPDEs requires estimates that remain valid in the presence of an unbounded linear operator in an infinite-dimensional state space.

For parabolic SPDEs, exponential methods exploit the semigroup generated by the linear operator. The accelerated approximation introduced by Jentzen and Kloeden \cite{JK09} retains the stochastic convolution, and its strong convergence for additive Wiener noise was further studied by Wang and Qi \cite{WQ15}. For fractional noise, Wang, Qi, and Jiang \cite{WQJ17} established sharp regularity estimates and convergence rates for spectral Galerkin and linear implicit Euler approximations. Hong and Huang \cite{HH20} proved strong order one for an exponential integrator for the stochastic heat equation by combining Malliavin calculus with semigroup estimates. Cai, Gan, and Wang \cite{CGW24} investigated weak convergence of a fully discrete spectral Galerkin and exponential Euler approximation. These results establish approximation rates but leave the leading error process of the temporal scheme studied here to be identified.

For additive Wiener noise, Hong, Jin, Wang, and Yang \cite{HJWY24} identified the asymptotic error distribution of the accelerated exponential Euler method through finite-dimensional approximations that are uniform in the temporal discretization. Their result gives convergence in distribution at each fixed time, and the limiting equation involves an additional Wiener process independent of the original noise. For multiplicative Wiener noise, Hong, Jin, and Wang \cite{HJW25} obtained limiting equations involving infinitely many independent Wiener processes. Fractional noise with $H>\frac12$ leads to a different stochastic forcing in the error limit. Its correlated increments and lack of a semimartingale structure also require a different argument for identifying that limit.

We consider the stochastic evolution equation \eqref{eq:SPDE} and its accelerated exponential Euler approximation $X^m$ defined in \eqref{eq:Xm}, where $m$ is the number of time steps. The method retains the stochastic convolution and freezes the nonlinear drift on each time interval. We first establish strong error estimates in $L^p(\Omega;\V)$ for every finite $p\ge1$, with rates that account for the spatial regularity of the noise and the action of the nonlinear drift. These estimates hold uniformly over time for the continuous interpolation and include a regime of first-order convergence (see Theorem \ref{tm:rate}). They provide the moment bounds and remainder estimates needed to study the normalized error, including cases where the noise coefficient itself is not Hilbert--Schmidt.

Our main result, Theorem \ref{tm:asymptotic-error}, identifies the limit in this first-order regime:
\[
m(X^m-X)\to U \quad\text{in probability in }C([0,T];\V).
\]
The process $U$ is the unique continuous solution of the random linear integral equation \eqref{eq:limit-U}. This gives an asymptotic expansion of the numerical error uniformly over the whole time interval. The limiting forcing contains both the accumulated drift error and a fractional-noise contribution. The latter is represented by a divergence integral together with an explicit covariance correction, which records the dependence between the random integrand and the noise. The limit is constructed from the original fractional noise and initial value.  Under the initial moment assumptions of the theorem, convergence also holds in lower moments at each fixed time.

The first difficulty is to control the stochastic terms in the error expansion while retaining the effect of correlations between different time steps. Estimates based only on the temporal H\"older regularity of the solution do not recover the first-order rate or identify the limiting noise contribution. At the same time, the unbounded operator and limited spatial regularity produce singular semigroup kernels. We derive Malliavin derivative estimates for the exact and numerical solutions and combine semigroup smoothing with divergence estimates. The divergence product formula separates the stochastic error into a divergence term and a covariance correction. Controlling both terms is essential for obtaining the strong estimates and identifying the limiting equation.

A second difficulty is to connect finite-dimensional averaging with the infinite-dimensional error process. An argument at a fixed spectral dimension is insufficient unless the approximation error can be controlled uniformly in the number of time steps. We construct auxiliary error processes and establish such uniform spectral approximation estimates. For each fixed spatial dimension, Young integration and an averaging argument identify the limit of the accumulated fractional-noise increments. Passing to the infinite-dimensional limit then requires convergence of both the divergence term and its covariance correction, as well as control of the deterministic integral terms. This combines a pathwise interpretation after projection with a Malliavin-calculus representation suited to the infinite-dimensional setting.

The final difficulty is to obtain convergence of the entire error path. Bounds that are uniform in time on fixed-time moments do not by themselves imply convergence in $C([0,T];\V)$, and the spectral approximation must remain uniform when time increments are considered. We prove uniform H\"older moment estimates for the auxiliary processes and the remainders. Together with a finite time-grid argument, these estimates allow us to pass from convergence at fixed times to convergence in path space and to remove the spectral projection. Localization extends the path-space convergence to independent random initial values under finite moment assumptions.

The paper is organized as follows. Section~\ref{sec:2} introduces the setting and establishes solution regularity. Section~\ref{sec:3} develops the Malliavin derivative estimates and proves the strong convergence rates. Section~\ref{sec:asymptotic} identifies the limiting equation and establishes convergence of the normalized error in path space.

\section{Preliminaries and assumptions}
\label{sec:2}
\subsection{Mathematical setting and assumptions}
Let $\V$ and $\U$ be separable real Hilbert spaces equipped with inner products $\langle\cdot,\cdot\rangle$ and $\langle\cdot,\cdot\rangle_U$, and the corresponding induced norms $\|\cdot\|$ and $\|\cdot\|_U$, respectively.

\begin{assumption}\label{as:A}
Let $-A:Dom(-A)\subset\V\to\V$ be a linear, positive self-adjoint operator with compact inverse.
\end{assumption}
Then $A$ generates an analytic semigroup $\{E(t)=e^{tA}\}_{t\ge0}$. For $r\in\R$, define the Hilbert space $\V^r:=Dom((-A)^{\frac r2})$ equipped with the norm $\|v\|_r=\|(-A)^{\frac r2}v\|$, with the usual completion when $r<0$. It holds that
\begin{align}
\label{eq:semi1} \|(-A)^\theta E(t)\|_{\cL(\V)}&\le C t^{-\theta},\qquad t>0,\ \theta\ge0, \\ \label{eq:semi2} \|(-A)^{-\theta}(E(t)-I)\|_{\cL(\V)}&\le C t^\theta, \qquad t\ge0,\ 0\le\theta\le1,
\end{align}
where $\mathcal{L}(V)$ denotes the space of bounded linear operators on $V$, $I$ denotes the identity operator on $V$, and the constant $C$ is independent of $t$. Here and throughout, $C$ denotes a generic positive constant whose value may vary from line to line.

Denote by $\mathcal{L}_2(U,V)$ the space of all Hilbert--Schmidt operators from $U$ into $V$ with the scalar product and norm
\[
\langle\Gamma_1,\Gamma_2\rangle_{\mathcal{L}_2(U,V)}:=\sum_{n=1}^{\infty}\langle\Gamma_1g_n,\Gamma_2g_n\rangle,\qquad\|\Gamma\|_{\mathcal{L}_2(U,V)} := \left( \sum_{n=1}^{\infty}\|\Gamma g_n\|^2 \right)^{\frac12}
\]
with $\{g_n\}_{n\in\mathbb N}$ being a complete orthonormal basis of $U$. For $\Gamma\in\mathcal{L}(V)$ and $\Gamma_1,\Gamma_2\in\mathcal{L}_2(U,V)$, the following inequalities hold:
\begin{align*}
\left| \langle \Gamma_1,\Gamma_2\rangle_{\mathcal{L}_2(U,V)} \right| &\leq \|\Gamma_1\|_{\mathcal{L}_2(U,V)} \|\Gamma_2\|_{\mathcal{L}_2(U,V)},\\ \|\Gamma\Gamma_1\|_{\mathcal{L}_2(U,V)} &\leq \|\Gamma\|_{\mathcal{L}(V)} \|\Gamma_1\|_{\mathcal{L}_2(U,V)}.
\end{align*}

Let $W^H=\{W^H(t)\}_{t\in[0,T]}$ be a cylindrical fractional Brownian motion on $\U$ with Hurst parameter $H\in(\frac12,1)$ given by $W^H(t)=\sum_{n=1}^{\infty}g_n \beta_n^H(t)$, where $\{\beta_n^H(t)\}_{n\in \mathbb N}$ is a sequence of independent real-valued fractional Brownian motions each with the Hurst parameter $H\in (\frac 12,1)$. Its canonical Hilbert space $\HH$ is the completion of step functions with respect to the scalar product $\langle\mathbf1_{[0,s]},\mathbf1_{[0,t]}\rangle_{\HH} :=\tfrac12(s^{2H}+t^{2H}-|t-s|^{2H})$. For $f,g\in\HH$, the inner product is given by
\begin{equation}
\label{eq:HH} \langle f,g\rangle_{\HH} =\alpha_H\int_0^T\int_0^T f(u)g(v)|u-v|^{2H-2} du dv, \qquad \alpha_H:=H(2H-1)
\end{equation}
provided that the double integral is absolutely convergent (cf. \cite[Section 5.1.3]{N06}).

Consider the following SPDE
\begin{equation}
\label{eq:SPDE} \left\{
\begin{aligned}
dX(t)&=(AX(t)+F(X(t)))dt+\Phi dW^H(t),\qquad t\in(0,T],\\ X(0)&=X_0,
\end{aligned}
\right.
\end{equation}
whose mild solution, if exists, admits the following form
\begin{equation}
\label{eq:mild} X(t)=E(t)X_0+\int_{0}^{t}E(t-s)F(X(s)) ds+\int_{0}^{t} E(t-s) \Phi dW^H(s),\qquad \mathbb{P}\text{-a.s.}
\end{equation}
for any $t\in(0,T]$. Hereafter, we omit the qualifier ``$\mathbb{P}$-a.s." whenever no confusion arises.

The nonlinear term $F$ and the operator $\Phi$ satisfy the following assumptions.

\begin{assumption}\label{as:Phi}
Let $\Phi:\U\to \V$ satisfy $\|(-A)^{\frac{\beta-1}{2}}\Phi\|_{\mathcal{L}_2(\U,\V)}<\infty$ for some $\beta\in(1-2H,1]$.
\end{assumption}

For convenience, set $\rho:=2H+\beta-1$, where $\beta$ is as in Assumption \ref{as:Phi}. Then $\rho\in(0,2)$.

\begin{assumption}\label{as:F}
The map $F:\V\to\V$ is globally Lipschitz and G\^ateaux differentiable at every $u\in\V$. There exist constants $\delta\in[0,\rho]$, $\gamma\ge0$, $\eta\in[1,2H)$, $\sigma\in[0,1-\beta]$, and $C>0$ such that the following conditions hold.

The first derivative $DF$ is uniformly bounded and strongly operator continuous on $\V$ satisfying
\begin{equation}
\label{eq:F} \|F(u)\|\le C(1+\|u\|),\qquad \|DF(u)\|_{\cL(\V)}\le C,\qquad u\in\V.
\end{equation}
Moreover, on $\V^\rho$, it also satisfies
\[
\|DF(u)-DF(v)\|_{\cL(\V)}\le C\|u-v\|_\rho,\qquad u,v\in\V^\rho.
\]
For $u\in\V^\rho$, the following negative-norm estimate holds:
\begin{equation}
\label{eq:F1} \|DF(u)v\|_{-\delta}\le C(1+\|u\|_\rho)\|v\|_{-\rho-\gamma},\qquad v\in\V.
\end{equation}

For every $u,v,w\in\V$, the directional derivative
\[
D^2F(u)(v,w):=\lim_{\varepsilon\to0}\frac{DF(u+\varepsilon v)w-DF(u)w}{\varepsilon}
\]
exists in $\V^{-\eta}$. For each $u\in\V$, the map $(v,w)\mapsto D^2F(u)(v,w)$ is bilinear and satisfies
\begin{equation}
\label{eq:F2} \|D^2F(u)(v,w)\|_{-\eta}\le C\|v\|\|w\|_{-\sigma},\qquad v,w\in\V.
\end{equation}
\end{assumption}

By density, \eqref{eq:F1} and \eqref{eq:F2} yield unique continuous extensions
\[
DF(u):\V^{-\rho-\gamma}\to\V^{-\delta},\qquad u\in\V^\rho,
\]
and
\[
D^2F(u):\V\times\V^{-\sigma}\to\V^{-\eta},\qquad u\in\V,
\]
respectively. We use the same notation for these extensions, and estimates \eqref{eq:F1} and \eqref{eq:F2} remain valid on the extended domains.

\begin{remark}
The restriction $\beta\le1$ in Assumption \ref{as:Phi} does not exclude spatially smoother noise. Indeed, if $\|(-A)^{(\beta_*-1)/2}\Phi\|_{\HS(\U,\V)}<\infty$ for some $\beta_*>1$, then
\[
\|\Phi\|_{\HS(\U,\V)}\le \|(-A)^{-(\beta_*-1)/2}\|_{\cL(\V)}\|(-A)^{(\beta_*-1)/2}\Phi\|_{\HS(\U,\V)}<\infty.
\]
Thus our results also apply to such noise by taking $\beta=1$, provided that the drift and initial value satisfy the corresponding assumptions with $\rho=2H$ and $\sigma=0$.
\end{remark}

\subsection{Stochastic convolution}
For $0\le s\le t\le T$, define the stochastic convolution starting from time $s$ by
\[
O_{s,t}=\int_s^t E(t-r) \Phi dW^H(r).
\]
In particular, when $s=0$, we write $O_t:=O_{0,t}$. A direct calculation yields $O_{s,t}=O_t-E(t-s)O_s$.

\begin{lemma}\label{lm:O}
Under Assumptions \ref{as:A} and \ref{as:Phi}, the stochastic convolution $O_{s,t}$ is well-defined in $L^2(\Omega;\V)$ for every $0\le s\le t\le T$. Moreover, for any $\tilde{\delta}\in[0,\rho]$ and $p\in[1,\infty)$, there exists a constant $C>0$ independent of $s$ and $t$, such that
\begin{align}
\label{eq:O} &\|O_t\|_{L^p(\Omega;V^{\tilde{\delta}})}\leq C \| (-A)^{\frac{\beta-1}{2}} \Phi \|_{\mathcal{L}_2(\U,\V)},\qquad t\in[0,T],\\ &\|O_{s,t}\|_{L^p(\Omega;V^{\tilde{\delta}})}\leq C \| (-A)^{\frac{\beta-1}{2}} \Phi \|_{\mathcal{L}_2(\U,\V)}(t-s)^{\frac{\rho-\tilde{\delta}}{2}},\qquad 0\le s\le t\le T.\label{eq:Ost}
\end{align}
For any $\tilde{\sigma}\in[0,1-\beta]$ and $p\in[1,\infty)$, there exists a constant $C>0$ independent of $s$ and $t$, such that
\begin{align}
\label{eq:neg-Ost} \|O_{s,t}\|_{L^p(\Omega;\V^{-\tilde{\sigma}})}+\|O_t - O_s\|_{L^p(\Omega;\V^{-\tilde{\sigma}})} \le C\|(-A)^{\frac{\beta-1}{2}}\Phi \|_{\mathcal{L}_2(\U,\V)}(t-s)^{\frac{\rho+\tilde{\sigma}}{2}},\qquad 0\le s\le t\le T.
\end{align}
\end{lemma}

\begin{proof}
The proof for the case $p=2$ follows from \cite{WQJ17}. Indeed, the well-definedness of $O_{s,t}$ in $L^2(\Omega;\V)$ follows from \cite[Proposition 3.7]{WQJ17} and the identity $O_{s,t}=O_t-E(t-s)O_s$.

For $\tilde{\delta}\in[0,\rho]$, by the spatial regularity estimate \cite[(3.19)]{WQJ17} and the continuous embedding $\V^\rho\hookrightarrow\V^{\tilde{\delta}}$,
\[
\|O_t\|_{L^2(\Omega;\V^{\tilde{\delta}})} \le C\|O_t\|_{L^2(\Omega;\V^\rho)} \le C \|(-A)^{\frac{\beta-1}{2}}\Phi \|_{\mathcal{L}_2(\U,\V)}.
\]
Moreover, the stochastic convolution estimate over $[s,t]$ established in \cite[(3.24)]{WQJ17} gives
\[
\|O_{s,t}\|_{L^2(\Omega;\V^{\tilde{\delta}})} \le C \|(-A)^{\frac{\beta-1}{2}}\Phi \|_{\mathcal{L}_2(\U,\V)} (t-s)^{\frac{\rho-\tilde{\delta}}{2}}.
\]

For $\tilde{\sigma}\in[0,1-\beta]$, applying \eqref{eq:Ost} with $\Phi$ and $\beta$ replaced by $\tilde{\Phi}:=(-A)^{-\frac{\tilde{\sigma}}2}\Phi$ and $\tilde{\beta}:=\beta+\tilde{\sigma}$, respectively, and with $\tilde{\delta}=0$ yields
\[
\|O_{s,t}\|_{L^2(\Omega;\V^{-\tilde{\sigma}})}\le C \|(-A)^{\frac{\beta-1}{2}}\Phi \|_{\mathcal{L}_2(\U,\V)} (t-s)^{\frac{\rho+\tilde{\sigma}}{2}}.
\]
The above replacements satisfy the required assumptions, since $\tilde{\beta}\in(1-2H,1]$ and
\[
(-A)^{\frac{\tilde{\beta}-1}{2}} \tilde{\Phi} = (-A)^{\frac{\beta-1}{2}}\Phi\in\mathcal{L}_2(\U,\V).
\]
Then the identity $O_t-O_s=O_{s,t}+(E(t-s)-I)O_s$, \eqref{eq:semi2}, and \eqref{eq:O} with $\tilde\delta=\rho$, imply
\begin{align*}
&\quad\|O_{s,t}\|_{L^2(\Omega;\V^{-\tilde{\sigma}})} +\|O_t-O_s\|_{L^2(\Omega;\V^{-\tilde{\sigma}})}\\ &\le 2\|O_{s,t}\|_{L^2(\Omega;\V^{-\tilde{\sigma}})}+\|(-A)^{-\frac{\rho+\tilde{\sigma}}2}(E(t-s)-I)\|_{\mathcal L(V)}\|O_s\|_{L^2(\Omega;\V^{\rho})}\\ &\le C \|(-A)^{\frac{\beta-1}{2}}\Phi \|_{\mathcal{L}_2(\U,\V)} (t-s)^{\frac{\rho+\tilde{\sigma}}{2}}.
\end{align*}

Finally, since $\Phi$ is deterministic, the stochastic convolutions and their increments are centered Gaussian random variables in the respective Hilbert spaces. For any centered Gaussian random variable $Z$ in a Hilbert space $\V$, the Kahane--Khintchine inequality holds (cf. \cite[Proposition 2.7]{vN10}):
\begin{align}
\label{eq:KK} \|Z\|_{L^p(\Omega;\V)} \le C_p\|Z\|_{L^2(\Omega;\V)},\qquad p\ge1.
\end{align}
Applying this inequality to the preceding estimates completes the proof.
\end{proof}

\subsection{Well-posedness and regularity}
The problem \eqref{eq:SPDE} is well-posed in the above setting.

\begin{theorem}\label{tm:X}
Let Assumptions \ref{as:A}, \ref{as:Phi}, and \ref{as:F}~\eqref{eq:F} hold. Fix $p\in[1,\infty)$, and suppose that $X_0\in L^p(\Omega;V^{\rho})$ is independent of $W^H$. Then the problem \eqref{eq:SPDE} possesses a unique mild solution $X\in C([0,T];L^p(\Omega;V))$ determined by \eqref{eq:mild}.

For any $\tilde{\delta}\in[0,\rho]$, there exists a constant $C>0$ independent of $s$ and $t$, such that
\begin{align*}
\| X(t) \|_{L^p(\Omega;V^{\rho})} &\leq C \Big( 1 + \| X_0 \|_{L^p(\Omega; V^{\rho})} \Big), \qquad t \in [0, T],\\ \|X(t)-X(s)\|_{L^p(\Omega;V^{\tilde\delta})}&\leq C(1+\|X_0\|_{L^p(\Omega;V^{\rho})})(t-s)^{\frac{\rho-\tilde{\delta}}{2}},\qquad 0\le s\le t\le T.
\end{align*}

For any $\tilde{\sigma}\in[0,1-\beta]$, there exists a constant $C>0$ independent of $s$ and $t$, such that,
\begin{equation*}
\|X(t) - X(s)\|_{L^p(\Omega;\V^{-\tilde{\sigma}})} \le C\bigl(1 + \|X_0\|_{L^p(\Omega;\V^\rho)}\bigr)(t-s)^{\frac{\rho+\tilde{\sigma}}{2}},\qquad 0\le s\le t\le T.
\end{equation*}
\end{theorem}

\begin{proof}
By the semigroup estimates \eqref{eq:semi1} and \eqref{eq:semi2}, the global Lipschitz continuity of $F$ in \eqref{eq:F}, and Lemma \ref{lm:O}, the standard fixed-point argument used in \cite[Theorem 3.5]{WQJ17} yields a unique mild solution $X\in C([0,T];L^p(\Omega;V))$.

The mild formulation \eqref{eq:mild}, the linear growth of $F$ in \eqref{eq:F}, and Gronwall's inequality imply
\[
\sup_{r\in[0,T]} \left( \|X(r)\|_{L^p(\Omega;V)} +\|F(X(r))\|_{L^p(\Omega;V)} \right) \le C\bigl(1+\|X_0\|_{L^p(\Omega;V^\rho)}\bigr).
\]
Thus, by Minkowski's inequality, the semigroup smoothing estimate \eqref{eq:semi1} with $\rho<2$, and Lemma \ref{lm:O},
\[
\begin{aligned}
\|X(t)\|_{L^p(\Omega;V^\rho)} &\le C\|X_0\|_{L^p(\Omega;V^\rho)} +C\int_0^t(t-r)^{-\frac{\rho}2} \|F(X(r))\|_{L^p(\Omega;V)}dr+C \\ &\le C\bigl(1+\|X_0\|_{L^p(\Omega;V^\rho)}\bigr).
\end{aligned}
\]

For $0\le s<t\le T$, the semigroup property gives
\[
X(t)-X(s) =(E(t-s)-I)X(s) +\int_s^tE(t-r)F(X(r))dr+O_{s,t}.
\]
For $\tilde\delta\in[0,\rho]$, the preceding moment bounds, the semigroup estimates, and Lemma \ref{lm:O} yield
\[
\begin{aligned}
&\|X(t)-X(s)\|_{L^p(\Omega;V^{\tilde\delta})} \\ \le&~ C(t-s)^{\frac{\rho-\tilde\delta}2} \|X(s)\|_{L^p(\Omega;V^\rho)} +C\int_s^t(t-r)^{-\frac{\tilde\delta}2} \|F(X(r))\|_{L^p(\Omega;V)} dr +C(t-s)^{\frac{\rho-\tilde\delta}2} \\ \le&~ C\bigl(1+\|X_0\|_{L^p(\Omega;V^\rho)}\bigr) \left((t-s)^{\frac{\rho-\tilde\delta}2} +(t-s)^{1-\frac{\tilde\delta}2}\right) \\ \le&~ C\bigl(1+\|X_0\|_{L^p(\Omega;V^\rho)}\bigr) (t-s)^{\frac{\rho-\tilde\delta}2},
\end{aligned}
\]
since $\tilde\delta\le\rho<2$ and $t-s\le T$.

Finally, for $\tilde\sigma\in[0,1-\beta]$, since $0<\frac{\rho+\tilde\sigma}2\le H<1$, we have
\[
\begin{aligned}
\|(E(t-s)-I)X(s)\|_{L^p(\Omega;V^{-\tilde\sigma})} &\le \|(-A)^{-\frac{\rho+\tilde\sigma}2}(E(t-s)-I)\|_{\mathcal{L}(V)} \|X(s)\|_{L^p(\Omega;V^\rho)} \\ &\le C(t-s)^{\frac{\rho+\tilde\sigma}2} \bigl(1+\|X_0\|_{L^p(\Omega;V^\rho)}\bigr).
\end{aligned}
\]
The boundedness of $(-A)^{-\frac{\tilde\sigma}2}$ and $E(t)$, together with Lemma \ref{lm:O}, therefore gives
\[
\begin{aligned}
\|X(t)-X(s)\|_{L^p(\Omega;V^{-\tilde\sigma})} &\le C\bigl(1+\|X_0\|_{L^p(\Omega;V^\rho)}\bigr) \left((t-s)^{\frac{\rho+\tilde\sigma}2}+(t-s)\right) \\ &\le C\bigl(1+\|X_0\|_{L^p(\Omega;V^\rho)}\bigr) (t-s)^{\frac{\rho+\tilde\sigma}2},
\end{aligned}
\]
which completes the proof.
\end{proof}

\section{Strong convergence and auxiliary estimates}
\label{sec:3}
This section establishes strong convergence rates for the accelerated exponential Euler scheme and provides auxiliary estimates for the subsequent asymptotic error analysis. We first derive Malliavin derivative estimates for the exact and numerical solutions. These estimates are then used to prove strong convergence at the grid points and for the continuous interpolation. Finally, we establish uniform moment bounds and temporal regularity estimates for the numerical solution.

By applying the accelerated exponential Euler scheme to \eqref{eq:SPDE}, we have
\begin{eqnarray*}
X_{k+1}^{m}=e^{h A}X_{k}^{m}+A^{-1}(E(h)-I)F(X_{k}^{m})+\int_{t_k}^{t_{k+1}} E(t_{k+1}-s) \Phi dW^H(s),
\end{eqnarray*}
with the temporal step-size $h:=\frac T m$ for some $m\in \mathbb N^+$ and $t_k=kh$ for $k=0,\cdots,m.$ Its time-continuous interpolation is defined by
\begin{eqnarray}
\label{eq:Xm} X^m(t)=E(t)X_0+\int_0^t E(t-s) F(X^m(\kappa_m(s))) ds+\int_0^t E(t-s) \Phi dW^H(s),\qquad t\in(0,T]
\end{eqnarray}
such that $X^m(t_k)=X^m_k$ for $k=0,\cdots,m$, where $\kappa_m(s):=\lfloor\frac s h\rfloor h$ with $\lfloor\cdot\rfloor$ denoting the floor function.

\subsection{Estimates of Malliavin derivatives}
Let $D$ denote the Malliavin derivative associated with the isonormal Gaussian process indexed by $\HH\otimes\U$, and let $\deltaH$ denote its $L^2$-adjoint, called the divergence operator. For an orthonormal basis $\{g_n\}_{n\ge1}$ of $\U$, we write $D_u^nX(t):=(D_uX(t))g_n$. The same convention applies to $X^m$ and other Malliavin-differentiable random variables.

For $1\le p<\infty$, let $\mathbb D^{1,p}(\V)$ denote the completion of smooth $\V$-valued cylindrical random variables under the norm
\[
\|Z\|_{\mathbb D^{1,p}(\V)} := \left( \mathbb E\left[\|Z\|^{p}\right] + \mathbb E\left[\|DZ\|_{\HH\otimes\mathcal{L}_2(\U,\V)}^{p}\right] \right)^{\frac1p},
\]
where we identify $(\HH\otimes\U)\otimes\V$ with $\HH\otimes\mathcal{L}_2(\U,\V)$ and
\[
\|DZ\|_{\HH\otimes\mathcal{L}_2(\U,\V)} =\left(\alpha_H\int_0^T\int_0^T \langle D_uZ,D_vZ\rangle_{\mathcal{L}_2(\U,\V)} |u-v|^{2H-2}dudv\right)^{\frac12},
\]
whenever $DZ$ has a function representative for which this integral is absolutely convergent. For $Z=Z(t)$ with $D_uZ(t)=0$ for $u>t$, the integration domain reduces to $[0,t]^2$. Here $D$ is understood as the closure of the Malliavin derivative initially defined on smooth cylindrical random variables (cf. \cite[Section 1.2]{N06}).

For a random initial value independent of $W^H$, we use the product-space extension of $\mathbb D^{1,p}(\V)$, with Malliavin differentiation acting only on the Gaussian noise component. Thus $X_0\in L^p(\Omega;\V)$ implies $X_0\in\mathbb D^{1,p}(\V)$ and $DX_0=0$.

\begin{lemma}[Malliavin chain rule for Lipschitz G\^ateaux differentiable maps]
\label{lm:chain}
Let $G_0,G_1$ be real separable Hilbert spaces, $1<p<\infty$, and $\mathfrak H:=\HH\otimes\U$. Suppose that $\Psi:G_0\to G_1$ is globally Lipschitz with constant $L$ and is G\^ateaux differentiable everywhere, with derivative $\Psi'(x)\in\mathcal L(G_0,G_1)$. If $Z\in\mathbb D^{1,p}(G_0)$, then $\Psi(Z)\in\mathbb D^{1,p}(G_1)$ and
\[
D(\Psi(Z))=\Psi'(Z)DZ,\qquad \|D(\Psi(Z))\|_{\mathfrak H\otimes G_1} \le L\|DZ\|_{\mathfrak H\otimes G_0}\quad\text{a.s.},
\]
where $\Psi'(Z)$ acts on the $G_0$-valued part of $DZ$.
\end{lemma}

\begin{proof}
The Lipschitz composition result in \cite[Section 3.4]{PV14} gives $\Psi(Z)\in\mathbb D^{1,p}(G_1)$. We identify its derivative. On the canonical Gaussian space, let $\tau_{\varepsilon e}$ denote the Cameron--Martin shift in direction $e\in\mathfrak H$, and write $\Delta_\varepsilon^eY:=(Y\circ\tau_{\varepsilon e}-Y)/\varepsilon$. For any real separable Hilbert space $G$ and $Y\in\mathbb D^{1,p}(G)$,
\[
\Delta_\varepsilon^eY\to (DY)e \quad\text{in }L^q(\Omega;G),\qquad 1\le q<p.
\]
Indeed, \cite[Theorem 4.1]{MPR17} gives this for each finite orthogonal projection of $Y$. For smooth cylindrical $Y$, the fundamental theorem of calculus and the Cameron--Martin formula (\cite[Proposition 3.1]{MPR17}) give
\[
\sup_{0<|\varepsilon|\le1} \|\Delta_\varepsilon^eY\|_{L^q(\Omega;G)} \le C_{e,p,q}\|DY\|_{L^p(\Omega;\mathfrak H\otimes G)}.
\]
This bound extends by density to $\mathbb D^{1,p}(G)$ and controls the projection tails uniformly in $\varepsilon$, proving the Hilbert-valued assertion.

Set $v_\varepsilon:=\Delta_\varepsilon^eZ$ and $v:=(DZ)e$. Then
\begin{align*}
\|\Delta_\varepsilon^e(\Psi(Z))-\Psi'(Z)v\|_{G_1} \le L\|v_\varepsilon-v\|_{G_0}+\left\| \frac{\Psi(Z+\varepsilon v)-\Psi(Z)}{\varepsilon}-\Psi'(Z)v \right\|_{G_1}.
\end{align*}
The first term tends to zero in probability, and the second tends to zero almost surely by G\^ateaux differentiability. The same difference-quotient result applied to $\Psi(Z)$ identifies its limit as $(D(\Psi(Z)))e$. Uniqueness of limits in probability and a countable dense set of directions $e$ yield the claimed identity. Since $\|\Psi'(x)\|_{\mathcal L(G_0,G_1)}\le L$, the norm bound follows. Directional difference quotients and finite-rank expansions ensure measurability of $\Psi'(Z)DZ$.
\end{proof}

\begin{proposition}\label{prop:Malliavin-X}
Let Assumptions \ref{as:A}--\ref{as:F} hold. Fix $p\in(1,\infty)$, and suppose that $X_0\in L^p(\Omega;V)$ is independent of $W^H$. Then $X(t)\in\mathbb D^{1,p}(V)$ for every $t\in[0,T]$. Its Malliavin derivative admits, for almost every $u\in[0,t)$, that
\begin{equation}
\label{eq:DX} D_uX(t) =E(t-u)\Phi +\int_u^t E(t-r)DF(X(r))D_uX(r) dr
\end{equation}
in $\mathcal{L}_2(\U,\V)$, while $D_uX(t)=0$ for almost every $u>t$. Moreover,
\begin{align}
\label{eq:DXpoint} &\|D_uX(t)\|_{L^p(\Omega;\mathcal{L}_2(\U,\V))} \le C(t-u)^{\frac{\beta-1}2}, \qquad 0\le u<t\le T,\\\label{eq:DXsup} &\sup_{t\in[0,T]} \|D_{\cdot}X(t)\|_{L^p(\Omega;\mathcal H\otimes \mathcal{L}_2(\U,\V))} \le C,
\end{align}
where the constants $C$ are independent of $u$ and $t$.
\end{proposition}
\begin{proof}
Since $X_0\in L^p(\Omega;V)$ is independent of $W^H$, thus $X_0\in\mathbb D^{1,p}(V)$ and $DX_0=0$. Throughout the proof, the generic constant $C$ may depend on $T$, $H$, $\beta$, the bound for $DF$, and $\|(-A)^{\frac{\beta-1}{2}}\Phi\|_{\mathcal{L}_2(\U,\V)}$, but is independent of $n$, $m$, $u$, $t$, and $X_0$.

Consider the Picard iterates
\[
X^{(0)}(t)=E(t)X_0,\qquad X^{(n+1)}(t) =E(t)X_0+\int_0^tE(t-r)F(X^{(n)}(r)) dr+O_t,\qquad n\in\mathbb N.
\]
Lemma \ref{lm:O} and \eqref{eq:F} imply that $X^{(n)}\to X$ in $C([0,T];L^p(\Omega;V))$ and $O_t\in\mathbb D^{1,p}(\V)$ with
\[
D_uO_t=\mathbf 1_{[0,t)}(u)E(t-u)\Phi.
\]
Since $F$ is globally Lipschitz and everywhere G\^ateaux differentiable, we may apply Lemma \ref{lm:chain} with $\Psi=F$ at each Picard step. It follows by induction that $X^{(n)}(t)\in\mathbb D^{1,p}(V)$ and, for almost every $u<t$,
\begin{align}
\label{eq:DXn} D_uX^{(n+1)}(t) =\int_u^tE(t-r)DF(X^{(n)}(r))D_uX^{(n)}(r) dr+E(t-u)\Phi,
\end{align}
whereas $D_uX^{(n)}(t)=0$ for almost every $u>t$. At each Picard step, differentiation under the time integral is justified by the derivative estimates below, and the closedness of the Malliavin derivative.

For $u<t$, the semigroup smoothing estimate \eqref{eq:semi1} gives
\[
\|E(t-u)\Phi\|_{\mathcal{L}_2(\U,\V)} \le C(t-u)^{\frac{\beta-1}{2}} \|(-A)^{\frac{\beta-1}{2}}\Phi\|_{\mathcal{L}_2(\U,\V)} \le C(t-u)^{\frac{\beta-1}{2}},
\]
which, together with the boundedness of $E(t)$ and $DF$, leads to
\[
\|D_uX^{(n+1)}(t)\|_{\mathcal{L}_2(\U,\V)} \le C\int_u^t\|D_uX^{(n)}(r)\|_{\mathcal{L}_2(\U,\V)} dr+C(t-u)^{\frac{\beta-1}{2}}.
\]
Since $\frac{\beta-1}{2}>-H>-1$, the integral form of Gronwall's inequality, applied to $\max\limits_{0\le n\le N}\|D_uX^{(n)}(t)\|_{\mathcal{L}_2(\U,\V)}$ and followed by $N\to\infty$, yields
\begin{align}
\label{eq:DXn-norm} \sup_{n\ge0}\|D_uX^{(n)}(t)\|_{\mathcal{L}_2(\U,\V)} \le C\left((t-u)^{\frac{\beta-1}{2}} +\int_u^t e^{C(t-r)}(r-u)^{\frac{\beta-1}{2}} dr\right) \le C(t-u)^{\frac{\beta-1}{2}}.
\end{align}
Consequently, by the integral representation of the norm in $\HH\otimes\mathcal{L}_2(\U,\V)$ and the Cauchy--Schwarz inequality, we derive
\begin{align}
\label{eq:DXn-sup} \sup_{n\ge0} \|D_{\cdot}X^{(n)}(t)\|_{\mathcal H\otimes \mathcal{L}_2(\U,\V)}^2 &\le C\int_0^t\int_0^t (t-u)^{\frac{\beta-1}{2}}(t-v)^{\frac{\beta-1}{2}} |u-v|^{2H-2} du dv\notag\\ &=Ct^{\rho} \int_0^1\int_0^1 x^{\frac{\beta-1}2}y^{\frac{\beta-1}2}|x-y|^{2H-2} dx dy \le Ct^\rho,
\end{align}
where in the last step we used
\begin{align}
\label{eq:double-int} \int_0^1\int_0^1 x^{\frac{\beta-1}2}y^{\frac{\beta-1}2}|x-y|^{2H-2} dx dy &=2\int_0^1\int_0^y x^{\frac{\beta-1}2}y^{\frac{\beta-1}2}(y-x)^{2H-2} dx dy\notag\\ &=2\left(\int_0^1y^{2H+\beta-2} dy\right)\left(\int_0^1 \tilde x^{\frac{\beta-1}2}(1-\tilde x)^{2H-2} d\tilde x\right)\notag\\ &=\frac{2}{\rho} B\left(\frac{\beta+1}2,2H-1\right)
\end{align}
with $B(\cdot,\cdot)$ denoting the beta function, under conditions $\rho=2H+\beta-1>0$, $\frac{\beta+1}2>0$, and $2H-1>0$, or equivalently, $H>\frac12$ and $\beta>1-2H$.

We next justify convergence of the derivatives. The expression \eqref{eq:DXn}, together with the boundedness of $E(t)$ and $DF$, gives
\[
\begin{aligned}
&\|D_uX^{(n+1)}(t)-D_uX^{(m+1)}(t)\|_ {L^2(\Omega;\mathcal{L}_2(\U,\V))}\\ \le&~ C\int_u^t \|D_uX^{(n)}(r)-D_uX^{(m)}(r)\|_ {L^2(\Omega;\mathcal{L}_2(\U,\V))} dr\\ &+C\int_u^t(r-u)^{\frac{\beta-1}{2}} \|DF(X^{(n)}(r))-DF(X^{(m)}(r))\|_ {L^2(\Omega;\mathcal{L}(V))} dr,
\end{aligned}
\]
Using $\|(-A)^{\rho/2}E(t)\|_{\mathcal L(\V)}\le Ct^{-\rho/2}$ and the global Lipschitz bound, the mild equations give
\[
\|X^{(n)}-X\|_{C([0,T];L^p(\Omega;\V^\rho))} \le C\sup_{t\in[0,T]}\int_0^t(t-r)^{-\rho/2} \|F(X^{(n-1)}(r))-F(X(r))\|_{L^p(\Omega;\V)}\,dr \to0.
\]
By the Lipschitz estimate for $DF$ on $\V^\rho$ in Assumption \ref{as:F} and the uniform boundedness of $DF$ in \eqref{eq:F}, this implies
\[
\lim_{n,m\to\infty} \|DF(X^{(n)}(r))-DF(X^{(m)}(r))\|_{L^2(\Omega;\mathcal{L}(V))}=0.
\]
By the dominated convergence theorem and the estimate \eqref{eq:DXn-norm},
\[
\begin{aligned}
&\limsup_{n,m\to\infty} \|D_uX^{(n)}(t)-D_uX^{(m)}(t)\|_ {L^2(\Omega;\mathcal{L}_2(\U,\V))}\\ \le&~ C\int_u^t \limsup_{n,m\to\infty} \|D_uX^{(n)}(r)-D_uX^{(m)}(r)\|_ {L^2(\Omega;\mathcal{L}_2(\U,\V))} dr,
\end{aligned}
\]
which, according to Gronwall's inequality, implies for almost every $u<t$ that
\[
\lim_{n,m\to\infty} \|D_uX^{(n)}(t)-D_uX^{(m)}(t)\|_ {L^2(\Omega;\mathcal{L}_2(\U,\V))}=0.
\]
By the covariance formula, Cauchy--Schwarz inequality, and the dominated convergence theorem,
\[
\begin{aligned}
&\mathbb E\left[\|D_{\cdot}X^{(n)}(t)-D_{\cdot}X^{(m)}(t)\|_ {\mathcal H\otimes \mathcal{L}_2(\U,\V)}^2\right]\\ \le&~\alpha_H\int_0^t\int_0^t \|D_uX^{(n)}(t)-D_uX^{(m)}(t)\|_ {L^2(\Omega;\mathcal{L}_2(\U,\V))}\\ &\qquad\times \|D_vX^{(n)}(t)-D_vX^{(m)}(t)\|_ {L^2(\Omega;\mathcal{L}_2(\U,\V))} |u-v|^{2H-2} du dv \to0.
\end{aligned}
\]
Thus $\{D_{\cdot}X^{(n)}(t)\}_{n\in\mathbb N}$ is a Cauchy sequence in $L^p(\Omega;\mathcal H\otimes \mathcal{L}_2(\U,\V))$ for every finite $p\ge1$. Indeed, for $p\le2$, this follows from monotonicity of probability-space norms. For $p>2$, it follows from
\begin{align*}
\mathbb E\left[\|D_{\cdot}X^{(n)}(t)-D_{\cdot}X^{(m)}(t)\|_ {\mathcal H\otimes \mathcal{L}_2(\U,\V)}^p\right] \le(Ct^{\frac{\rho}2})^{p-2}\mathbb E\left[\|D_{\cdot}X^{(n)}(t)-D_{\cdot}X^{(m)}(t)\|_ {\mathcal H\otimes \mathcal{L}_2(\U,\V)}^2\right]\to0
\end{align*}
according to \eqref{eq:DXn-sup}.

Since $X^{(n)}(t)\to X(t)$ in $L^p(\Omega;V)$, the closedness of the Malliavin derivative gives
\[
X(t)\in\mathbb D^{1,p}(V),\qquad D_{\cdot}X^{(n)}(t)\to D_{\cdot}X(t) \quad\text{in } L^p(\Omega;\mathcal H\otimes \mathcal{L}_2(\U,\V)).
\]

For fixed $t$, the limits of $D_uX^{(n)}(t)$ in $L^2(\Omega;\mathcal{L}_2(\U,\V))$, defined for almost every $u<t$ and extended by zero for $u\ge t$, represent $D_{\cdot}X(t)$. Indeed, this identification follows from dominated convergence in the integral representation of the $\HH\otimes\mathcal{L}_2(\U,\V)$-norm. Passing to the limit in
\[
D_uX^{(n+1)}(t) =E(t-u)\Phi+ \int_u^t E(t-r)DF(X^{(n)}(r))D_uX^{(n)}(r) dr
\]
then yields \eqref{eq:DX} for almost every $u<t$, while the chosen representative vanishes for $u\ge t$.

The preceding convergence statements imply convergence in probability. Hence the uniform almost-sure bounds in \eqref{eq:DXn-norm} and \eqref{eq:DXn-sup} also hold for their limits
\[
\|D_uX(t)\|_{\mathcal{L}_2(\U,\V)} \le C(t-u)^{\frac{\beta-1}{2}},\qquad \|D_{\cdot}X(t)\|_{\mathcal H\otimes \mathcal{L}_2(\U,\V)} \le Ct^{\frac{\rho}2},
\]
which finishes the proof by taking $p$th moments.
\end{proof}

Estimates for the Malliavin derivative of the numerical solution can be obtained similarly to the exact solution case.

\begin{proposition}\label{prop:Malliavin-Xm}
Let Assumptions \ref{as:A}--\ref{as:F} hold. Fix $p\in(1,\infty)$, and suppose that $X_0\in L^p(\Omega;\V)$ is independent of $W^H$. Then $X^m(\kappa_m(s))\in\mathbb D^{1,p}(\V)$ for every $m\in\mathbb N^+$ and $s\in[0,T]$. Its Malliavin derivative admits, for almost every $u\in[0,\kappa_m(s))$, that
\begin{equation}
\label{eq:DXm} D_uX^m(\kappa_m(s)) ={}E(\kappa_m(s)-u)\Phi +\int_0^{\kappa_m(s)} E(\kappa_m(s)-r)DF(X^m(\kappa_m(r))) D_uX^m(\kappa_m(r))dr
\end{equation}
in $\mathcal{L}_2(\U,\V)$, while $D_uX^m(\kappa_m(s))=0$ for $u>\kappa_m(s)$. Moreover, for almost every $u\in[0,\kappa_m(s))$,
\begin{equation}
\label{eq:DXm-norm} \|D_uX^m(\kappa_m(s))\|_ {L^p(\Omega;\mathcal{L}_2(\U,\V))} \le C\left[ (\kappa_m(s)-u)^{\frac{\beta-1}{2}} +h\mathbf 1_{\{\overline\kappa_m(u)<\kappa_m(s)\}} (\overline\kappa_m(u)-u)^{\frac{\beta-1}{2}} \right],
\end{equation}
where $\overline\kappa_m(u)=t_{\ell}$ for $u\in(t_{\ell-1},t_{\ell})$, and
\begin{equation}
\label{eq:DXm-sup} \sup_{m\ge1}\max_{0\le k\le m} \|D_{\cdot}X^m(t_k)\|_ {L^p(\Omega;\HH\otimes\mathcal{L}_2(\U,\V))} \le C.
\end{equation}
Here, $C$ is independent of $m$, $s$, $u$, and $X_0$.
\end{proposition}

\begin{proof}
At each grid point, \eqref{eq:Xm} reads
\[
X^m(t_k)=E(t_k)X_0+ \sum_{i=0}^{k-1}\int_{t_i}^{t_{i+1}} E(t_k-r)F(X^m(t_i)) dr+O_{t_k}.
\]
As in Proposition \ref{prop:Malliavin-X}, $O_{t_k}\in\mathbb D^{1,p}(\V)$ and $D_uO_{t_k}=\mathbf 1_{[0,t_k)}(u)E(t_k-u)\Phi$. Applying Lemma \ref{lm:chain} with $\Psi=F$ at each time step gives by induction $X^m(t_k)\in\mathbb D^{1,p}(\V)$ and
\[
D_uX^m(t_k) =E(t_k-u)\Phi+ \sum_{i=0}^{k-1}\int_{t_i}^{t_{i+1}} E(t_k-r)DF(X^m(t_i))D_uX^m(t_i) dr, \qquad u<t_k,
\]
with the representative chosen to vanish for $u\ge t_k$. This proves \eqref{eq:DXm}.

Fix $u\in(t_{\ell-1},t_\ell)$. For $k\ge\ell$, the smoothing estimate \eqref{eq:semi1}, the boundedness of $DF$ in \eqref{eq:F}, and the vanishing of $D_uX^m(t_i)$ for $i<\ell$ imply
\[
\|D_uX^m(t_k)\|_{\mathcal{L}_2(\U,\V)} \le C(t_k-u)^{\frac{\beta-1}{2}} +Ch\sum_{i=\ell}^{k-1} \|D_uX^m(t_i)\|_{\mathcal{L}_2(\U,\V)}.
\]
Iterating this inequality and using $(1+Ch)^m\le e^{CT}$ yields
\[
\|D_uX^m(t_k)\|_{\mathcal{L}_2(\U,\V)} \le C(t_k-u)^{\frac{\beta-1}{2}} +Ch\sum_{i=\ell}^{k-1}(t_i-u)^{\frac{\beta-1}{2}}.
\]
For $k=\ell$, the sum is empty. For $k>\ell$, since $-1<\frac{\beta-1}{2}\le0$,
\[
h\sum_{i=\ell+1}^{k-1}(t_i-u)^{\frac{\beta-1}{2}} \le\int_{t_\ell}^{t_{k-1}}(r-u)^{\frac{\beta-1}{2}}dr \le C(t_k-u)^{\frac{\beta+1}{2}} \le CT(t_k-u)^{\frac{\beta-1}{2}}.
\]
Separating the term $i=\ell$ and taking $p$th moments proves \eqref{eq:DXm-norm}.

For \eqref{eq:DXm-sup}, similar to the proof of Proposition \ref{prop:Malliavin-X}, we estimate the derivative in $\HH\otimes\mathcal{L}_2(\U,\V)$. By Lemma \ref{lm:O} and the boundedness of $E(t)$ and $DF$ in $\mathcal{L}(\V)$, the derivative recursion gives
\[
\|D_{\cdot}X^m(t_k)\|_{\HH\otimes\mathcal{L}_2(\U,\V)} \le Ct_k^{\frac{\rho}2} +Ch\sum_{i=0}^{k-1} \|D_{\cdot}X^m(t_i)\|_{\HH\otimes\mathcal{L}_2(\U,\V)}.
\]
Since $DX_0=0$, the discrete Gronwall inequality yields
\[
\|D_{\cdot}X^m(t_k)\|_{\HH\otimes\mathcal{L}_2(\U,\V)} \le Ct_k^{\frac{\rho}2}
\]
uniformly in $m$ and $k$. Taking $p$th moments proves \eqref{eq:DXm-sup}.
\end{proof}

\subsection{Strong convergence rate}
We first establish two integral estimates for the strong convergence analysis, which will be used repeatedly in the following sections and are therefore stated separately. Lemma \ref{lm:int} controls weighted integrals over a single time step, while Lemma \ref{lm:Gj} estimates the covariance correction arising from the divergence decomposition.

\begin{lemma}\label{lm:int}
Let $a\in(-1,0]$ and $b>-1$. Then there exists a constant $C>0$, depending only on $a$ and $b$, such that
\[
\int_r^{(\kappa_m(r)+h)\wedge t} (t-s)^a(s-r)^bds \le Ch^{b+1}(t-r)^a, \qquad 0\le r<t\le T,
\]
where $s\wedge t:=\min\{s,t\}$.
\end{lemma}

\begin{proof}
If $(\kappa_m(r)+h)\wedge t-r\le \frac{t-r}2$, then
\[
\int_r^{(\kappa_m(r)+h)\wedge t} (t-s)^a(s-r)^bds \le \left(\frac{t-r}2\right)^a \int_r^{(\kappa_m(r)+h)\wedge t}(s-r)^b ds \le \frac1{2^a(b+1)}h^{b+1}(t-r)^a.
\]
If $(\kappa_m(r)+h)\wedge t-r> \frac{t-r}2$, then $t-r<2h$. By setting $x=\frac{s-r}{t-r}$, we obtain
\[
\begin{aligned}
\int_r^{(\kappa_m(r)+h)\wedge t} (t-s)^a(s-r)^b ds &\le \int_r^t(t-s)^a(s-r)^b ds \\ &= (t-r)^{a+b+1} \int_0^1(1-x)^a x^b dx \\ &\le Ch^{b+1}(t-r)^a.
\end{aligned}
\]
Here in the last step we used the boundedness of the beta function $\int_0^1(1-x)^{a}x^{b}dx$ since $a,b>-1$.
\end{proof}

\begin{lemma}\label{lm:Gj}
Let $\beta\in(1-2H,1]$ and $\sigma\in[0,1-\beta]$. For $j\in\mathbb N$, denote
\[
\mathcal{G}_j(s):=\int_0^{t_j}\int_{t_j}^{s} (t_j-u)^{\frac{\beta-1}2}(s-r)^{\frac{\beta-1}{2}+\frac{\sigma}{2}}(r-u)^{2H-2}dr du,\qquad s\in(t_j,t_{j+1}).
\]
Then there exists a constant $C>0$ independent of $j$ and $s$, such that
\begin{equation*}
\max_{j=0,1,\cdots,k-1}\sup_{s\in(t_j,t_{j+1})}\mathcal{G}_j(s) \le Ch^{\min\left\{\rho+\frac{\sigma}{2},\frac{\beta+\sigma+1}{2}\right\}}\left(1+{\bf 1}_{\{\beta=3-4H\}}|\log h|\right),
\end{equation*}
where $\mathbf{1}_{\{\beta = 3-4H\}}$ equals $1$ if $\beta = 3-4H$ and $0$ otherwise.
\end{lemma}

\begin{proof}
By letting $x:=\frac{t_j-u}{s-t_j}$ and $y:=\frac{r-t_j}{s-t_j}$, we obtain
\begin{align*}
\mathcal{G}_j(s)&=\int_0^{\frac{t_j}{s-t_j}}\int_0^1((s-t_j)x)^{\frac{\beta-1}2}(s-t_j-(s-t_j)y)^{\frac{\beta-1}{2}+\frac{\sigma}{2}}((s-t_j)(x+y))^{2H-2}(s-t_j)^2dydx\\ &=(s-t_j)^{\rho+\frac{\sigma}{2}}\int_0^{\frac{t_j}{s-t_j}}\int_0^1x^{\frac{\beta-1}2}(1-y)^{\frac{\beta-1}{2}+\frac{\sigma}{2}}(x+y)^{2H-2}dydx.
\end{align*}

Since $\mathcal{G}_0(s)\equiv0$, we only consider the case $j\ge1$, which yields $\frac{t_j}{s-t_j}>\frac{t_j}h\ge1$ for $s\in(t_j,t_{j+1})$. We then split the last integral above for $x\in(0,1)$ and $x\in(1,\frac{t_j}{s-t_j})$.

For $x\in(0,1)$, it holds
\begin{align*}
&\int_0^1\int_0^1x^{\frac{\beta-1}2}(1-y)^{\frac{\beta-1}{2}+\frac{\sigma}{2}}(x+y)^{2H-2}dydx\\ \le&~\int_0^1\int_0^{\frac12}x^{\frac{\beta-1}2}\left(\frac12\right)^{\frac{\beta-1}{2}+\frac{\sigma}{2}}(x+y)^{2H-2}dydx+\int_0^1\int_{\frac12}^1x^{\frac{\beta-1}2}(1-y)^{\frac{\beta-1}{2}+\frac{\sigma}{2}}\left(\frac12\right)^{2H-2}dydx\\ \le&~ C\int_0^1x^{\frac{\beta-1}2}\left[\left(x+\frac12\right)^{2H-1}-x^{2H-1}\right]dx+C\int_0^1x^{\frac{\beta-1}2}dx\int_{\frac12}^1(1-y)^{\frac{\beta-1}{2}+\frac{\sigma}{2}}dy\le C.
\end{align*}

For $x\in(1,\frac{t_j}{s-t_j})$, the fact $x+y>x$ leads to
\begin{align*}
\int_1^{\frac{t_j}{s-t_j}}\int_0^1x^{\frac{\beta-1}2}(1-y)^{\frac{\beta-1}{2}+\frac{\sigma}{2}}(x+y)^{2H-2}dydx &\le\int_1^{\frac{t_j}{s-t_j}}x^{\frac{\beta-(3-4H)}2-1}dx\int_0^1(1-y)^{\frac{\beta-1}{2}+\frac{\sigma}{2}}dy\\ &\le \left\{
\begin{aligned}
&C\left(\frac{t_j}{s-t_j}\right)^{\max\{0,\frac{\beta-(3-4H)}2\}},\qquad&&\beta\neq3-4H,\\ &C\log \frac{t_j}{s-t_j},\qquad&&\beta=3-4H,
\end{aligned}
\right.
\end{align*}
where the constant $C$ is independent of $j$ and $s$. The above estimate is sharp by noting also that $x+y<2x$ for $x>1$ and $y\in(0,1)$.

We then derive for $\beta\neq3-4H$ that
\begin{align*}
\mathcal{G}_j(s)\le C\left((s-t_j)^{\rho+\frac{\sigma}{2}}+(s-t_j)^{\min\left\{\rho+\frac{\sigma}{2},\frac{\beta+\sigma+1}2\right\}}\right)\le Ch^{\min\left\{\rho+\frac{\sigma}{2},\frac{\beta+\sigma+1}2\right\}}
\end{align*}
and for $\beta=3-4H$, equivalently $\rho=2H+\beta-1=\frac{\beta+1}2$, that
\[
\mathcal{G}_j(s)\le (s-t_j)^{\rho+\frac{\sigma}{2}}\left(C+C|\log(s-t_j)|\right)\le Ch^{\rho+\frac{\sigma}{2}}(1+|\log h|),
\]
where in the last step we used the fact $z^{\rho+\frac{\sigma}2}(1+|\log z|)\le Ch^{\rho+\frac{\sigma}{2}}(1+|\log h|)$ for all $z\in(0,h)$. It then completes the proof.
\end{proof}

Combining these integral estimates with the regularity and Malliavin derivative bounds established above, we now derive the strong convergence rate at the grid points.

\begin{theorem}\label{tm:rate-k}
Let Assumptions \ref{as:A}--\ref{as:F} hold. Fix $p\in[1,\infty)$. Suppose that $X_0\in L^{2p}(\Omega;\V^\rho)$ is independent of $W^H$, and let $X(t_k)$ and $X^m(t_k)$ be defined by \eqref{eq:mild} and \eqref{eq:Xm}, respectively. Denote
\begin{equation*}
q:=\min\left\{1,\rho+\frac{\gamma}2,\rho+\frac{\sigma}2,\frac{\beta+\sigma+1}2\right\}.
\end{equation*}
Then, there exists a constant $C > 0$, independent of $m$, $h$, and $X_0$, such that
\begin{equation}
\label{eq:rate-k} \sup_{0\le k\le m} \|X^m(t_k) - X(t_k)\|_{L^p(\Omega;\V)} \le C\bigl(1 + \|X_0\|_{L^{2p}(\Omega;\V^\rho)}^2\bigr) h^q \Bigl(1 + \mathbf{1}_{\{\beta = 3-4H\}} |\log h|\Bigr).
\end{equation}
\end{theorem}

\begin{proof}
Subtracting \eqref{eq:mild} from \eqref{eq:Xm}, and using the boundedness of $E$ and $DF$, yields
\begin{align}
\label{eq:error-split} &\quad\|X^m(t_k)-X(t_k)\|_{L^p(\Omega;\V)}\notag\\ &=\left\|\sum_{j=0}^{k-1}\int_{t_j}^{t_{j+1}}E(t_k-s) \bigl(F(X^m(t_j))-F(X(s))\bigr) ds\right\|_{L^p(\Omega;\V)}\notag\\ &\le\left\|\sum_{j=0}^{k-1}\int_{t_j}^{t_{j+1}}E(t_k-s) \bigl(F(X(s))-F(X(t_j))\bigr) ds\right\|_{L^p(\Omega;\V)} +Ch\sum_{j=0}^{k-1}\|X^m(t_j)-X(t_j)\|_{L^p(\Omega;\V)}\notag\\ &=:I_k+Ch\sum_{j=0}^{k-1}\|X^m(t_j)-X(t_j)\|_{L^p(\Omega;\V)}.
\end{align}

By Assumption \ref{as:F} and the mean-value inequality, the map $r\mapsto DF(u+rz)z$ is Lipschitz continuous from $[0,1]$ to $\V^{-\eta}$ for every $u,z\in\V$. The fundamental theorem of calculus therefore yields Taylor's formula in $\V^{-\eta}$. For $s\in[t_j,t_{j+1})$, we obtain
\begin{align*}
F(X(s))-F(X(t_j)) &=DF(X(t_j))(E(s-t_j)-I)X(t_j)+DF(X(t_j))\int_{t_j}^sE(s-r)F(X(r)) dr\\ &\quad+DF(X(t_j))\int_{t_j}^sE(s-r)\Phi dW^H(r)\\ &\quad+\int_0^1(1-\lambda)D^2F \bigl(X(t_j)+\lambda(X(s)-X(t_j))\bigr)\bigl(X(s)-X(t_j),X(s)-X(t_j)\bigr) d\lambda,
\end{align*}
which leads to
\begin{align}
I_k&\le\left\|\sum_{j=0}^{k-1}\int_{t_j}^{t_{j+1}} E(t_k-s)DF(X(t_j))(E(s-t_j)-I)X(t_j) ds \right\|_{L^p(\Omega;\V)}\notag\\ &\qquad+\left\|\sum_{j=0}^{k-1}\int_{t_j}^{t_{j+1}} E(t_k-s)DF(X(t_j))\int_{t_j}^sE(s-r)F(X(r)) dr ds \right\|_{L^p(\Omega;\V)}\notag\\ &\qquad+\left\|\sum_{j=0}^{k-1}\int_{t_j}^{t_{j+1}} E(t_k-s)DF(X(t_j))\int_{t_j}^sE(s-r)\Phi dW^H(r) ds \right\|_{L^p(\Omega;\V)}\notag\\ &\qquad+\Biggl\|\sum_{j=0}^{k-1}\int_{t_j}^{t_{j+1}} E(t_k-s)\int_0^1(1-\lambda)D^2F \bigl(X(t_j)+\lambda(X(s)-X(t_j))\bigr)\notag\\[-1mm] &\hspace{31mm}\bigl(X(s)-X(t_j),X(s)-X(t_j)\bigr) d\lambda ds \Biggr\|_{L^p(\Omega;\V)}\notag\\ &=: I_{1,k}+I_{2,k}+I_{3,k}+I_{4,k}. \label{eq:I-decomp}
\end{align}

{\it Step 1. Estimates of $I_{1,k}, I_{2,k},$ and $I_{4,k}$.} For $I_{1,k}$, using \eqref{eq:semi1}, the assumption \eqref{eq:F1}, the fact
\begin{align*}
\|(E(t)-I)v\|_{-(\rho+\gamma)} = \|(-A)^{-\left(\rho+\frac{\gamma}{2}\right)}(E(t)-I)(-A)^{\frac{\rho}2}v\| \le Ct^{\min\left\{1,\rho+\frac{\gamma}{2}\right\}}\|v\|_\rho,
\end{align*}
and Theorem \ref{tm:X}, we obtain
\begin{align}
\label{eq:I1k} I_{1,k} &\le C\sum_{j=0}^{k-1}\int_{t_j}^{t_{j+1}} (t_k-s)^{-\frac{\delta}2} \|DF(X(t_j))(E(s-t_j)-I)X(t_j)\|_{L^p(\Omega;\V^{-\delta})} ds \notag\\ &\le C\sum_{j=0}^{k-1}\int_{t_j}^{t_{j+1}} (t_k-s)^{-\frac{\delta}2} (s-t_j)^{\min\left\{1,\rho+\frac{\gamma}{2}\right\}}ds\bigl(1+\|X(t_j)\|_{L^{2p}(\Omega;\V^\rho)}^2\bigr)\notag\\ &\le C\left(1+\|X_0\|_{L^{2p}(\Omega;\V^\rho)}^2\right) h^{\min\left\{1,\rho+\frac{\gamma}{2}\right\}},
\end{align}
where in the last step we used the fact that the time integral is uniformly bounded since $\delta\le\rho<2$.

For $I_{2,k}$, the boundedness of $E$ and $DF$ in \eqref{eq:F}, and Theorem \ref{tm:X} yield
\begin{align}
\label{eq:I2k} I_{2,k} \le C\sum_{j=0}^{k-1}\int_{t_j}^{t_{j+1}}\int_{t_j}^s \bigl(1+\|X(r)\|_{L^p(\Omega;\V)}\bigr) dr ds \le C(1+\|X_0\|_{L^p(\Omega;\V)})h.
\end{align}

For the Taylor remainder $I_{4,k}$, the property \eqref{eq:semi1}, the assumption \eqref{eq:F2}, and Theorem \ref{tm:X} give
\begin{align}
\label{eq:I4k} I_{4,k} &\le C\sum_{j=0}^{k-1}\int_{t_j}^{t_{j+1}} (t_k-s)^{-\frac{\eta}2} \|X(s)-X(t_j)\|_{L^{2p}(\Omega;\V)} \|X(s)-X(t_j)\|_{L^{2p}(\Omega;\V^{-\sigma})} ds \notag\\ &\le C\bigl(1+\|X_0\|_{L^{2p}(\Omega;\V^\rho)}^2\bigr) \sum_{j=0}^{k-1}\int_{t_j}^{t_{j+1}} (t_k-s)^{-\frac{\eta}2}(s-t_j)^{\frac{\rho}{2}}(s-t_j)^{\frac{\rho+\sigma}{2}} ds \notag\\ &\le C\bigl(1+\|X_0\|_{L^{2p}(\Omega;\V^\rho)}^2\bigr) \sum_{j=0}^{k-1}\int_{t_j}^{t_{j+1}} (t_k-s)^{-\frac{\eta}2}(s-t_j)^{\rho+\frac{\sigma}{2}} ds \notag\\ &\le C\bigl(1+\|X_0\|_{L^{2p}(\Omega;\V^\rho)}^2\bigr) h^{\rho+\frac{\sigma}{2}},
\end{align}
where the last time integral is finite because $\eta<2$.

{\it Step 2. Estimate of $I_{3,k}$.} Set $p_*:=\max\{2,p\}\le2p$ and denote
\[
B_j(s):=E(t_k-s)DF(X(t_j)),\qquad f_{j,s}(r):=\1_{[t_j,s]}(r)E(s-r)\Phi.
\]
For fixed $t_j<s<t_{j+1}$, set
\[
\mathcal K:=\HH\otimes\HS(\U,\V),\qquad \Psi_{j,s}(x):=E(t_k-s)DF(x)f_{j,s},\quad x\in\V.
\]
By \eqref{eq:semi1}, \eqref{eq:F}, and \eqref{eq:double-int}, $\|\Psi_{j,s}(x)\|_{\mathcal K}\le C(s-t_j)^{\rho/2}$, so this map is well-defined. The semigroup estimates, the directional differentiability in Assumption \ref{as:F}, \eqref{eq:F2}, the mean-value inequality, and \eqref{eq:double-int} with $\beta+\sigma$ in place of $\beta$ give
\begin{align*}
\|(-A)^{-\sigma/2}E(s-r)\Phi\|_{\HS(\U,\V)} &\le C(s-r)^{(\beta+\sigma-1)/2},\qquad t_j\le r<s,\\* \|\Psi_{j,s}(x)-\Psi_{j,s}(y)\|_{\mathcal K} &\le C(t_k-s)^{-\eta/2}(s-t_j)^{(\rho+\sigma)/2}\|x-y\|.
\end{align*}
Thus $\Psi_{j,s}$ is globally Lipschitz. We next verify that $\Psi_{j,s}$ has a directional derivative in the norm of $\mathcal K$. Fix $x,v\in\V$. For $\theta\in\mathbb R\setminus\{0\}$, define the difference quotient
\[
\Delta_\theta :=\frac{\Psi_{j,s}(x+\theta v)-\Psi_{j,s}(x)}{\theta}\in\mathcal K.
\]
Define the candidate derivative kernel $L_{x,v}$ by
\[
L_{x,v}(r)g_n :=E(t_k-s)D^2F(x)\bigl(v,f_{j,s}(r)g_n\bigr), \qquad t_j\le r<s,\quad n\ge1,
\]
and set $L_{x,v}(r)=0$ for $r\notin[t_j,s)$. For every fixed $r\in[t_j,s)$ and $n\ge1$, the directional differentiability in Assumption \ref{as:F} and the boundedness of $E(t_k-s):\V^{-\eta}\to\V$ imply $\Delta_\theta(r)g_n\to L_{x,v}(r)g_n$ in $\V$ as $\theta\to0$. Moreover, \eqref{eq:F2} and the mean-value inequality give
\[
\|\Delta_\theta(r)g_n\|+\|L_{x,v}(r)g_n\| \le C(t_k-s)^{-\eta/2}\|v\|\, \|(-A)^{-\sigma/2}E(s-r)\Phi g_n\|,
\]
for all $\theta\ne0$, $t_j\le r<s$, and $n\ge1$. For each fixed $r<s$, the sequence of squared coordinate bounds is summable because
\[
\sum_{n\ge1}\|(-A)^{-\sigma/2}E(s-r)\Phi g_n\|^2 =\|(-A)^{-\sigma/2}E(s-r)\Phi\|_{\HS(\U,\V)}^2<\infty.
\]
The coordinate estimate and the summability show that $L_{x,v}(r)\in\HS(\U,\V)$ for each fixed $r\in[t_j,s)$. Dominated convergence for the series $\sum_{n\ge1}\|\Delta_\theta(r)g_n-L_{x,v}(r)g_n\|^2$ therefore proves convergence in the Hilbert--Schmidt norm:
\[
d_\theta(r):=\|\Delta_\theta(r)-L_{x,v}(r)\|_{\HS(\U,\V)} \to0\qquad(\theta\to0).
\]
Write $\alpha=(\beta+\sigma-1)/2$. The coordinate bounds and the preceding semigroup estimate also give
\[
d_\theta(r)+\|L_{x,v}(r)\|_{\HS(\U,\V)} \le C(t_k-s)^{-\eta/2}\|v\|(s-r)^\alpha, \qquad t_j\le r<s,
\]
uniformly in $\theta\ne0$. Equation \eqref{eq:double-int}, with $\beta+\sigma$ in place of $\beta$, yields
\[
\int_{t_j}^s\int_{t_j}^s (s-r)^\alpha(s-r')^\alpha|r-r'|^{2H-2}\,dr'\,dr =C(s-t_j)^{\rho+\sigma}<\infty,
\]
since $\rho+\sigma>0$. The bound for $\|L_{x,v}(r)\|_{\HS(\U,\V)}$ and the covariance norm formula \eqref{eq:HH} show that $L_{x,v}\in\mathcal K$. Applying dominated convergence to $d_\theta(r)d_\theta(r')|r-r'|^{2H-2}$ on $[t_j,s]^2$ gives
\[
\|\Delta_\theta-L_{x,v}\|_{\mathcal K}^2 \le\alpha_H\int_{t_j}^s\int_{t_j}^s d_\theta(r)d_\theta(r')|r-r'|^{2H-2}\,dr'\,dr \to0\qquad(\theta\to0).
\]
Thus the $\mathcal K$-valued directional derivative of $\Psi_{j,s}$ at $x$ in direction $v$ is $L_{x,v}$. The bilinearity of $D^2F(x)$ makes $v\mapsto L_{x,v}$ linear, and the same covariance estimate gives
\[
\|L_{x,v}\|_{\mathcal K} \le C(t_k-s)^{-\eta/2}(s-t_j)^{(\rho+\sigma)/2}\|v\|.
\]
Consequently, $\Psi_{j,s}$ is G\^ateaux differentiable with bounded linear derivative $\Psi'_{j,s}(x)v=L_{x,v}$. Since $X(t_j)\in\mathbb D^{1,p_*}(\V)$ by Proposition \ref{prop:Malliavin-X}, Lemma \ref{lm:chain}, with $G_0=\V$ and $G_1=\mathcal K$, yields $B_j(s)f_{j,s}\in\mathbb D^{1,p_*}(\mathcal K)$ and $D(B_j(s)f_{j,s})=\Psi'_{j,s}(X(t_j))DX(t_j)$. In coordinates, for $t_j\le r<s<t_{j+1}$, this is
\[
D_u^i\bigl(B_j(s)f_{j,s}\bigr)(r)g_n = \1_{[t_j,s]}(r)E(t_k-s)D^2F(X(t_j)) \bigl(D_u^iX(t_j),E(s-r)\Phi g_n\bigr).
\]
By \eqref{eq:F}, \eqref{eq:F2}, \eqref{eq:DXpoint}, and \eqref{eq:double-int},
\begin{align*}
\|B_j(s)f_{j,s}\|_ {L^{p_*}(\Omega;\HH\otimes\HS(\U,\V))} &\le C(s-t_j)^{\frac{\rho}{2}},\\ \|D(B_j(s)f_{j,s})\|_ {L^{p_*}(\Omega;(\HH\otimes\U)^{\otimes2}\otimes\V)} &\le C(t_k-s)^{-\frac{\eta}{2}} t_j^{\frac{\rho}{2}} (s-t_j)^{\frac{\rho+\sigma}{2}}.
\end{align*}
Since $\eta<2$, it follows that
\begin{align*}
\sum_{j=0}^{k-1}\int_{t_j}^{t_{j+1}} \|B_j(s)f_{j,s}\|_ {\mathbb D^{1,p_*}(\HH\otimes\HS(\U,\V))}ds &\le C\sum_{j=0}^{k-1}\int_{t_j}^{t_{j+1}} \left( (s-t_j)^{\frac{\rho}{2}} +(t_k-s)^{-\frac{\eta}{2}} t_j^{\frac{\rho}{2}} (s-t_j)^{\frac{\rho+\sigma}{2}} \right)ds\\ &\le CTh^{\frac{\rho}{2}} +CT^{\frac{\rho}{2}}h^{\frac{\rho+\sigma}{2}} \int_0^{t_k}(t_k-s)^{-\frac{\eta}{2}}ds <\infty.
\end{align*}
Hence, $ \mathcal U_k:=\sum_{j=0}^{k-1}\int_{t_j}^{t_{j+1}} B_j(s)f_{j,s}ds $ is well-defined in $\mathbb D^{1,p_*}(\HH\otimes\HS(\U,\V))$, and $D$ can be interchanged with the time integral. The divergence estimate (cf. \cite[Proposition 6.10]{M10})
\begin{equation}
\label{eq:div-Lp} \|\deltaH(v)\|_{L^r(\Omega;\V)} \le C\left( \|v\|_{L^r(\Omega;\HH\otimes\HS(\U,\V))} +\|Dv\|_{L^r(\Omega;(\HH\otimes\U)^{\otimes2}\otimes\V)} \right), \qquad 1<r<\infty,
\end{equation}
for $v\in\mathbb D^{1,r}(\HH\otimes\HS(\U,\V))$ also allows $\deltaH$ to be interchanged with the time integral.

By the product formula for the divergence operator (cf. \cite[Proposition 1.3.3]{N06}), we obtain
\begin{align}
\label{eq:I3split} \mathcal I_{3,k} &:=\sum_{j=0}^{k-1}\int_{t_j}^{t_{j+1}} E(t_k-s)DF(X(t_j)) \int_{t_j}^sE(s-r)\Phi dW^H(r)ds\notag\\ &=\sum_{j=0}^{k-1}\int_{t_j}^{t_{j+1}} B_j(s)\deltaH(f_{j,s})ds\notag\\ &=\deltaH(\mathcal U_k)+\mathcal C_k,
\end{align}
where
\begin{align*}
\mathcal C_k &:=\sum_{j=0}^{k-1}\int_{t_j}^{t_{j+1}} \operatorname{Tr}_{\HH\otimes\U} \langle DB_j(s),f_{j,s}\rangle ds\\ &=\alpha_H\sum_{j=0}^{k-1} \int_{t_j}^{t_{j+1}}\int_0^{t_j}\int_{t_j}^s \sum_{n\ge1}E(t_k-s)D^2F(X(t_j)) \bigl(D_u^nX(t_j),E(s-r)\Phi g_n\bigr) |r-u|^{2H-2}drduds.
\end{align*}

{\it Step 2.1. Estimate of $\deltaH(\mathcal U_k)$.} To estimate the $p$th moment of $\mathcal I_{3,k}$, we use \eqref{eq:div-Lp}. Since $L^{p_*}(\Omega;\V)\hookrightarrow L^p(\Omega;\V)$, it suffices to estimate the two terms on the right-hand side of \eqref{eq:div-Lp} with $r=p_*$.

For $\mathcal U_k$, we first consider the case $\beta<1$. Then $\delta\le\rho<2H$, and \eqref{eq:HH} yields
\begin{align*}
\|\mathcal U_k\|_ {L^{p_*}(\Omega;\HH\otimes\HS(\U,\V))}^2 &=\alpha_H\left\| \int_0^T\int_0^T \langle\mathcal U_k(r_1),\mathcal U_k(r_2)\rangle_{\HS(\U,\V)} |r_1-r_2|^{2H-2}dr_1dr_2 \right\|_{L^{p_*/2}(\Omega)}\\ &\le\alpha_H\sum_{j,l=0}^{k-1} \int_{t_l}^{t_{l+1}}\int_{t_j}^{t_{j+1}} \int_{t_l}^{s_2}\int_{t_j}^{s_1} \bigl\|E(t_k-s_1)DF(X(t_j))E(s_1-r_1)\Phi \bigr\|_{L^{p_*}(\Omega;\HS(\U,\V))}\\ &\qquad\times \bigl\|E(t_k-s_2)DF(X(t_l))E(s_2-r_2)\Phi \bigr\|_{L^{p_*}(\Omega;\HS(\U,\V))} |r_1-r_2|^{2H-2}dr_1dr_2ds_1ds_2.
\end{align*}
By \eqref{eq:semi1}, \eqref{eq:F1}, and Theorem \ref{tm:X},
\begin{align*}
&\bigl\|E(t_k-s_1)DF(X(t_j))E(s_1-r_1)\Phi \bigr\|_{L^{p_*}(\Omega;\HS(\U,\V))}\\ \le&~ C(t_k-s_1)^{-\frac{\delta}{2}} \bigl(1+\|X(t_j)\|_{L^{p_*}(\Omega;\V^\rho)}\bigr) \bigl\| (-A)^{\frac{1-\beta-\rho-\gamma}{2}}E(s_1-r_1) \bigr\|_{\cL(\V)} \bigl\|(-A)^{\frac{\beta-1}{2}}\Phi\bigr\|_{\HS(\U,\V)}\\ \le&~ C\bigl(1+\|X_0\|_{L^{p_*}(\Omega;\V^\rho)}\bigr) (t_k-s_1)^{-\frac{\delta}{2}} (s_1-r_1)^{\min\left\{0,\frac{\beta+\rho+\gamma-1}{2}\right\}}.
\end{align*}
Changing the order of integration and applying Lemma \ref{lm:int} with $t=t_k$ give
\begin{align}
\label{eq:Ubound} \|\mathcal U_k\|_ {L^{p_*}(\Omega;\HH\otimes\HS(\U,\V))}^2 &\le C \bigl(1+\|X_0\|_{L^{p_*}(\Omega;\V^\rho)}\bigr)^2 \sum_{j,l=0}^{k-1} \int_{t_l}^{t_{l+1}}\int_{t_j}^{t_{j+1}}\notag\\ &\qquad\times \left( \int_{r_1}^{t_{j+1}} (t_k-s_1)^{-\frac{\delta}{2}} (s_1-r_1)^{\min\left\{0,\frac{\beta+\rho+\gamma-1}{2}\right\}}ds_1 \right)\notag\\ &\qquad\times \left( \int_{r_2}^{t_{l+1}} (t_k-s_2)^{-\frac{\delta}{2}} (s_2-r_2)^{\min\left\{0,\frac{\beta+\rho+\gamma-1}{2}\right\}}ds_2 \right) |r_1-r_2|^{2H-2}dr_1dr_2\notag\\ &\le C \bigl(1+\|X_0\|_{L^{p_*}(\Omega;\V^\rho)}\bigr)^2 h^{2\min\left\{0,\frac{\beta+\rho+\gamma-1}{2}\right\}+2} \notag\\ &\qquad\times \int_0^{t_k}\int_0^{t_k} (t_k-r_1)^{-\frac{\delta}{2}} (t_k-r_2)^{-\frac{\delta}{2}} |r_1-r_2|^{2H-2}dr_1dr_2\notag\\ &\le C \bigl(1+\|X_0\|_{L^{p_*}(\Omega;\V^\rho)}\bigr)^2 h^{\min\left\{2,\beta+\rho+\gamma+1\right\}},
\end{align}
where the last integral is bounded since $\delta<2H$ similar to the estimate in \eqref{eq:double-int}.

When $\beta=1$, Assumption \ref{as:Phi} gives $\Phi\in\HS(\U,\V)$. For almost every $r\in[0,t_k)$,
\[
\mathcal U_k(r) = \int_r^{\kappa_m(r)+h} E(t_k-s)DF(X(\kappa_m(r)))E(s-r)\Phi ds.
\]
The boundedness of $E$ and $DF$ implies
\[
\|\mathcal U_k(r)\|_{\HS(\U,\V)} \le Ch\|\Phi\|_{\HS(\U,\V)}.
\]
Consequently,
\begin{align*}
\|\mathcal U_k\|_ {L^{p_*}(\Omega;\HH\otimes\HS(\U,\V))}^2 &\le C\alpha_Hh^2\|\Phi\|_{\HS(\U,\V)}^2 \int_0^{t_k}\int_0^{t_k} |r_1-r_2|^{2H-2}dr_1dr_2\\ &=Ch^2\|\Phi\|_{\HS(\U,\V)}^2t_k^{2H}\le Ch^2.
\end{align*}
Thus \eqref{eq:Ubound} also holds for $\beta=1$, including the case $\delta=\rho=2H$.

To estimate $D\mathcal U_k$, we use the index $i$ for the Malliavin direction and the index $n$ for the cylindrical-noise direction of the integrand. Applying the chain rule, we obtain
\begin{align*}
D_u^i\mathcal U_k(r)g_n = \sum_{j=0}^{k-1}\int_{t_j}^{t_{j+1}} \1_{[t_j,s]}(r)E(t_k-s)D^2F(X(t_j)) \bigl(D_u^iX(t_j),E(s-r)\Phi g_n\bigr)ds.
\end{align*}
Using \eqref{eq:HH}, \eqref{eq:F2}, and Minkowski's inequality, we obtain
\begin{align*}
&\qquad \|D\mathcal U_k\|_ {L^{p_*}(\Omega;(\HH\otimes\U)^{\otimes2}\otimes\V)}^2\\ &=\alpha_H^2\left\| \int_{[0,T]^4}\sum_{i,n\ge1} \left\langle D_{u_1}^i\mathcal U_k(r_1)g_n, D_{u_2}^i\mathcal U_k(r_2)g_n \right\rangle |u_1-u_2|^{2H-2}|r_1-r_2|^{2H-2} du_1du_2dr_1dr_2 \right\|_{L^{p_*/2}(\Omega)}\\ &\le C\sum_{j,l=0}^{k-1} \int_{t_l}^{t_{l+1}}\int_{t_j}^{t_{j+1}} (t_k-s_1)^{-\frac{\eta}{2}} (t_k-s_2)^{-\frac{\eta}{2}} \int_{t_l}^{s_2}\int_{t_j}^{s_1} \int_0^{t_l}\int_0^{t_j} \left\| \sum_{i\ge1} \|D_{u_1}^iX(t_j)\| \|D_{u_2}^iX(t_l)\| \right\|_{L^{p_*/2}(\Omega)}\\ &\qquad\quad\times \sum_{n\ge1} \|E(s_1-r_1)\Phi g_n\|_{-\sigma} \|E(s_2-r_2)\Phi g_n\|_{-\sigma} |u_1-u_2|^{2H-2}|r_1-r_2|^{2H-2} du_1du_2dr_1dr_2ds_1ds_2.
\end{align*}
By Cauchy--Schwarz inequality and $\sigma\in[0,1-\beta]$,
\begin{align*}
&\sum_{n\ge1} \|E(s_1-r_1)\Phi g_n\|_{-\sigma} \|E(s_2-r_2)\Phi g_n\|_{-\sigma}\\ \le&~ \|(-A)^{-\frac{\sigma}{2}}E(s_1-r_1)\Phi\|_{\HS(\U,\V)} \|(-A)^{-\frac{\sigma}{2}}E(s_2-r_2)\Phi\|_{\HS(\U,\V)}\\ \le&~ C\|(-A)^{\frac{\beta-1}{2}}\Phi\|_{\HS(\U,\V)}^2 (s_1-r_1)^{\frac{\beta-1+\sigma}{2}} (s_2-r_2)^{\frac{\beta-1+\sigma}{2}}.
\end{align*}
Moreover, H\"older's inequality and \eqref{eq:DXpoint} yield
\begin{align*}
&\int_0^{t_l}\int_0^{t_j} \left\| \sum_{i\ge1} \|D_{u_1}^iX(t_j)\| \|D_{u_2}^iX(t_l)\| \right\|_{L^{p_*/2}(\Omega)} |u_1-u_2|^{2H-2}du_1du_2\\ \le&~ \int_0^{t_l}\int_0^{t_j} \|D_{u_1}X(t_j)\|_{L^{p_*}(\Omega;\HS(\U,\V))} \|D_{u_2}X(t_l)\|_{L^{p_*}(\Omega;\HS(\U,\V))} |u_1-u_2|^{2H-2}du_1du_2\\ \le&~ C\int_0^{t_l}\int_0^{t_j} (t_j-u_1)^{\frac{\beta-1}{2}} (t_l-u_2)^{\frac{\beta-1}{2}} |u_1-u_2|^{2H-2}du_1du_2 \le C,
\end{align*}
since $\beta>1-2H$ as in the proof of \eqref{eq:Ubound}. Applying Lemma \ref{lm:int} with $t=t_k$, we derive
\begin{align}
\label{eq:DUbound} \|D\mathcal U_k\|_ {L^{p_*}(\Omega;(\HH\otimes\U)^{\otimes2}\otimes\V)}^2 &\le C\sum_{j,l=0}^{k-1} \int_{t_l}^{t_{l+1}}\int_{t_j}^{t_{j+1}} (t_k-s_1)^{-\frac{\eta}{2}} (t_k-s_2)^{-\frac{\eta}{2}}\notag\\ &\qquad\times \int_{t_l}^{s_2}\int_{t_j}^{s_1} (s_1-r_1)^{\frac{\beta-1+\sigma}{2}} (s_2-r_2)^{\frac{\beta-1+\sigma}{2}} |r_1-r_2|^{2H-2}dr_1dr_2ds_1ds_2 \notag\\ &=C\sum_{j,l=0}^{k-1} \int_{t_l}^{t_{l+1}}\int_{t_j}^{t_{j+1}} \left( \int_{r_1}^{t_{j+1}} (t_k-s_1)^{-\frac{\eta}{2}} (s_1-r_1)^{\frac{\beta-1+\sigma}{2}}ds_1 \right)\notag\\ &\qquad\times \left( \int_{r_2}^{t_{l+1}} (t_k-s_2)^{-\frac{\eta}{2}} (s_2-r_2)^{\frac{\beta-1+\sigma}{2}}ds_2 \right) |r_1-r_2|^{2H-2}dr_1dr_2 \notag\\ &\le Ch^{\beta+\sigma+1} \int_0^{t_k}\int_0^{t_k} (t_k-r_1)^{-\frac{\eta}{2}} (t_k-r_2)^{-\frac{\eta}{2}} |r_1-r_2|^{2H-2}dr_1dr_2 \notag\\ &\le Ch^{\beta+\sigma+1},
\end{align}
since $\eta<2H$ as in the proof of \eqref{eq:Ubound}. Combining \eqref{eq:div-Lp}--\eqref{eq:DUbound}, we obtain
\begin{equation}
\label{eq:deltaU} \|\deltaH(\mathcal U_k)\|_{L^p(\Omega;\V)} \le C \bigl(1+\|X_0\|_{L^{p_*}(\Omega;\V^\rho)}\bigr) h^{\min\left\{ 1,\frac{\beta+\rho+\gamma+1}{2}, \frac{\beta+\sigma+1}{2} \right\}}.
\end{equation}

{\it Step 2.2. Estimate of $\mathcal C_k$.} Based on the property \eqref{eq:semi1}, the assumption \eqref{eq:F2}, the estimate for the Malliavin derivative given in \eqref{eq:DXpoint}, and Lemma \ref{lm:Gj}, we derive
\begin{align*}
\|\mathcal C_k\|_{L^{p_*}(\Omega;V)} &\le C\sum_{j=0}^{k-1}\int_{t_j}^{t_{j+1}}\int_0^{t_j}\int_{t_j}^{s} (t_k-s)^{-\frac{\eta}2}\|D_uX(t_j)\|_{L^{p_*}(\Omega;\mathcal{L}_2(\U,\V))}\\ &\qquad\qquad\times\|(-A)^{-\frac{\sigma}2}E(s-r)\Phi\|_{\mathcal{L}_2(\U,\V)}|r-u|^{2H-2} dr du ds\\ &\le C\sum_{j=0}^{k-1}\int_{t_j}^{t_{j+1}}(t_k-s)^{-\frac{\eta}2} \left(\int_0^{t_j}\int_{t_j}^{s} (t_j-u)^{\frac{\beta-1}2}(s-r)^{\frac{\beta-1+\sigma}{2}}(r-u)^{2H-2}dr du\right)ds\\ &\le Ch^{\min\left\{\rho+\frac{\sigma}{2},\frac{\beta+\sigma+1}{2}\right\}}\left(1+{\bf 1}_{\{\beta=3-4H\}}|\log h|\right)\sum_{j=0}^{k-1}\int_{t_j}^{t_{j+1}}(t_k-s)^{-\frac{\eta}2}ds\\ &\le Ch^{\min\left\{\rho+\frac{\sigma}{2},\frac{\beta+\sigma+1}{2}\right\}}\left(1+{\bf 1}_{\{\beta=3-4H\}}|\log h|\right)
\end{align*}
since $\eta<2H<2$. By $L^{p_*}(\Omega;\V)\hookrightarrow L^p(\Omega;\V)$ and $p_*\le 2p$, the above estimate, together with \eqref{eq:I3split} and \eqref{eq:deltaU}, yields
\begin{equation}
\label{eq:I3kfinal} I_{3,k}=\|\mathcal I_{3,k}\|_{L^p(\Omega;V)} \le C\bigl(1+\|X_0\|_{L^{2p}(\Omega;\V^\rho)}\bigr)h^{\min\left\{1,\frac{\beta+\rho+\gamma+1}2,\rho+\frac{\sigma}{2},\frac{\beta+\sigma+1}{2}\right\}} \left(1+{\bf 1}_{\{\beta=3-4H\}}|\log h|\right).
\end{equation}

Combining \eqref{eq:I-decomp}--\eqref{eq:I4k} and \eqref{eq:I3kfinal} with \eqref{eq:error-split} and utilizing the fact $\rho+\frac{\gamma}2<\frac{\beta+\rho+\gamma+1}2$, we derive
\begin{align*}
\|X^m(t_k)-X(t_k)\|_{L^p(\Omega;\V)}\le&~ C\bigl(1+\|X_0\|_{L^{2p}(\Omega;\V^\rho)}^2\bigr)h^{\min\left\{1,\rho+\frac{\gamma}2,\rho+\frac{\sigma}2,\frac{\beta+\sigma+1}2\right\}}\left(1+{\bf 1}_{\{\beta=3-4H\}}|\log h|\right)\notag\\ &+Ch\sum_{j=0}^{k-1}\|X^m(t_j)-X(t_j)\|_{L^p(\Omega;\V)},
\end{align*}
which completes the proof of \eqref{eq:rate-k} based on the discrete Gronwall lemma.
\end{proof}

Theorem \ref{tm:rate-k} establishes the strong convergence rate at the grid points. We next extend this result to the continuous interpolation $X^m$ defined in \eqref{eq:Xm} and show that the same error bound holds uniformly for all $t\in[0,T]$.

\begin{theorem}\label{tm:rate}
Under the assumptions of Theorem \ref{tm:rate-k}, there exists a constant $C > 0$, independent of $m$, $h$, and $X_0$, such that
\begin{equation*}
\sup_{t\in[0,T]}\| X^m(t)-X(t)\|_{L^p(\Omega;\V)} \le C_p\bigl(1 + \|X_0\|_{L^{2p}(\Omega;\V^\rho)}^2\bigr) h^q \Bigl(1 + \mathbf{1}_{\{\beta = 3-4H\}} |\log h|\Bigr),
\end{equation*}
where $q$ is defined in Theorem \ref{tm:rate-k}.
\end{theorem}

\begin{proof}
Fix $t\in[t_k,t_{k+1}]$. Then \eqref{eq:mild}, \eqref{eq:Xm}, and the semigroup property give
\[
X^m(t)-X(t)=E(t-t_k)\left(X^m(t_k)-X(t_k)\right) +\int_{t_k}^tE(t-s) \bigl(F(X^m(t_k))-F(X(s))\bigr) ds.
\]
By \eqref{eq:F} and Theorem \ref{tm:X},
\begin{align*}
\|X^m(t)-X(t)\|_{L^p(\Omega;\V)} &\le \|X^m(t_k)-X(t_k)\|_{L^p(\Omega;\V)}\\ &\quad+C\int_{t_k}^t\bigl(\|X^m(t_k)-X(t_k)\|_{L^p(\Omega;\V)}+\|X(t_k)-X(s)\|_{L^p(\Omega;\V)}\bigr)ds\\ &\le C_p(1+h)\|X^m(t_k)-X(t_k)\|_{L^p(\Omega;\V)} +C\bigl(1+\|X_0\|_{L^p(\Omega;\V^\rho)}\bigr)h^{1+\frac{\rho}2}.
\end{align*}
The assertion follows from Theorem \ref{tm:rate-k} since $q\le1$.
\end{proof}

\subsection{Uniform moment bounds and temporal regularity of the numerical solution} Having established the strong convergence rate, we next derive uniform moment bounds and temporal regularity estimates for the numerical solution $X^m$. These estimates will be used to control the remainder terms and spatial approximation errors in the asymptotic error analysis of the next section.

\begin{lemma}\label{lm:Xm-bound}
Let Assumptions \ref{as:A}, \ref{as:Phi}, and \ref{as:F}~\eqref{eq:F} hold. Fix $p\in[1,\infty)$ and $\tilde\delta\in[0,\rho]$, and suppose that $X_0\in L^p(\Omega;\V^{\tilde\delta})$ is independent of $W^H$. Then
\[
\sup_{m\ge1}\sup_{t\in[0,T]} \|X^m(t)\|_{L^p(\Omega;\V^{\tilde\delta})} \le C\bigl(1+\|X_0\|_{L^p(\Omega;\V^{\tilde\delta})}\bigr),
\]
where $C$ is independent of $m$ and $t$.
\end{lemma}

\begin{proof}
By \eqref{eq:Xm}, \eqref{eq:F}, and Lemma \ref{lm:O}, the nondecreasing function $M_m(t):=\sup_{r\in[0,t]}\|X^m(r)\|_{L^p(\Omega;\V)}$ satisfies
\[
M_m(t)\le C\bigl(1+\|X_0\|_{L^p(\Omega;\V)}\bigr) +C\int_0^tM_m(s) ds.
\]
Gronwall's inequality therefore gives
\begin{equation}
\label{supXm} \sup_{m\ge1}\sup_{t\in[0,T]}\|X^m(t)\|_{L^p(\Omega;\V)} \le C\bigl(1+\|X_0\|_{L^p(\Omega;\V)}\bigr).
\end{equation}
Applying \eqref{eq:semi1} and \eqref{eq:O} to \eqref{eq:Xm}, and then using \eqref{supXm}, yields
\begin{align*}
\|X^m(t)\|_{L^p(\Omega;\V^{\tilde\delta})} &\le \|X_0\|_{L^p(\Omega;\V^{\tilde\delta})} +C\int_0^t(t-s)^{-\frac{\tilde\delta}2} \bigl(1+\|X^m(\kappa_m(s))\|_{L^p(\Omega;\V)}\bigr) ds+C\\ &\le C\bigl(1+\|X_0\|_{L^p(\Omega;\V^{\tilde\delta})}\bigr),
\end{align*}
since $\tilde\delta\le\rho<2$.
\end{proof}

The moment bound yields the same increment estimates as for the exact solution, with constants uniform in the time discretization.

\begin{lemma}\label{lm:Xm-holder}
Under the assumptions of Lemma \ref{lm:Xm-bound}, suppose further that $X_0\in L^p(\Omega;\V^\rho)$. For $\tilde\delta\in[0,\rho]$, $\tilde\sigma\in[0,1-\beta]$, and $0\le s<t\le T$, we have
\begin{align*}
\sup_{m\ge1}\|X^m(t)-X^m(s)\|_{L^p(\Omega;\V^{\tilde\delta})} &\le C\bigl(1+\|X_0\|_{L^p(\Omega;\V^\rho)}\bigr) (t-s)^{\frac{\rho-\tilde\delta}2},\\ \sup_{m\ge1}\|X^m(t)-X^m(s)\|_{L^p(\Omega;\V^{-\tilde\sigma})} &\le C\bigl(1+\|X_0\|_{L^p(\Omega;\V^\rho)}\bigr) (t-s)^{\frac{\rho+\tilde\sigma}2}.
\end{align*}
\end{lemma}

\begin{proof}
The identity
\[
X^m(t)-X^m(s)=(E(t-s)-I)X^m(s) +\int_s^tE(t-r)F(X^m(\kappa_m(r))) dr+O_{s,t}
\]
has the same form as the increment identity in the proof of Theorem \ref{tm:X}. Its proof applies with Lemma \ref{lm:Xm-bound} in place of the moment bound of the exact solution, using \eqref{eq:semi1}--\eqref{eq:semi2} and \eqref{eq:Ost}--\eqref{eq:neg-Ost}. For the negative norm, the needed exponent satisfies $0<\frac{\rho+\tilde\sigma}2\le H<1$.
\end{proof}

\section{Asymptotic error distribution}
\label{sec:asymptotic}
Building on the strong error estimates in Section \ref{sec:3}, we next investigate the asymptotic distribution of the numerical error in a particular parameter regime.

Throughout this section, let Assumptions \ref{as:A}--\ref{as:F} hold with fixed $\beta,\sigma$, and $\gamma$ satisfying in addition
\begin{equation}
\label{eq:new-index} \beta\in(3-4H,1],\qquad \sigma=1-\beta,\qquad \rho+\frac{\gamma}{2}\ge1.
\end{equation}

The following examples illustrate Assumptions \ref{as:A}--\ref{as:F} and \eqref{eq:new-index} in two different noise regimes.
\begin{example}[$\beta=1$]\label{ex:local-drift}
Let $\mathcal O=(0,1)$, $\V=\U=L^2(\mathcal O)$, and $A$ be the Dirichlet Laplacian with domain $H^2(\mathcal O)\cap H_0^1(\mathcal O)$, which satisfies Assumption \ref{as:A} apparently. 

For $1\le s\le2$, we have $\V^s=H^s(\mathcal O)\cap H_0^1(\mathcal O)$, with equivalent norms. Choose
\[
H\in\Bigl(\frac12,\frac34\Bigr),\qquad \beta=1,\qquad \rho=\delta=2H,\qquad \gamma=\sigma=0,\qquad \eta=1,
\]
and take $\Phi=(-A)^{-1/2}$. Then $\|(-A)^{\frac{\beta-1}{2}}\Phi\|_{\HS(\U,\V)}^2 =\sum_{k=1}^{\infty}\lambda_k^{-1}<\infty$, which verifies Assumption \ref{as:Phi}. 

Define the nonlinear Nemytskii operator
\[
F(u)(x):=f(u(x)),\qquad f(z):=\arctan z, \qquad g(z):=f'(z)=\frac{1}{1+z^2}.
\]
Since $|f'|\le1$ and $f(0)=0$, $F$ is globally Lipschitz, $\|F(u)\|\le\|u\|$, and
\[
DF(u)v=g(u)v,\qquad \|DF(u)\|_{\cL(\V)}\le1.
\]
The derivative exists in the G\^ateaux sense by dominated convergence, so \eqref{eq:F} holds.

To verify \eqref{eq:F1} with a constant independent of $u$, we use the boundedness of the absolute-value map on $H^s$ for $1<s<3/2$ (cf. \cite[Section~1]{MN19}). A bounded Sobolev extension from $\mathcal O$ to $\mathbb R$ gives
\[
\bigl\||z|\bigr\|_{H^s(\mathcal O)} \le C_s\|z\|_{H^s(\mathcal O)}, \qquad z\in H^s(\mathcal O),\quad 1<s<\frac32.
\]
For the above $g$, integration by parts yields
\[
g(z)=\frac12\int_{\mathbb R}|z-a|g''(a)\,da, \qquad g''(a)=\frac{6a^2-2}{(1+a^2)^3}, \qquad \int_{\mathbb R}(1+|a|)|g''(a)|\,da<\infty.
\]
Consequently, since $1<\rho<3/2$, for every $u\in\V^\rho$,
\begin{align*}
\|g(u)\|_{H^\rho(\mathcal O)} &\le\frac12\int_{\mathbb R}|g''(a)| \bigl\||u-a|\bigr\|_{H^\rho(\mathcal O)}\,da\\* &\le C\int_{\mathbb R}|g''(a)| \bigl(\|u\|_{H^\rho(\mathcal O)}+|a|\bigr)\,da \le C(1+\|u\|_\rho).
\end{align*}
The one-dimensional Sobolev multiplication estimate, and the zero boundary trace of $\varphi\in\V^\rho$, now imply
\[
\|g(u)\varphi\|_\rho \le C\|g(u)\|_{H^\rho(\mathcal O)}\|\varphi\|_\rho \le C(1+\|u\|_\rho)\|\varphi\|_\rho.
\]
Defining multiplication on $\V^{-\rho}$ by duality therefore gives
\[
\|DF(u)v\|_{-\rho} =\sup_{\|\varphi\|_\rho=1}|\langle v,g(u)\varphi\rangle| \le C(1+\|u\|_\rho)\|v\|_{-\rho}, \qquad v\in\V^{-\rho},
\]
which proves \eqref{eq:F1}.

Moreover, $L^1(\mathcal O)\hookrightarrow\V^{-1}=(H_0^1(\mathcal O))^*$ and boundedness of $g'$ give
\[
D^2F(u)(v,w)=g'(u)vw,\qquad \|D^2F(u)(v,w)\|_{-1} \le C\|g'(u)vw\|_{L^1(\mathcal O)} \le C\|v\|\,\|w\|,
\]
which is \eqref{eq:F2} for $\eta=1$ and $\sigma=0$. For $u,v,w\in\V$, dominated convergence, with dominating function $\|g'\|_\infty|vw|\in L^1(\mathcal O)$, shows that
\[
\frac{DF(u+\varepsilon v)w-DF(u)w}{\varepsilon} \to g'(u)vw \quad\text{in }L^1(\mathcal O)\hookrightarrow\V^{-1}.
\]
This verifies the directional differentiability required in Assumption \ref{as:F}.  If $u_n\to u$ in $\V$, continuity and boundedness of $g$ give $g(u_n)v\to g(u)v$ in $\V$ for each fixed $v\in\V$, proving the strong operator continuity of $DF$ on $\V$. Finally, the embedding $\V^\rho\hookrightarrow L^\infty(\mathcal O)$ yields
\[
\|DF(u)-DF(v)\|_{\cL(\V)} \le\|g'\|_\infty\|u-v\|_{L^\infty(\mathcal O)} \le C\|u-v\|_\rho, \qquad u,v\in\V^\rho,
\]
which verifies the Lipschitz estimate for $DF$ on $\V^\rho$. Hence Assumption \ref{as:F} is satisfied.

The chosen parameters satisfy
\[
\beta=1>3-4H,\qquad \sigma=1-\beta=0, \qquad \rho+\frac\gamma2=2H>1,
\]
and therefore fulfill \eqref{eq:new-index}.
\end{example}

\begin{example}[$\beta<1$]\label{ex:nonlocal-drift}
Retain $V=U=L^2(0,1)$ and the Dirichlet Laplacian $A$ as in Example \ref{ex:local-drift}. Set $\Phi=I$ and choose
\[
H=\frac34,\qquad \beta=\frac14,\qquad \rho=\sigma=\frac34,\qquad \gamma=\frac12,\qquad \delta=0,\qquad \eta=1.
\]
Since $\lambda_k=\pi^2k^2$,
\[
\|(-A)^{(\beta-1)/2}\Phi\|_{\mathcal L_2(U,V)}^2 =\sum_{k\ge1}\lambda_k^{-3/4}<\infty.
\]
Let $e_1$ be the normalized first eigenfunction and define
\[
F(u)=f(\langle u,e_1\rangle)e_1,\qquad f(z)=\arctan z.
\]
Then
\begin{align*}
&DF(u)v=f'(\langle u,e_1\rangle)\langle v,e_1\rangle e_1,\\
&D^2F(u)(v,w) =f''(\langle u,e_1\rangle) \langle v,e_1\rangle\langle w,e_1\rangle e_1.
\end{align*}
Because $e_1\in V^r$ for every $r\ge0$, duality and the boundedness of $f'$ and $f''$ give
\[
\|DF(u)v\|\le C\|v\|_{-5/4}, \qquad \|D^2F(u)(v,w)\|_{-1} \le C\|v\|\|w\|_{-3/4}.
\]
Moreover,
\[
\|DF(u)-DF(v)\|_{\mathcal L(V)} \le C\|u-v\|\le C\|u-v\|_{3/4}, \qquad u,v\in V^{3/4}.
\]
Thus Assumption \ref{as:F} holds, with the extensions following \eqref{eq:F1}--\eqref{eq:F2} defined by the same formulas. Finally,
\[
\beta>3-4H=0,\qquad \sigma=1-\beta,\qquad \rho+\frac{\gamma}{2}=1,
\]
so \eqref{eq:new-index} is satisfied.
\end{example}

Fix $p\in[2,\infty)$ and assume that $X_0\in L^{2p}(\Omega;\V^\rho)$ is independent of $W^H$. Then Theorem \ref{tm:rate} applies with $q=1$ and without a logarithmic factor. More precisely, since $h=\frac{T}{m}$, it gives
\[
\sup_{t\in[0,T]} \|X^m(t)-X(t)\|_{L^p(\Omega;\V)}\le Cm^{-1}.
\]
Therefore, we define the normalized error process
\[
U^m(t):=m\bigl(X^m(t)-X(t)\bigr),\qquad t\in[0,T].
\]
To identify its limiting distribution, we first derive a linear integral equation for $U^m$ up to vanishing remainder terms. We then establish a spectral approximation uniform in $m$, identify the limit for each fixed spatial dimension, and pass to the infinite-dimensional limit in $C([0,T];\V)$.

Unless otherwise stated, constants may depend on the fixed data, $p$, and $\|X_0\|_{L^{2p}(\Omega;\V^\rho)}$, but are independent of $m$, $n$, and the time variables.

\subsection{Decomposition of the normalized error and an auxiliary process $\tilde U^m$}

Subtracting \eqref{eq:mild} from \eqref{eq:Xm} and expanding $F$ at $X(s)$ and $X^m(\kappa_m(s))$ gives
\begin{align}
\label{eq:Um} U^m(t)=&~ m \int_0^t E(t-s)\big(F(X^m(s))-F(X(s))\big)ds \notag\\ &\qquad - m \int_0^t E(t-s)\big(F(X^m(s))-F(X^m(\kappa_m(s)))\big)ds\notag\\ =&\int_0^tE(t-s)DF(X(s))U^m(s) ds+R^m_1(t)\notag\\ &-B^m_1(t)-B^m_2(t)-B^m_3(t)-R^m_2(t),
\end{align}
where
\begin{align*}
B^m_1(t)&:=m\int_0^tE(t-s)DF(X^m(\kappa_m(s))) (E(s-\kappa_m(s))-I)X^m(\kappa_m(s))ds,\\ B^m_2(t)&:=m\int_0^tE(t-s)DF(X^m(\kappa_m(s))) \int_{\kappa_m(s)}^sE(s-r)F(X^m(\kappa_m(s))) dr ds,\\ B^m_3(t)&:=m\int_0^tE(t-s)DF(X^m(\kappa_m(s)))O_{\kappa_m(s),s} ds,\\ R^m_1(t)&:=m\int_0^tE(t-s)\int_0^1(1-\lambda) D^2F(X(s)+\lambda (X^m(s)-X(s)))\\ &\qquad\qquad(X^m(s)-X(s),X^m(s)-X(s))d\lambda ds,\\ R^m_2(t)&:=m\int_0^tE(t-s)\int_0^1(1-\lambda) D^2F(X^m(\kappa_m(s))+\lambda(X^m(s)-X^m(\kappa_m(s))))\\ &\qquad\qquad(X^m(s)-X^m(\kappa_m(s)),X^m(s)-X^m(\kappa_m(s))) d\lambda ds.
\end{align*}

In the following, we will show that the remainders $R^m_1$ and $R^m_2$ make a negligible contribution to the limiting distribution of $U^m$, and hence only terms $B^m_1$, $B^m_2$, and $B^m_3$ need to be retained.

For the remainder estimates below, set
\begin{equation}
\label{eq:nu} \nu:=\rho+\frac{\sigma}{2}-1 =\frac{\beta-(3-4H)}{2}>0,
\end{equation}
where positivity follows from \eqref{eq:new-index}.

\begin{lemma}\label{lm:RB}
Under the standing assumptions of this section, the following estimates hold
\begin{align*}
&\sup_{t\in[0,T]}\|R^m_1(t)\|_{L^{\frac{p}2}(\Omega;\V)}\le Ch,\\ &\sup_{t\in[0,T]}\|R^m_2(t)\|_{L^p(\Omega;\V)}\le Ch^\nu,\\ &\sup_{t\in[0,T]}\sum_{i=1,2,3}\|B^m_i(t)\|_{L^p(\Omega;\V)} \le C
\end{align*}
uniformly in $m$, where $\nu$ is defined in \eqref{eq:nu}.
\end{lemma}

\begin{proof}
By \eqref{eq:semi1}, \eqref{eq:F2}, the embedding $\V\hookrightarrow\V^{-\sigma}$, and Theorem \ref{tm:rate}, it holds
\[
\|R^m_1(t)\|_{L^{\frac{p}2}(\Omega;\V)} \le Cm\int_0^t(t-s)^{-\frac{\eta}{2}} \|X^m(s)-X(s)\|_{L^p(\Omega;\V)}^2 ds \le Cmh^2\le Ch.
\]
Similarly, \eqref{eq:semi1}, \eqref{eq:F2}, and Lemma \ref{lm:Xm-holder} yield
\begin{align*}
\|R^m_2(t)\|_{L^p(\Omega;\V)} &\le Cm\int_0^t(t-s)^{-\frac{\eta}{2}}\|X^m(s)-X^m(\kappa_m(s))\|_{L^{2p}(\Omega;\V)}\|X^m(s)-X^m(\kappa_m(s))\|_{L^{2p}(\Omega;\V^{-\sigma})}ds\\ &\le Cmh^{\frac{\rho}{2}+\frac{\rho+\sigma}{2}}\int_0^t(t-s)^{-\frac{\eta}{2}}ds\le Ch^{\nu},
\end{align*}
where $\eta\in[1,2H)$ is given in Assumption \ref{as:F} and $\sigma=1-\beta$ is assumed in \eqref{eq:new-index}.

The argument for \eqref{eq:I1k}, with Lemma \ref{lm:Xm-bound} in place of Theorem \ref{tm:X}, gives
\[
\|B^m_1(t)\|_{L^p(\Omega;\V)} \le Cmh^{\min\{1,\rho+\frac{\gamma}{2}\}}\le C.
\]
Similarly, \eqref{eq:I2k} applies to $X^m(\kappa_m(s))$ by Lemma \ref{lm:Xm-bound}, and yields
\[
\|B^m_2(t)\|_{L^p(\Omega;\V)}\le Cmh\le C.
\]
For $B^m_3(t)$, \eqref{eq:Um} yields
\[
\begin{aligned}
B^m_3(t) ={}&m\int_0^tE(t-s)\bigl(F(X^m(s))-F(X(s))\bigr)ds-U^m(t)\\ &-B^m_1(t)-B^m_2(t)-R^m_2(t).
\end{aligned}
\]
Hence, by \eqref{eq:F}, Theorem \ref{tm:rate}, and the preceding bounds for $B^m_1$, $B^m_2$, and $R^m_2$,
\[
\begin{aligned}
\|B^m_3(t)\|_{L^p(\Omega;\V)} \le Cm\int_0^t \|X^m(s)-X(s)\|_{L^p(\Omega;\V)}ds +m\|X^m(t)-X(t)\|_{L^p(\Omega;\V)} +C+Ch^\nu \le C,
\end{aligned}
\]
uniformly in $m\ge1$ and $t\in[0,T]$, which completes the proof.
\end{proof}

By eliminating the negligible terms $R^m_1$ and $R^m_2$, we define an auxiliary process $\tilde U^m$ as the solution of the following equation:
\begin{equation}
\label{eq:tildeUm} \tilde U^m(t)=\int_0^tE(t-s)DF(X(s))\tilde U^m(s) ds -B^m_1(t)-B^m_2(t)-B^m_3(t).
\end{equation}
Based on the boundedness of $DF$ and Lemma \ref{lm:RB}, a standard fixed-point argument yields a unique solution in the class of strongly measurable processes that are uniformly bounded in time in $L^p(\Omega;\V)$. Moreover, Gronwall's inequality gives
\begin{align}
\label{eq:tildeUm-bound} \sup_{m\ge1}\sup_{t\in[0,T]} \|\tilde U^m(t)\|_{L^p(\Omega;\V)}\le C.
\end{align}

\begin{lemma}\label{lm:U-tildeU}
Under the standing assumptions of this section, it holds
\[
\sup_{t\in[0,T]}\|U^m(t)-\tilde U^m(t)\|_{L^{\frac{p}2}(\Omega;\V)} \le C(h+h^\nu)\to0
\]
as $m\to\infty$, where $\nu$ is defined in \eqref{eq:nu}.
\end{lemma}

\begin{proof}
Subtracting \eqref{eq:tildeUm} from \eqref{eq:Um} gives
\begin{align*}
U^m(t)-\tilde U^m(t) =\int_0^tE(t-s)DF(X(s)) \bigl(U^m(s)-\tilde U^m(s)\bigr) ds +R^m_1(t)-R^m_2(t).
\end{align*}
Then \eqref{eq:F}, Lemma \ref{lm:RB}, the embedding $L^p(\Omega;\V)\hookrightarrow L^{p/2}(\Omega;\V)$, and Gronwall's inequality verify the proof.
\end{proof}

\subsection{Uniform spatial approximation}
This subsection is devoted to constructing a finite-dimensional approximation $\tilde U^{m,n}$ of $\tilde U^m$. By Assumption \ref{as:A}, choose an orthonormal eigenbasis $\{e_j\}_{j\ge1}$ of $\V$ such that $-Ae_j=\lambda_je_j$, where $0<\lambda_1\le\lambda_2\le\cdots$ and $\lambda_j\to\infty$. Define a finite-dimensional space
\[
\V_n:=\operatorname{span}\{e_1,\ldots,e_n\},
\]
and operators $A_n:=A|_{\V_n}$ and $E_n(t):=e^{tA_n}$. Let $P_n:V\to V_n$ be the projection operator defined by $ P_nv:=\sum_{j=1}^n\langle v,e_j\rangle e_j$. Thus, $E_n(t)P_n=E(t)P_n$. The spectral representation gives, for $\epsilon\ge0$,
\begin{equation}
\label{eq:spec-tail} \|(-A)^{-\epsilon}(I-P_n)\|_{\cL(\V)} =\lambda_{n+1}^{-\epsilon}.
\end{equation}
Set
\begin{equation}
\label{eq:noise-tail} \varepsilon_n:=\|(I-P_n)(-A)^{\frac{\beta-1}{2}}\Phi\|_{\HS(\U,\V)} =\left(\sum_{j>n}\Big\|\big((-A)^{\frac{\beta-1}{2}}\Phi\big)^*e_j\Big\|_U^2\right)^{\frac{1}{2}}\to0,
\end{equation}
where the convergence follows from $(-A)^{\frac{\beta-1}{2}}\Phi\in\HS(\U,\V)$. Define the projection of the stochastic convolution onto $V_n$ by
\[
O_{s,t}^n:=P_nO_{s,t} =\int_s^tE_n(t-r)P_n\Phi dW^H(r).
\]

We define a finite-dimensional process $\tilde U^{m,n}$ in $V_n$ as the solution of the following equation
\begin{equation}
\label{eq:tildeUmn} \tilde U^{m,n}(t)=\int_0^tE_n(t-s)P_nDF(X(s))\tilde U^{m,n}(s) ds -\sum_{i=1}^3I_i^{m,n}(t),
\end{equation}
where
\begin{align*}
I_1^{m,n}(t)&:=m\int_0^tE_n(t-s)P_nDF(X^m(\kappa_m(s))) \bigl(E_n(s-\kappa_m(s))-I_n\bigr)P_nX^m(\kappa_m(s)) ds,\\ I_2^{m,n}(t)&:=m\int_0^tE_n(t-s)P_nDF(X^m(\kappa_m(s))) \int_{\kappa_m(s)}^sE_n(s-r)P_nF(X^m(\kappa_m(s))) dr ds,\\ I_3^{m,n}(t)&:=m\int_0^tE_n(t-s)P_nDF(X^m(\kappa_m(s))) O_{\kappa_m(s),s}^n ds
\end{align*}
with $I_n$ being the identity operator in $V_n$.

Before establishing the existence and uniqueness of $\tilde U^{m,n}$ and its convergence to $\tilde U^m$ uniformly in $m$, we first give the estimates for the spectral projection errors associated with the stochastic integral term $B^m_3$ in \eqref{eq:tildeUm} and $I_3^{m,n}$, together with uniform bound for $I_3^{m,n}$.

\begin{lemma}\label{lm:proj}
For $Q,R\in\{I,P_n,I-P_n\}$, define
\[
\mathcal J_m(Q,R;t):= m\int_0^t E(t-s)QDF(X^m(\kappa_m(s))) \int_{\kappa_m(s)}^s E(s-r)R\Phi dW^H(r) ds.
\]
Let $\tau:=\delta{\bf 1}_{\{\beta<1\}}$ and $0<\vartheta< \min\left\{H-\frac{\tau}2,H-\frac\eta2\right\}$, where $\delta$ and $\eta$ are given in Assumption \ref{as:F}. Under the standing assumptions of this section, the following estimates hold uniformly in $m$ and $n$:
\begin{align*}
&\sup_{t\in[0,T]}\|\mathcal J_m(I-P_n,I;t)\|_{L^p(\Omega;\V)} \le C\lambda_{n+1}^{-\vartheta},\\ &\sup_{t\in[0,T]}\|\mathcal J_m(P_n,I-P_n;t)\|_{L^p(\Omega;\V)} \le C\varepsilon_n,\\ &\sup_{t\in[0,T]}\|\mathcal J_m(P_n,P_n;t)\|_{L^p(\Omega;\V)} \le C.
\end{align*}
\end{lemma}

\begin{proof}
{\it Step 1. Estimate of $\mathcal J_m(I-P_n,I;t)$.} Applying the divergence decomposition used in \eqref{eq:I3split}, we obtain
\begin{equation}
\label{eq:Jm} \mathcal J_m(I-P_n,I;t) =m\bigl(\deltaH(\mathcal U_t)+\mathcal C_t\bigr),
\end{equation}
where, for $r<t$,
\[
\mathcal U_t(r):= \int_r^{(\kappa_m(r)+h)\wedge t} E(t-s)(I-P_n)DF(X^m(\kappa_m(s)))E(s-r)\Phi ds,
\]
with $\mathcal U_t(r)=0$ for $r\ge t$, and
\begin{align*}
\mathcal C_t :={}&\alpha_H\int_0^t\int_0^{\kappa_m(s)} \int_{\kappa_m(s)}^s \sum_{j\ge1}E(t-s)(I-P_n) D^2F(X^m(\kappa_m(s)))\\ &\hspace{22mm}\bigl( D_u^jX^m(\kappa_m(s)),E(s-r)\Phi g_j \bigr)|r-u|^{2H-2} dr du ds.
\end{align*}
By Lemma \ref{lm:Xm-bound} and Proposition \ref{prop:Malliavin-Xm}, a similar argument as in the proof of Theorem \ref{tm:rate-k} shows that $\mathcal U_t\in \mathbb D^{1,p}(\HH\otimes\HS(\U,\V))$ and justifies \eqref{eq:Jm}. We next estimate $\deltaH(\mathcal U_t)$ and $\mathcal C_t$.

For the divergence term $\deltaH(\mathcal U_t)$, we need to estimate $\mathcal U_t$ and $D\mathcal U_t$ according to \eqref{eq:div-Lp}.

When $\beta<1$, the parameter conditions in \eqref{eq:new-index} imply $ \rho+\gamma-\sigma =(2\rho+\gamma)-(\rho+\sigma) \ge2-2H>0. $ Since $\sigma=1-\beta$, it follows that
\[
\|E(s-r)\Phi\|_{\mathcal L_2(\U,\V^{-\rho-\gamma})} = \|(-A)^{-\frac{\rho+\gamma-\sigma}{2}} E(s-r)(-A)^{\frac{\beta-1}{2}}\Phi\|_{\mathcal L_2(\U,\V)} \le C\|(-A)^{\frac{\beta-1}{2}}\Phi\|_{\mathcal L_2(\U,\V)}.
\]
Thus \eqref{eq:F1} and Lemma \ref{lm:Xm-bound} yield
\begin{align*}
&\quad\|E(t-s)(I-P_n)DF(X^m(\kappa_m(s)))E(s-r)\Phi \|_{L^{p}(\Omega;\mathcal L_2(\U,\V))}\\ &\le C\lambda_{n+1}^{-\vartheta} (t-s)^{-\frac{\tau}2-\vartheta}\|DF(X^m(\kappa_m(s)))E(s-r)\Phi \|_{L^{p}(\Omega;\mathcal L_2(\U,\V^{-\tau}))}\\ &\le C\lambda_{n+1}^{-\vartheta} (t-s)^{-\frac{\tau}2-\vartheta}\Big(1+\sup_{t\in[0,T]}\|X^m(t)\|_{L^{p}(\Omega;V^{\rho})}\Big)\|E(s-r)\Phi \|_{\mathcal L_2(\U,\V^{-\rho-\gamma})}\\ &\le C\lambda_{n+1}^{-\vartheta} (t-s)^{-\frac{\tau}2-\vartheta},
\end{align*}
where in the first step we used the fact $ \|E(t-s)(I-P_n)v\| \le C\lambda_{n+1}^{-\vartheta} (t-s)^{-\frac{\tau}2-\vartheta}\|v\|_{-\tau} $ according to \eqref{eq:semi1} and \eqref{eq:spec-tail}. When $\beta=1$, the same estimate holds with $\tau=0$ by the boundedness of $DF$ and $\Phi\in\mathcal L_2(\U,\V)$. Then Lemma \ref{lm:int} with $b=0$ gives
\[
\|\mathcal U_t(r)\|_{L^{p}(\Omega;\mathcal L_2(\U,\V))} \le Ch\lambda_{n+1}^{-\vartheta} (t-r)^{-\frac{\tau}2-\vartheta}
\]
under the condition $\frac{\tau}2+\vartheta<H$, which, together with \eqref{eq:double-int}, yields
\begin{align*}
\|\mathcal U_t\|_{L^{p}(\Omega;\HH\otimes\mathcal L_2(\U,\V))}^2 &\le\alpha_H\int_0^t\int_0^t\|\mathcal U_t(r)\|_{L^{p}(\Omega;\HS(\U,\V))}\|\mathcal U_t(v)\|_{L^{p}(\Omega;\HS(\U,\V))}|r-v|^{2H-2}drdv\\ &\le Ch^2\lambda_{n+1}^{-2\vartheta} \int_0^t\int_0^t (t-r)^{-\frac{\tau}2-\vartheta} (t-v)^{-\frac{\tau}2-\vartheta} |r-v|^{2H-2} dr dv \notag\\ &\le Ch^2\lambda_{n+1}^{-2\vartheta}.
\end{align*}

To estimate the Malliavin derivative of $\mathcal U_t$, the chain rule gives
\begin{align*}
D_u^i\mathcal U_t(r)g_j = \int_r^{(\kappa_m(r)+h)\wedge t} E(t-s)(I-P_n)D^2F(X^m(\kappa_m(s)))\bigl( D_u^iX^m(\kappa_m(s)),E(s-r)\Phi g_j \bigr) ds.
\end{align*}
Using \eqref{eq:F2}, \eqref{eq:semi1} with smoothing exponent $\frac{\eta}{2}+\vartheta$, and the fact that
\[
\|(-A)^{-\frac{\sigma}{2}}E(s-r)\Phi \|_{\mathcal L_2(\U,\V)} = \|E(s-r)(-A)^{\frac{\beta-1}{2}}\Phi \|_{\mathcal L_2(\U,\V)} \le C\|(-A)^{\frac{\beta-1}{2}}\Phi \|_{\mathcal L_2(\U,\V)},
\]
we obtain
\begin{align*}
\|D\mathcal U_t(r)\|_ {L^{p}(\Omega;(\HH\otimes\U)\otimes\mathcal L_2(\U,\V))} &\le C\lambda_{n+1}^{-\vartheta} \int_r^{(\kappa_m(r)+h)\wedge t} (t-s)^{-\frac{\eta}{2}-\vartheta} \|D_{\cdot}X^m(\kappa_m(s))\|_ {L^{p}(\Omega;\HH\otimes\mathcal L_2(\U,\V))}\,ds\\ &\le Ch\lambda_{n+1}^{-\vartheta} (t-r)^{-\frac{\eta}{2}-\vartheta},
\end{align*}
where the last step follows from \eqref{eq:DXm-sup} and Lemma \ref{lm:int}. Applying the $\HH$ inner product formula yields
\begin{align*}
\|D\mathcal U_t\|_{L^{p}(\Omega;(\HH\otimes\U)^{\otimes2}\otimes\V)}^2 &\le Ch^2\lambda_{n+1}^{-2\vartheta} \int_0^t\int_0^t (t-r)^{-\frac\eta2-\vartheta} (t-v)^{-\frac\eta2-\vartheta} |r-v|^{2H-2} dr dv \notag\\ &\le Ch^2\lambda_{n+1}^{-2\vartheta},
\end{align*}
since $\frac\eta2+\vartheta<H$. Consequently, \eqref{eq:div-Lp} gives
\[
\|\deltaH(\mathcal U_t)\|_{L^{p}(\Omega;\V)} \le Ch\lambda_{n+1}^{-\vartheta}.
\]

For the covariance correction $\mathcal C_t$, the assumption \eqref{eq:F2} and the Cauchy--Schwarz inequality imply
\begin{align*}
\|\mathcal C_t\|_{L^{p}(\Omega;\V)} \le C\lambda_{n+1}^{-\vartheta} \int_0^t(t-s)^{-\frac\eta2-\vartheta} \int_0^{\kappa_m(s)}\int_{\kappa_m(s)}^s \|D_uX^m(\kappa_m(s))\|_ {L^{p}(\Omega;\mathcal L_2(\U,\V))}(r-u)^{2H-2} dr du ds.
\end{align*}
Taking $s\in[t_j,t_{j+1})$ and applying \eqref{eq:DXm-norm} give
\begin{align*}
&\int_0^{t_j}\int_{t_j}^s \|D_uX^m(t_j)\|_ {L^{p}(\Omega;\mathcal L_2(\U,\V))}(r-u)^{2H-2} dr du\\ \le {}&C\int_0^{t_j}\int_{t_j}^s (t_j-u)^{\frac{\beta-1}{2}} (r-u)^{2H-2} dr du+ Ch\sum_{\ell=1}^{j-1} \int_{t_{\ell-1}}^{t_\ell} (t_\ell-u)^{\frac{\beta-1}{2}} \int_{t_j}^s(r-u)^{2H-2} dr du\\ \le{}&Ch+Ch^{\frac{\beta+1}2+2H}\sum_{\ell=1}^{j-1}(j-\ell)^{2H-2}\le Ch,
\end{align*}
where in the second step the sum is empty for $j\le1$ and we used Lemma \ref{lm:Gj} with $\sigma=1-\beta$ and $\beta>3-4H$. We also used the facts $r-u\ge(j-\ell)h$ and
\[
\sum_{\ell=1}^{j-1}(j-\ell)^{2H-2}=\sum_{k=1}^{j-1}k^{2H-2}\le\int_0^{j-1}x^{2H-2}dx\le C_p\left(\frac{T}{h}\right)^{2H-1}.
\]
Therefore,
\begin{equation*}
\|\mathcal C_t\|_{L^{p}(\Omega;\V)} \le Ch\lambda_{n+1}^{-\vartheta} \int_0^t(t-s)^{-\frac\eta2-\vartheta} ds \le Ch\lambda_{n+1}^{-\vartheta}.
\end{equation*}
Together with \eqref{eq:Jm} and $mh=T$, it proves
\[
\|\mathcal J_m(I-P_n,I;t)\|_{L^{p}(\Omega;\V)} \le C_p\lambda_{n+1}^{-\vartheta}.
\]

{\it Step 2. Estimates of $\mathcal J_m(P_n,I-P_n;t)$ and $\mathcal J_m(P_n,P_n;t)$.} Apply the same decomposition as \eqref{eq:Jm} by replacing $I-P_n$ by $P_n$ and $\Phi$ by $(I-P_n)\Phi$. Utilizing \eqref{eq:noise-tail} and the fact that $P_n$ is a contraction, we derive for the corresponding $\deltaH(\mathcal U_t)$ and $\mathcal C_t$ that
\begin{align*}
\|\deltaH(\mathcal U_t)\|_{L^{p}(\Omega;\V)} &\le C\left( \|\mathcal U_t\|_{L^{p}(\Omega;\HH\otimes\mathcal L_2(\U,\V))} +\|D\mathcal U_t\|_{L^{p}(\Omega;(\HH\otimes\U)^{\otimes2}\otimes\V)} \right)\le Ch\varepsilon_n,\\ \|\mathcal C_t\|_{L^{p}(\Omega;\V)}&\le Ch\varepsilon_n
\end{align*}
under the condition $\frac{\tau}2<H$ and $\frac{\eta}2<H$. It then implies
\[
\|\mathcal J_m(P_n,I-P_n;t)\|_{L^{p}(\Omega;\V)} \le C\varepsilon_n.
\]

Finally, the fact $B^m_3(t)=\mathcal J_m(I,I;t)$ and the bilinearity of $\mathcal J_m(Q,R;t)$ with respect to $Q$ and $R$ give
\[
\mathcal J_m(P_n,P_n;t) =B^m_3(t)-\mathcal J_m(I-P_n,I;t) -\mathcal J_m(P_n,I-P_n;t).
\]
The preceding estimates and Lemma \ref{lm:RB} yield
\[
\|\mathcal J_m(P_n,P_n;t)\|_{L^p(\Omega;\V)} \le C\bigl(1+\lambda_{n+1}^{-\vartheta} +\varepsilon_n\bigr) \le C,
\]
uniformly in $m,n$ and $t\in[0,T]$.
\end{proof}

The boundedness of $\mathcal J_m(P_n,P_n;\cdot)$ given in the above lemma allows us to establish the well-posedness of \eqref{eq:tildeUmn} and a bound uniform in $m$ and $n$.

\begin{lemma}\label{lm:tildeUmn}
Under the standing assumptions of this section, equation \eqref{eq:tildeUmn} has a unique solution in the class of strongly measurable $\V_n$-valued processes uniformly bounded in time in $L^p(\Omega;\V)$. Moreover,
\[
\sup_{m,n\ge1}\sup_{t\in[0,T]} \|\tilde U^{m,n}(t)\|_{L^p(\Omega;\V)}\le C.
\]
\end{lemma}

\begin{proof}
For $I_1^{m,n}$ and $I_2^{m,n}$ in \eqref{eq:tildeUmn}, the arguments for $B^m_1$ and $B^m_2$ in Lemma \ref{lm:RB} apply, since $P_n$ commutes with $E(t)$ and fractional powers of $-A$, and is a contraction on each $\V^r$. For $I_3^{m,n}=\mathcal J_m(P_n,P_n;\cdot)$, Lemma \ref{lm:proj} applies. Hence,
\[
\sup_{m,n\ge1}\sup_{t\in[0,T]} \sum_{i=1}^3\|I_i^{m,n}(t)\|_{L^p(\Omega;\V)} \le C.
\]
Together with the boundedness of $DF$, this yields, by a standard fixed-point argument, a unique solution in the class of strongly measurable $\V_n$-valued processes uniformly bounded in time in $L^p(\Omega;\V)$.

Taking $L^p(\Omega;\V)$-norms in \eqref{eq:tildeUmn} gives
\[
\|\tilde U^{m,n}(t)\|_{L^p(\Omega;\V)} \le C+C\int_0^t \|\tilde U^{m,n}(s)\|_{L^p(\Omega;\V)}\,ds.
\]
The uniform bound then follows from Gronwall's inequality.
\end{proof}

We now compare $\tilde U^{m,n}$ with $\tilde U^m$. The uniform bound above controls the projection error in the linear term, while Lemma \ref{lm:proj} handles the stochastic forcing. The remaining deterministic terms can be estimated using the spatial regularity of $X^m$.

\begin{lemma}\label{lm:tildeUmn-tildeUm}
Under the standing assumptions of this section, then
\[
\lim_{n\to\infty}\sup_{m\ge1}\sup_{t\in[0,T]} \|\tilde U^{m,n}(t)-\tilde U^m(t)\|_{L^p(\Omega;\V)}=0.
\]
\end{lemma}

\begin{proof}
Subtracting \eqref{eq:tildeUm} from \eqref{eq:tildeUmn} and using $E_n(t)P_n=E(t)P_n$, we obtain
\begin{equation*}
\tilde U^{m,n}(t)-\tilde U^m(t) =\int_0^tE(t-s)DF(X(s)) \bigl(\tilde U^{m,n}(s)-\tilde U^m(s)\bigr) ds+\sum_{i=0}^3S_i^{m,n}(t),
\end{equation*}
where
\[
S_0^{m,n}(t):= \int_0^tE(t-s)(P_n-I)DF(X(s)) \tilde U^{m,n}(s) ds,\qquad S_i^{m,n}:=B^m_i-I_i^{m,n}, \qquad i=1,2,3.
\]
By the boundedness of $DF$ and Gronwall's inequality, it suffices to show that these four projection errors $S_i^{m,n}(t)$, $i=0,1,2,3$, converge to zero in $L^p(\Omega;\V)$, uniformly in $m$ and $t\in[0,T]$.

For $S_0^{m,n}$, by \eqref{eq:semi1}, \eqref{eq:spec-tail}, and Lemma \ref{lm:tildeUmn}, it holds for any $\epsilon_0\in(0,1)$ that
\begin{equation}
\label{eq:S0} \sup_{m\ge1}\sup_{t\in[0,T]}\|S_0^{m,n}(t)\|_{L^p(\Omega;\V)} \le C\lambda_{n+1}^{-\epsilon_0} \sup_{t\in[0,T]}\int_0^t(t-s)^{-\epsilon_0} ds \le C\lambda_{n+1}^{-\epsilon_0}.
\end{equation}

For $S_1^{m,n}$, we split it into two terms
\begin{align*}
S_1^{m,n}(t) &=m\int_0^tE(t-s)(I-P_n)DF(X^m(\kappa_m(s))) (E(s-\kappa_m(s))-I)X^m(\kappa_m(s)) ds\\ &\quad+m\int_0^tE(t-s)P_nDF(X^m(\kappa_m(s))) (E(s-\kappa_m(s))-I)(I-P_n)X^m(\kappa_m(s)) ds\\ &=:S_{1,1}^{m,n}(t)+S_{1,2}^{m,n}(t).
\end{align*}
Based on \eqref{eq:semi1}, \eqref{eq:F1}, \eqref{eq:spec-tail}, Lemma \ref{lm:Xm-bound}, and the estimate $\|(E(t)-I)v\|_{-\rho-\gamma} \le C_pt\|v\|_\rho$ for $v\in\V^\rho$ since $\rho+\frac\gamma2\ge1$, we derive for any $\epsilon_1\in(0,1-\frac\delta2)$ that
\[
\|S_{1,1}^{m,n}\|_{L^p(\Omega;\V)}\le Cm\lambda_{n+1}^{-\epsilon_1} \int_0^t(t-s)^{-\frac\delta2-\epsilon_1} (s-\kappa_m(s)) ds \le C\lambda_{n+1}^{-\epsilon_1}.
\]
A similar argument applied to $S_{1,2}^{m,n}$ yields
\begin{align*}
\|S_{1,2}^{m,n}\|_{L^p(\Omega;\V)}&\le Cm\int_0^t(t-s)^{-\frac{\delta}2}(s-\kappa_m(s))\left(1+\|X^m(\kappa_m(s))\|_{L^{2p}(\Omega;\V^\rho)}\right)\\ &\qquad\times\|(I-P_n)X^m(\kappa_m(s))\|_{L^{2p}(\Omega;\V^\rho)}ds\\ &\le Ca_n,
\end{align*}
where
\[
a_n:= \sup_{m\ge1}\sup_{t\in[0,T]} \|(I-P_n)X^m(t)\|_{L^{2p}(\Omega;\V^\rho)}.
\]
Choosing $\epsilon\in(0,1-\frac{\rho}2)$, multiplying both sides of \eqref{eq:Xm} by $I-P_n$, and utilizing Lemma \ref{lm:Xm-bound} give
\begin{align}
\label{eq:an} a_n &\le \|(I-P_n)X_0\|_{L^{2p}(\Omega;\V^\rho)} +C\lambda_{n+1}^{-\epsilon} \sup_{t\in[0,T]}\int_0^t(t-s)^{-\frac\rho2-\epsilon} ds +C\varepsilon_n\notag\\ &\le \|(I-P_n)X_0\|_{L^{2p}(\Omega;\V^\rho)} +C\lambda_{n+1}^{-\epsilon} +C\varepsilon_n \to0,\qquad n\to\infty.
\end{align}
Here the initial tail tends to zero by strong convergence of $P_n$ on $\V^\rho$ and dominated convergence. Combining the above estimates yields
\begin{equation}
\label{eq:S1} \sup_{m\ge1}\sup_{t\in[0,T]}\|S_1^{m,n}(t)\|_{L^p(\Omega;\V)} \le C\bigl(\lambda_{n+1}^{-\epsilon_1}+a_n\bigr) \to0,\qquad n\to\infty.
\end{equation}

For $S_2^{m,n}$, similarly write
\begin{align*}
S_2^{m,n}(t) ={}&m\int_0^tE(t-s)(I-P_n)DF(X^m(\kappa_m(s))) \int_{\kappa_m(s)}^s E(s-r)F(X^m(\kappa_m(s))) dr ds\\ &+m\int_0^tE(t-s)P_nDF(X^m(\kappa_m(s))) \int_{\kappa_m(s)}^s E(s-r)(I-P_n)F(X^m(\kappa_m(s))) dr ds.
\end{align*}
Then \eqref{eq:F}, \eqref{eq:F1}, H\"older's inequality, Lemma \ref{lm:Xm-bound}, and \eqref{eq:spec-tail} with smoothing exponents $\epsilon=\frac12$ and $\epsilon=\frac{\rho+\gamma}{2}$, respectively, give
\begin{align}
\label{eq:S2} \|S_2^{m,n}(t)\|_{L^p(\Omega;\V)} &\le Cmh\left( \lambda_{n+1}^{-\frac12} \int_0^t(t-s)^{-\frac12} ds +\lambda_{n+1}^{-\frac{\rho+\gamma}{2}} \int_0^t(t-s)^{-\frac\delta2} ds \right) \notag\\ &\le C\left( \lambda_{n+1}^{-\frac12} +\lambda_{n+1}^{-\frac{\rho+\gamma}{2}} \right).
\end{align}

For $S_3^{m,n}$, the fact $S_3^{m,n}(t) =\mathcal J_m(I-P_n,I;t) +\mathcal J_m(P_n,I-P_n;t)$ and Lemma \ref{lm:proj} imply
\begin{equation}
\label{eq:S3} \sup_{m\ge1}\sup_{t\in[0,T]}\|S_3^{m,n}(t)\|_{L^p(\Omega;\V)} \le C\bigl(\lambda_{n+1}^{-\vartheta}+\varepsilon_n\bigr) \to0,\qquad n\to\infty
\end{equation}
with $\vartheta$ being given in Lemma \ref{lm:proj}.

Combining the estimates \eqref{eq:S0}--\eqref{eq:S3} completes the proof by applying Gronwall's inequality.
\end{proof}

Lemmas \ref{lm:U-tildeU} and \ref{lm:tildeUmn-tildeUm} reduce the study of $U^m$ to the finite-dimensional auxiliary processes $\tilde U^{m,n}$, with an approximation error that vanishes uniformly in the time discretization. Identification of a limiting distribution requires, in addition, the limit of these processes for each fixed $n$.

\subsection{Limit for a fixed spatial dimension}
For a Banach space $(G,\|\cdot\|_G)$ and $\alpha\in(0,1)$, write $C^\alpha([0,T];G)$ for the space of $\alpha$-H\"older continuous functions, equipped with the norm
\[
\|f\|_{C^{\alpha}([0,T];G)} := \sup_{t\in[0,T]}\|f(t)\|_G +[f]_{C^{\alpha}([0,T];G)}
\]
with $[f]_{C^{\alpha}([0,T];G)}:=\sup_{0\le s<t\le T} \frac{\|f(t)-f(s)\|_G}{|t-s|^{\alpha}}$.

\begin{proposition}\label{prop:finite-limit}
Under the standing assumptions of this section, for every fixed $n\ge1$,
\[
\tilde U^{m,n}\to \tilde{U}^{\infty,n} \quad\hbox{in probability in }C([0,T];\V_n) \quad\hbox{as }m\to\infty,
\]
where $\tilde U^{m,n}$ is defined in \eqref{eq:tildeUmn} and $\tilde{U}^{\infty,n}$ is the unique continuous $\V_n$-valued solution of
\begin{align}
\label{eq:Un-limit} \tilde{U}^{\infty,n}(t)={}&\int_0^t E_n(t-s)P_nDF(X(s))\tilde{U}^{\infty,n}(s)ds -\frac{T}{2}\int_0^t E_n(t-s)P_nDF(X(s))A_nP_nX(s)ds\notag\\ &-\frac{T}{2}\int_0^t E_n(t-s)P_nDF(X(s))P_nF(X(s))ds\notag\\ &-\frac{T}{2}\int_0^t E_n(t-s)P_nDF(X(s))dW^{H,n}_\Phi(s)
\end{align}
with $ W^{H,n}_\Phi(t):=\sum_{j\ge1}P_n\Phi g_j\beta_j^H(t). $ If, additionally, $X_0\in L^\infty(\Omega;\V^\rho)$, the convergence also holds in $L^\ell(\Omega;C([0,T];\V))$ for every $\ell\in[1,\infty)$.
\end{proposition}

\begin{proof}
{\it Step 1. Bounded initial value case.} We first assume $X_0\in L^\infty(\Omega;\V^\rho)$ and fix an arbitrary $\ell\in[1,\infty)$. For every $0<a<\frac{\rho}{2}$, Theorem \ref{tm:X} and Lemma \ref{lm:Xm-holder} yield for any $r>\max\{\ell,(\frac{\rho}{2}-a)^{-1}\}$
\[
\|X(t)-X(s)\|_{L^r(\Omega;\V)}+\sup_{m\ge1}\|X^m(t)-X^m(s)\|_{L^r(\Omega;\V)}\le C|t-s|^{\frac{\rho}2}.
\]
Since $a<\frac{\rho}2-\frac1r$ and $\ell<r$, the Kolmogorov continuity theorem implies that $X$ and $X^m$ have continuous modifications, still denoted by the same symbols, which satisfy
\begin{equation}
\label{eq:holder-X} \|X\|_{L^\ell(\Omega;C^a([0,T];\V))} +\sup_{m\ge1}\|X^m\|_{L^\ell(\Omega;C^a([0,T];\V))}\le C.
\end{equation}

For $a_0\in(a,\frac{\rho}{2})$, \eqref{eq:mild}, \eqref{eq:Xm}, and \eqref{eq:F} give
\begin{align*}
\sup_{r\in[0,t]}\|X^m(r)-X(r)\|&\le C\int_0^t\|X^m(\kappa_m(s))-X(\kappa_m(s))\|ds+C\int_0^t\|X(\kappa_m(s))-X(s)\|ds\\ &\le C\int_0^t\sup_{r\in[0,s]}\|X^m(r)-X(r)\|ds+C \sup_{s\in[0,T]}\|X(\kappa_m(s))-X(s)\|\\ &\le C\int_0^t\sup_{r\in[0,s]}\|X^m(r)-X(r)\|ds+Ch^{a_0}[X]_{C^{a_0}([0,T];\V)}.
\end{align*}
Then Gronwall's inequality, together with \eqref{eq:holder-X}, leads to
\begin{align}
\label{eq:Xm-conv-C} \|X^m-X\|_{L^\ell(\Omega;C([0,T];\V))}\le Ch^{a_0}.
\end{align}

{\it Step 1.1. Convergence of $I_1^{m,n}$ and $I_2^{m,n}$.} For fixed $n$, we define continuous mappings $G_i:V\to V_n$ by
\[
G_1(x):=P_nDF(x)A_nP_nx,\quad G_2(x):=P_nDF(x)P_nF(x),\qquad
\]
which, due to the boundedness of $A_n$ and Assumption \ref{as:F}, satisfy
\begin{align}
\label{eq:Gi} \|G_i(x)\|\le C(1+\|x\|),\qquad\|G_i(x)-G_i(y)\|\le C(1+\|x\|+\|y\|)\|x-y\|
\end{align}
for $i=1,2$. Here and in what follows, the constant $C$ may depend on $n$ but is independent of $m$. Then one can rewrite $I_1^{m,n}$ and $I_2^{m,n}$ as
\begin{align*}
I_1^{m,n}(t)&=m\int_0^t(s-\kappa_m(s))E_n(t-s)G_1(X^m(\kappa_m(s))) ds+m\int_0^tE_n(t-s)P_nDF(X^m(\kappa_m(s)))\\ &\qquad\times\bigl(E_n(s-\kappa_m(s))-I_n-(s-\kappa_m(s))A_n\bigr)P_nX^m(\kappa_m(s)) ds\\ &=:I_{1,1}^{m,n}+I_{1,2}^{m,n}(t),\\ I_2^{m,n}(t)&=m\int_0^t(s-\kappa_m(s))E_n(t-s)G_2(X^m(\kappa_m(s))) ds+m\int_0^tE_n(t-s)P_nDF(X^m(\kappa_m(s)))\\ &\qquad\times\left[\int_{\kappa_m(s)}^sE_n(s-r)dr-(s-\kappa_m(s))I_n\right]P_nF(X^m(\kappa_m(s))) ds\\ &=:I_{2,1}^{m,n}+I_{2,2}^{m,n}(t).
\end{align*}

For $I_{i,2}^{m,n}$, $i=1,2$, due to the boundedness of $A_n$, it holds uniformly for $0\le v\le h$ that
\[
\|E_n(v)-I_n-vA_n\|_{\mathcal L(\V_n)}\le C v^2, \qquad \left\|\int_0^v E_n(r)dr-vI_n\right\|_{\mathcal L(\V_n)} \le Cv^2.
\]
The above estimates, together with \eqref{eq:F} and Lemma \ref{lm:Xm-bound}, yield
\begin{align}
\label{eq:I12-I22} &\|I_{1,2}^{m,n}\|_{L^\ell(\Omega;C([0,T];V))} +\|I_{2,2}^{m,n}\|_{L^\ell(\Omega;C([0,T];V))}\notag\\ \le{}& Cm\int_0^T(s-\kappa_m(s))^2 \left( \|X^m(\kappa_m(s))\|_{L^\ell(\Omega;V)} +\|F(X^m(\kappa_m(s)))\|_{L^\ell(\Omega;V)} \right)ds\notag\\ \le{}& Cm\int_0^T(s-\kappa_m(s))^2ds =Cm\sum_{j=0}^{m-1}\int_{t_j}^{t_{j+1}}(s-t_j)^2ds\le Ch\to0,\qquad m\to\infty.
\end{align}

For $I_{i,1}^{m,n}$, $i=1,2$, we claim that $I_{i,1}^{m,n}\to\frac{T}2\int_0^{\cdot}E_n(\cdot-s)G_i(X(s))ds$ in $L^\ell(\Omega;C([0,T];V))$ as $m\to\infty$. In fact, we consider the following error
\begin{align*}
I_{i,1}^{m,n}(t)-\frac{T}2\int_0^tE_n(t-s)G_i(X(s))ds ={}&m\int_0^t(s-\kappa_m(s))E_n(t-s) \bigl[G_i(X^m(\kappa_m(s)))-G_i(X(s))\bigr]ds\\ &+\int_0^t \left(m(s-\kappa_m(s))-\frac{T}{2}\right) E_n(t-s)G_i(X(s))ds\\ =:{}&\mathcal E_1^{m,n,i}(t)+\mathcal E_2^{m,n,i}(t),
\end{align*}
where, according to \eqref{eq:holder-X}, \eqref{eq:Xm-conv-C}, and \eqref{eq:Gi},
\begin{align*}
\|\mathcal E_1^{m,n,i}\|_{L^\ell(\Omega;C([0,T];V))} &\le C\left\| \sup_{s\in[0,T]} \|G_i(X^m(\kappa_m(s)))-G_i(X(s))\| \right\|_{L^\ell(\Omega)}\\ &\le C\left(1+\|X^m\|_{L^{2\ell}(\Omega;C([0,T];V))}+\|X\|_{L^{2\ell}(\Omega;C([0,T];V))}\right)\\ &\qquad\times\|X^m(\kappa_m(\cdot))-X(\cdot)\|_{L^{2\ell}(\Omega;C([0,T];V))}\\ &\le Ch^a\to0,\qquad m\to\infty.
\end{align*}
It then suffices to investigate the convergence of $\mathcal E_2^{m,n,i}$. Due to the facts
\begin{align}
\label{eq:m-T} \int_{t_j}^{t_{j+1}}\left(m(s-t_j)-\frac{T}{2}\right)ds=0, \qquad \left|m(s-\kappa_m(s))-\frac{T}{2}\right|\le\frac{T}{2}
\end{align}
and the expression
\begin{align*}
\mathcal E_2^{m,n,i}&=\sum_{j=0}^{\lfloor\frac{t}{h}\rfloor-1}\int_{t_j}^{t_{j+1}}\left(m(s-t_j)-\frac{T}{2}\right) \Bigl[ E_n(t-s)\bigl(G_i(X(s))-G_i(X(t_j))\bigr)\\ &\hspace{48mm} +\bigl(E_n(t-s)-E_n(t-t_j)\bigr)G_i(X(t_j)) \Bigr]ds\\ &\quad+\int_{t_{\lfloor\frac{t}{h}\rfloor}}^t \left(m(s-t_{\lfloor\frac{t}{h}\rfloor})-\frac{T}{2}\right)E_n(t-s)G_i(X(s))ds,
\end{align*}
we derive from \eqref{eq:holder-X} and \eqref{eq:Gi} that
\begin{align*}
\|\mathcal E_2^{m,n,i}\|_{L^\ell(\Omega;C([0,T];V))}\le{}&C\sum_{j=0}^{m-1}\int_{t_j}^{t_{j+1}}\left\|\left(1+\|X\|_{C([0,T];\V)}\right)[X]_{C^a([0,T];\V)}\right\|_{L^\ell(\Omega)}(s-t_j)^ads\\ &+C\sum_{j=0}^{m-1}\int_{t_j}^{t_{j+1}}\left(1+\|X\|_{L^\ell(\Omega;C([0,T];V))}\right)(s-t_j)ds\\ &+Ch\left(1+\|X\|_{L^\ell(\Omega;C([0,T];V))}\right)\\ \le{}&C(h^a+h)\to0,\qquad m\to\infty.
\end{align*}
It verifies the claim, and implies the following convergence in $L^\ell(\Omega;C([0,T];\V))$ combining with \eqref{eq:I12-I22}:
\begin{equation}
\label{eq:finite-drift-limits}
\begin{aligned}
I_1^{m,n}(t)&\to \frac{T}{2}\int_0^t E_n(t-s)P_nDF(X(s))A_nP_nX(s)ds,\qquad m\to\infty,\\ I_2^{m,n}(t)&\to \frac{T}{2}\int_0^t E_n(t-s)P_nDF(X(s))P_nF(X(s))ds,\qquad m\to\infty.
\end{aligned}
\end{equation}

{\it Step 1.2. Convergence of $I_3^{m,n}$.} Integration by parts gives
\[
O^n_{\kappa_m(s),s} =W^{H,n}_\Phi(s)-W^{H,n}_\Phi(\kappa_m(s)) +\int_{\kappa_m(s)}^s A_nE_n(s-r) \bigl(W^{H,n}_\Phi(r) -W^{H,n}_\Phi(\kappa_m(s))\bigr)dr.
\]
Define $Z^{m,n}(t) :=m\int_0^t \bigl(W^{H,n}_\Phi(s) -W^{H,n}_\Phi(\kappa_m(s))\bigr)ds$ for simplicity. Then
\begin{align}
\label{eq:I3-error-decomposition} &I_3^{m,n}(t) -\frac{T}{2}\int_0^t E_n(t-s)P_nDF(X(s))dW^{H,n}_\Phi(s) \notag\\ ={}&m\int_0^t E_n(t-s) P_nDF(X^m(\kappa_m(s))) \int_{\kappa_m(s)}^s A_nE_n(s-r) \bigl(W^{H,n}_\Phi(r) -W^{H,n}_\Phi(\kappa_m(s))\bigr)drds \notag\\ &+m\int_0^t E_n(t-s)P_n \bigl(DF(X^m(\kappa_m(s)))-DF(X^m(s))\bigr) \bigl(W^{H,n}_\Phi(s) -W^{H,n}_\Phi(\kappa_m(s))\bigr)ds \notag\\ &+\int_0^t E_n(t-s)P_n \bigl(DF(X^m(s))-DF(X(s))\bigr)dZ^{m,n}(s) \notag\\ &+\int_0^t E_n(t-s)P_nDF(X(s)) d\left(Z^{m,n}-\frac{T}{2}W^{H,n}_\Phi\right)(s) \notag\\ =:{}&R_1^{m,n}(t)+R_2^{m,n}(t) +R_3^{m,n}(t)+R_4^{m,n}(t).
\end{align}

We estimate the four terms in $L^\ell(\Omega;C([0,T];\V))$, for any finite $l\ge1$. The estimate
\begin{align}
\label{eq:W-Phi} \|W^{H,n}_\Phi(t)-W^{H,n}_\Phi(s)\|_{L^\ell(\Omega;\V_n)} \le C \|W^{H,n}_\Phi(t)-W^{H,n}_\Phi(s)\|_{L^2(\Omega;\V_n)} \le C|t-s|^H,
\end{align}
together with \eqref{eq:F} and boundedness of $A_n$ for fixed $n$, gives
\begin{align*}
\|R_1^{m,n}\|_{L^\ell(\Omega;C([0,T];\V))} &\le Cm\int_0^T \int_{\kappa_m(s)}^s (r-\kappa_m(s))^Hdrds \le Ch^H\to0,\qquad m\to\infty.
\end{align*}

Since $H+\frac{\rho}{2}=1+\nu>1$, choose $0<a<\frac{\rho}{2}$ and $0<b<H$ with $a+b>1$. Similar as above, \eqref{eq:F2}, \eqref{eq:holder-X}, \eqref{eq:W-Phi}, and the fact $|t-s|^H\le T^{H-b}|t-s|^b$ imply
\begin{align*}
\|R_2^{m,n}\|_{L^\ell(\Omega;C([0,T];\V))} &\le Cm\int_0^T \|X^m(s)-X^m(\kappa_m(s))\|_{L^{2\ell}(\Omega;\V)} \|W^{H,n}_\Phi(s)-W^{H,n}_\Phi(\kappa_m(s))\|_ {L^{2\ell}(\Omega;\V_n)}ds\\ &\le Cm\int_0^T (s-\kappa_m(s))^{a+b}ds\le Ch^{a+b-1}\to0,\qquad m\to\infty.
\end{align*}

To estimate $R_3^{m,n}$, we first establish the uniform increment estimate
\begin{equation}
\label{eq:Z-increments} \|Z^{m,n}(t)-Z^{m,n}(s)\|_{L^\ell(\Omega;\V_n)} \le C|t-s|^H.
\end{equation}
For $t-s\le h$, Minkowski's inequality and \eqref{eq:W-Phi} yield
\begin{align*}
\|Z^{m,n}(t)-Z^{m,n}(s)\|_{L^2(\Omega;\V_n)} \le m\int_s^t \|W^{H,n}_\Phi(r) -W^{H,n}_\Phi(\kappa_m(r))\|_{L^2(\Omega;\V_n)}dr\le C(t-s)^H.
\end{align*}
For $t-s>h$, we have, for almost every $r\in[0,T]$,
\[
0\le m\int_s^t\mathbf1_{[\kappa_m(u),u]}(r)du = m\bigl|[s,t]\cap[r,\kappa_m(r)+h]\bigr| \le T\mathbf1_{[\kappa_m(s),t]}(r)
\]
with $|[s,t]|$ denoting the length of the interval. Thus, by \eqref{eq:HH} and Fubini's theorem,
\begin{align*}
&\quad\|Z^{m,n}(t)-Z^{m,n}(s)\|_{L^2(\Omega;\V_n)}^2\\ &=H(2H-1)\|P_n\Phi\|_{\mathcal L_2(\U,\V)}^2m^2 \int_s^t\int_s^t \int_{\kappa_m(u)}^u\int_{\kappa_m(v)}^v |r-r'|^{2H-2}dr'drdvdu\\ &=H(2H-1)\|P_n\Phi\|_{\mathcal L_2(\U,\V)}^2 \int_0^T\int_0^T \left(m\int_s^t \mathbf1_{[\kappa_m(u),u]}(r)du\right)\left(m\int_s^t \mathbf1_{[\kappa_m(v),v]}(r')dv\right) |r-r'|^{2H-2}dr'dr\\ &\le C \int_{\kappa_m(s)}^t\int_{\kappa_m(s)}^t |r-r'|^{2H-2}dr'dr =\frac{C}{H(2H-1)} (t-\kappa_m(s))^{2H}\\ &\le C\bigl((t-s)+h\bigr)^{2H} \le C(t-s)^{2H}.
\end{align*}
Then the Kahane--Khintchine inequality \eqref{eq:KK} proves \eqref{eq:Z-increments} for all finite $l\ge1$. Similar to \eqref{eq:holder-X}, applying the Kolmogorov continuity theorem to \eqref{eq:Z-increments} and \eqref{eq:W-Phi} gives
\begin{align}
\label{eq:Ybound} \sup_{m\ge1} \|Z^{m,n}\|_{L^\ell(\Omega;C^{b}([0,T];\V_n))} +\|W^{H,n}_\Phi\|_{L^\ell(\Omega;C^{b}([0,T];\V_n))} \le C.
\end{align}
By Young's inequality (cf. \cite[(4.3)]{FH20})
\begin{equation}
\label{eq:Young} \left\|\int_u^v f(s)dz(s) -f(u)(z(v)-z(u))\right\| \le C_{a,b}[f]_{C^a}[z]_{C^b}(v-u)^{a+b},\qquad f\in C^a,~z\in C^b
\end{equation}
with $a+b>1$ and \eqref{eq:Ybound}, we obtain for $0<a<a_0<\frac{\rho}2$
\begin{align*}
\|R_3^{m,n}\|_{L^\ell(\Omega;C([0,T];\V))} &\le C \left\|P_n\bigl(DF(X^m)-DF(X)\bigr)|_{\V_n}\right\|_ {L^{2\ell}(\Omega;C^a([0,T];\mathcal L(\V_n)))}\|Z^{m,n}\|_{L^{2\ell}(\Omega;C^b([0,T];\V_n))}\\ &\le Ch^{a_0-a}\to0,\qquad m\to\infty,
\end{align*}
where in the last step we used the fact
\begin{align*}
&\left\|P_n\bigl(DF(X^m)-DF(X)\bigr)|_{\V_n}\right\|_ {L^{2\ell}(\Omega;C^a([0,T];\mathcal L(\V_n)))} \notag\\ \le{}& C \|X^m-X\|_{L^{4\ell}(\Omega;C([0,T];\V))}^{1-\frac{a}{a_0}} \left( \|X^m\|_{L^{4\ell}(\Omega;C^{a_0}([0,T];\V))} +\|X\|_{L^{4\ell}(\Omega;C^{a_0}([0,T];\V))} \right)^{\frac{a}{a_0}} \notag\\ \le{}& Ch^{a_0-a}
\end{align*}
Here the mean-value estimate for $DF(X^m)-DF(X)$ follows from the directional differentiability in Assumption \ref{as:F}, the bound \eqref{eq:F2}, and the finite-dimensional projection $P_n$. Based on \eqref{eq:F2}, \eqref{eq:holder-X}, \eqref{eq:Xm-conv-C}, and the H\"older interpolation inequality in Banach space $G$:
\begin{equation}
\label{eq:interpolation} \|v\|_{C^a([0,T];G)} \le C\|v\|_{C([0,T];G)}^{1-\frac{a}{a_0}} \|v\|_{C^{a_0}([0,T];G)}^{\frac{a}{a_0}}.
\end{equation}

For $R_4^{m,n}$, we first establish the convergence
\begin{equation}
\label{eq:averaged-noise} Z^{m,n}\to\frac{T}{2}W^{H,n}_\Phi \quad\hbox{in } L^\ell(\Omega;C^b([0,T];\V_n)).
\end{equation}
Set $w_m(r):=m(\kappa_m(r)+h-r)-\frac{T}{2}$. For $t\in[t_k,t_{k+1})$, stochastic Fubini theorem gives
\begin{align*}
Y^{m,n}(t):=Z^{m,n}(t)-\frac{T}{2}W^{H,n}_\Phi(t) =\int_0^t w_m(r)dW^{H,n}_\Phi(r)-m(t_{k+1}-t) \bigl(W^{H,n}_\Phi(t)-W^{H,n}_\Phi(t_k)\bigr).
\end{align*}
At $t=T$, the last term is absent. Hence, \eqref{eq:HH} and \eqref{eq:W-Phi} imply
\begin{align*}
\label{eq:Z-pointwise} \left\|Y^{m,n}(t)\right\|_{L^2(\Omega;\V_n)}^2 \le C \left|\int_0^t\int_0^t w_m(r)w_m(v)|r-v|^{2H-2}drdv\right| +Ch^{2H}.
\end{align*}
For a fixed $t$, set $K(r,v) = {\bf 1}_{[0,t]^2} (r,v)|r-v|^{2H-2}$. For any $\epsilon>0$, choose $N=N(\epsilon)\in\mathbb N$ sufficiently large and partition $[0,T]$ into $N$ intervals $J_1,\dots,J_N$ of equal length, so that the piecewise constant approximation
\[
K_\epsilon(r,v):=\sum_{i,j=1}^N c_{ij}\mathbf1_{J_i}(r)\mathbf1_{J_j}(v), \qquad c_{ij}:=\frac{1}{|J_i||J_j|} \int_{J_i}\int_{J_j}K(r,v)dvdr,
\]
satisfies $\|K-K_\epsilon\|_{L^1([0,T]^2)}<\epsilon$. Moreover, properties analogous to \eqref{eq:m-T} hold, and yield $\|w_m\|_{L^\infty(0,T)}\le\frac{T}{2}$ and $ \left|\int_J w_m(r)dr\right|\le Th$ for any interval $J\subset[0,T]$. Then it holds
\begin{align*}
\left|\int_0^T\int_0^Tw_m(r)w_m(v)K(r,v)drdv\right| &\le \int_0^T\int_0^T|w_m(r)w_m(v)||K-K_\epsilon|(r,v)drdv\\ &\quad+ \left|\int_0^T\int_0^Tw_m(r)w_m(v)K_\epsilon(r,v)drdv\right|\\ &\le\frac{T^2}{4}\epsilon +T^2h^2\sum_{i,j=1}^{N}|c_{ij}|.
\end{align*}
Letting $m\to\infty$ and then $\epsilon\to0$ yields $Y^{m,n}(t)\to0$ in $L^2(\Omega;\V_n)$ for every fixed $t\in[0,T]$. Since $Y^{m,n}(t)$ is a centered Gaussian random variable, the Kahane–Khintchine inequality \eqref{eq:KK} imply that $Y^{m,n}(t)\to0$ in $L^\ell(\Omega;\V_n)$ for every finite $l\ge1$.

For $N\ge1$, set $r_j=\frac{jT}{N}$, $j=0,\dots,N$. For every $t\in[r_j,r_{j+1}]$,
\[
\|Y^{m,n}(t)\| \le \|Y^{m,n}(r_j)\| +\left(\frac{T}{N}\right)^{b_0} [Y^{m,n}]_{C^{b_0}([0,T];\V_n)},
\]
which, together with \eqref{eq:Ybound}, yields
\begin{align}
\label{eq:mesh-upgrade} \|Y^{m,n}\|_{L^\ell(\Omega;C([0,T];\V_n))} &\le \left(\sum_{j=0}^N \mathbb E\left[\|Y^{m,n}(r_j)\|^{\ell}\right]\right)^{\frac1{\ell}} + \left(\frac{T}{N}\right)^{b_0} \|[Y^{m,n}]_{C^{b_0}([0,T];\V_n)}\|_{L^\ell(\Omega)} \notag\\ &\le (N+1)^{\frac1{\ell}} \max_{0\le j\le N} \|Y^{m,n}(r_j)\|_{L^\ell(\Omega;\V_n)} +CN^{-b_0}.
\end{align}
For fixed $N$, the pointwise convergence $Y^{m,n}(t)\to0$ in $L^\ell(\Omega;\V_n)$ established above therefore gives
\[
\limsup_{m\to\infty} \|Y^{m,n}\|_{L^\ell(\Omega;C([0,T];\V_n))} \le CN^{-b_0}.
\]
Letting $N\to\infty$, we obtain $Y^{m,n}\to0$ in $L^\ell(\Omega;C([0,T];\V_n))$ for every finite $l\ge1$.

Finally, for $0<b<b_0<H$, H\"older interpolation \eqref{eq:interpolation}, together with \eqref{eq:Ybound}, gives,
\begin{align*}
\|Y^{m,n}\|_{L^\ell(\Omega;C^b([0,T];\V_n))} &\le C \|Y^{m,n}\|_{L^{2\ell}(\Omega;C([0,T];\V_n))}^{1-\frac{b}{b_0}} \|Y^{m,n}\|_{L^{2\ell}(\Omega;C^{b_0}([0,T];\V_n))}^{\frac{b}{b_0}}\to0,
\end{align*}
which proves \eqref{eq:averaged-noise}.

Hence, \eqref{eq:holder-X}, \eqref{eq:Young}, and \eqref{eq:averaged-noise} yield
\begin{align*}
\|R_4^{m,n}\|_{L^\ell(\Omega;C([0,T];\V))} &\le C \|P_nDF(X)|_{\V_n}\|_ {L^{2\ell}(\Omega;C^a([0,T];\mathcal L(\V_n)))} \left\|Y^{m,n}\right\|_ {L^{2\ell}(\Omega;C^b([0,T];\V_n))}\\ &\to0,\qquad m\to\infty.
\end{align*}

Combining the estimates for $R_i^{m,n}$, $i=1,2,3,4$, with \eqref{eq:I3-error-decomposition}, we obtain
\begin{align}
\label{eq:finite-noise-limit} &\left\|I_3^{m,n} -\frac{T}{2}\int_0^\cdot E_n(\cdot-s)P_nDF(X(s))dW^{H,n}_\Phi(s) \right\|_{L^\ell(\Omega;C([0,T];\V_n))} \notag\\ \le{}& C\left(h^b+h^{a+b-1}+h^{a_0-a}\right) +C\left\|Z^{m,n}-\frac{T}{2}W^{H,n}_\Phi\right\|_ {L^{2\ell}(\Omega;C^b([0,T];\V_n))} \to0.
\end{align}

By \eqref{eq:Ybound} and \eqref{eq:Young}, the stochastic convolution in \eqref{eq:Un-limit} has continuous paths. Since the remaining inhomogeneous terms are also continuous and $DF$ is bounded by \eqref{eq:F} , a standard pathwise Picard iteration yields a unique continuous solution to \eqref{eq:Un-limit}. Subtracting \eqref{eq:Un-limit} from \eqref{eq:tildeUmn} and applying Gronwall's inequality to the supremum norm, \eqref{eq:finite-drift-limits} and \eqref{eq:finite-noise-limit} finish the proof for $X_0\in L^\infty(\Omega;\V^\rho)$.

{\it Step 2. General initial value case.} We now return to the general initial value $X_0$ and the exponent $p\in[2,\infty)$ fixed in the standing assumptions of this section. Set $X_{0,R}:=X_0\mathbf1_{A_R}$ with $A_R:=\{\|X_0\|_\rho\le R\}$ and denote by $\tilde U_R^{m,n}$ and $\tilde{U}_R^{\infty,n}$ the corresponding solutions. On $A_R$, it holds $X_R=X$, while the numerical recursion gives $X_R^m=X^m$. Applying Gronwall's inequality to the corresponding linear equations then yields $\widetilde U_R^{m,n}=\widetilde U^{m,n}$ and $\tilde{U}_R^{\infty,n}=\tilde{U}^{\infty,n}$ on $A_R$ almost surely. Then, for every $\epsilon>0$,
\begin{align*}
&\quad\mathbb P\left(\|\tilde U^{m,n}-\tilde{U}^{\infty,n}\|_{C([0,T];\V_n)}>\epsilon\right)\\ &=\mathbb P\left( \{\|\widetilde U_R^{m,n}-\tilde{U}_R^{\infty,n}\|_{C([0,T];\V_n)}>\epsilon\}\cap A_R \right) +\mathbb P\left( \{\|\widetilde U^{m,n}-\tilde{U}^{\infty,n}\|_{C([0,T];\V_n)}>\epsilon\}\cap A_R^c \right)\\ &\le\mathbb P\left(\|\tilde U_R^{m,n}-\tilde{U}_R^{\infty,n}\|_{C([0,T];\V_n)}>\epsilon\right)+\mathbb P(\|X_0\|_\rho>R) \notag\\ &\le \epsilon^{-p}\|\tilde U_R^{m,n}-\tilde{U}_R^{\infty,n}\|_{L^p(\Omega;C([0,T];\V_n))}^p+R^{-2p}\mathbb E\left[\|X_0\|_\rho^{2p}\right].
\end{align*}
Letting $m\to\infty$ and then $R\to\infty$ completes the proof.
\end{proof}

\subsection{The limiting equation and convergence in path space}
We now pass to the infinite-dimensional limit.

\begin{theorem}\label{tm:asymptotic-error}
Under the standing assumptions of this section, the process $U^m=m(X^m-X)$ converges in probability in $C([0,T];\V)$ to the unique continuous solution of
\begin{align}
\label{eq:limit-U} U(t)={}&\int_0^t E(t-s)DF(X(s))U(s)ds-\frac{T}{2}\int_0^t E(t-s)DF(X(s))AX(s)ds\notag\\ &-\frac{T}{2}\int_0^t E(t-s)DF(X(s))F(X(s))ds -\frac{T}{2}\mathcal I(t),
\end{align}
where
\begin{equation}
\label{eq:I} \mathcal I(t):=\deltaH(K_t) +\alpha_H\int_0^t\int_0^s \sum_{j\ge1}E(t-s)D^2F(X(s)) \bigl(D_u^jX(s),\Phi g_j\bigr)(s-u)^{2H-2}du ds
\end{equation}
with $K_t(s):=\mathbf1_{[0,t)}(s)E(t-s)DF(X(s))\Phi$. For each fixed $t\in[0,T]$, the convergence also holds in $L^\ell(\Omega;\V)$ for every $1\le\ell<p$, where $p$ is the exponent in the standing initial moment assumption.
\end{theorem}

\begin{proof}
We divide the proof into four steps. We first establish the well-posedness of \eqref{eq:limit-U} and the convergence of $\tilde{U}^{\infty,n}$ to $U$ in $L^p(\Omega;\V)$. For bounded initial value, we then establish uniform time increment estimates and use them to pass to path space and combine the approximation results with Proposition \ref{prop:finite-limit}. Finally, we remove the boundedness assumption by localization and prove the stated moment convergence.

{\it Step 1. Well-posedness of \eqref{eq:limit-U} and convergence of $\tilde{U}^{\infty,n}$ to $U$ in $L^p(\Omega;\V)$.} By \eqref{eq:new-index}, $\V^{\rho-2}\hookrightarrow\V^{-\rho-\gamma}$, since $\rho-2\ge-\rho-\gamma$. Thus \eqref{eq:semi1}, \eqref{eq:F}, \eqref{eq:F1}, and Theorem \ref{tm:X} give
\begin{align*}
\left\|\int_0^tE(t-s)DF(X(s))AX(s)ds \right\|_{L^p(\Omega;\V)} &\le C\int_0^t(t-s)^{-\frac{\delta}{2}} \|(1+\|X(s)\|_\rho)\|AX(s)\|_{-\rho-\gamma} \|_{L^p(\Omega)}ds\\ &\le C\int_0^t(t-s)^{-\frac{\delta}{2}} \bigl(1+\|X(s)\|_{L^{2p}(\Omega;\V^\rho)}^2\bigr)ds \le C,\\ \left\|\int_0^tE(t-s)DF(X(s))F(X(s))ds \right\|_{L^p(\Omega;\V)} &\le C\int_0^t \bigl(1+\|X(s)\|_{L^p(\Omega;\V)}\bigr)ds \le C.
\end{align*}

Let $\tau$ be as in Lemma \ref{lm:proj}. Using Assumption \ref{as:F}, Proposition \ref{prop:Malliavin-X}, and \eqref{eq:new-index}, with $\rho+\gamma+\beta-1=(2\rho+\gamma)-2H>0$, we obtain
\begin{align}
\|K_t(s)\|_{L^p(\Omega;\HS(\U,\V))} &\le C(t-s)^{-\frac{\tau}{2}},\label{eq:K-bound}\\ \|D_{\cdot}K_t(s)\|_ {L^p(\Omega;(\HH\otimes\U)\otimes\HS(\U,\V))} &\le C(t-s)^{-\frac{\eta}{2}} \|D_{\cdot}X(s)\|_ {L^p(\Omega;\HH\otimes\HS(\U,\V))} \|(-A)^{-\frac{\sigma}{2}}\Phi\|_{\HS(\U,\V)} \notag\\ &\le C(t-s)^{-\frac{\eta}{2}}.\label{eq:K-bounds}
\end{align}
To apply Lemma \ref{lm:chain} to the whole kernel, fix $t\in(0,T]$ and choose $\eta+1-2H<a<1$. Set
\[
G_0:=L^2\bigl((0,t),(t-s)^{-a}ds;\V\bigr),\qquad G_1:=\HH\otimes\HS(\U,\V).
\]
Theorem \ref{tm:X} and Proposition \ref{prop:Malliavin-X} give $\sup_{0<s<t}\|X(s)\|_{\mathbb D^{1,p}(\V)}<\infty$. We verify that the path $X_t:=X|_{(0,t)}$ belongs to $\mathbb D^{1,p}(G_0)$. Write $w_t(s):=(t-s)^{-a}$ and $\mathfrak H:=\HH\otimes\U$, and use jointly measurable versions of $X(s)$ and $DX(s)$. Since $p/2\ge1$, Minkowski's integral inequality gives
\begin{align*}
\|X_t\|_{L^p(\Omega;G_0)}^2 &=\left\|\int_0^t w_t(s)\|X(s)\|^2ds\right\|_{L^{p/2}(\Omega)}\\* &\le\int_0^t w_t(s)\|X(s)\|_{L^p(\Omega;\V)}^2ds \le C\frac{t^{1-a}}{1-a}<\infty.
\end{align*}
Here $a<1$ ensures that $\int_0^t w_t(s)ds=t^{1-a}/(1-a)$ is finite. Under the Hilbert-space identification
\[
\mathfrak H\otimes G_0 \simeq L^2\bigl((0,t),w_t(s)ds;\mathfrak H\otimes\V\bigr),
\]
define the candidate derivative $Y_t$ by $Y_t(s):=DX(s)$. The same inequality and Proposition \ref{prop:Malliavin-X} yield
\begin{align*}
\|Y_t\|_{L^p(\Omega;\mathfrak H\otimes G_0)}^2 &=\left\|\int_0^t w_t(s)\|DX(s)\|_{\mathfrak H\otimes\V}^2ds \right\|_{L^{p/2}(\Omega)}\\* &\le\int_0^t w_t(s) \|DX(s)\|_{L^p(\Omega;\mathfrak H\otimes\V)}^2ds \le C\frac{t^{1-a}}{1-a}<\infty.
\end{align*}
To identify $Y_t$ as the derivative of $X_t$, fix a deterministic $\varphi\in G_0$. Using $\sup_{0<s<t}\|X(s)\|_{\mathbb D^{1,p}(\V)}<\infty$, which follows from Theorem \ref{tm:X} and Proposition \ref{prop:Malliavin-X}, we obtain by Cauchy--Schwarz
\[
\int_0^t w_t(s) \|\langle X(s),\varphi(s)\rangle_{\V}\|_{\mathbb D^{1,p}(\mathbb R)}ds \le C\left(\int_0^t w_t(s)ds\right)^{1/2}\|\varphi\|_{G_0}<\infty.
\]
Consequently,
\[
Z_\varphi:=\langle X_t,\varphi\rangle_{G_0} =\int_0^t w_t(s)\langle X(s),\varphi(s)\rangle_{\V}ds
\]
is a Bochner integral in $\mathbb D^{1,p}(\mathbb R)$. The operator $D:\mathbb D^{1,p}(\mathbb R)\to L^p(\Omega;\mathfrak H)$ is continuous and linear, so
\[
DZ_\varphi =\int_0^t w_t(s)\langle DX(s),\varphi(s)\rangle_{\V}ds =\langle Y_t,\varphi\rangle_{G_0}.
\]
Here $\langle DX(s),\varphi(s)\rangle_{\V}$ and $\langle Y_t,\varphi\rangle_{G_0}$ take values in $\mathfrak H$, which gives
\[
X_t\in\mathbb D^{1,p}(G_0),\qquad DX_t=Y_t.
\]
Accordingly, $(DX_t)(s)=DX(s)$ for almost every $s\in(0,t)$ under the identification above.

For $x\in G_0$, define $(\Psi_t(x))(s):=E(t-s)DF(x(s))\Phi$ for $0<s<t$, and set the kernel equal to zero outside $(0,t)$. The estimate in \eqref{eq:K-bound} and $\tau<2H$ give $\Psi_t(0)\in G_1$. For every $v\in G_0$, weighted Cauchy--Schwarz with respect to $|s-r|^{2H-2}\,dr\,ds$ gives
\begin{align*}
&\; \int_0^t\int_0^t \frac{\|v(s)\|\,\|v(r)\|\,|s-r|^{2H-2}} {(t-s)^{\eta/2}(t-r)^{\eta/2}}\,dr\,ds \\ \le &\; \left[ \int_0^t\int_0^t (t-s)^{-a}\|v(s)\|^2 (t-r)^{a-\eta}|s-r|^{2H-2}\,dr\,ds \right]^{1/2}\\ &\;\times \left[ \int_0^t\int_0^t (t-r)^{-a}\|v(r)\|^2 (t-s)^{a-\eta}|s-r|^{2H-2}\,dr\,ds \right]^{1/2} \le M_t\|v\|_{G_0}^2,
\end{align*}
where $M_t:=\sup_{0<s<t}\int_0^t(t-r)^{a-\eta}|s-r|^{2H-2}dr$. Denote $q:=(\eta-a)+(2-2H)$. Then $q<1$ due to $a>\eta+1-2H$. Thus,
\begin{align*}
M_t &\le \int_0^t(t-r)^{-q}\,dr + \sup_{0<s<t}\int_0^t|s-r|^{-q}\,dr= \frac{t^{1-q}}{1-q} + \sup_{0<s<t}\frac{s^{1-q}+(t-s)^{1-q}}{1-q}\le \frac{3t^{1-q}}{1-q}.
\end{align*}
Together with \eqref{eq:F2}, the mean-value inequality, and $\|(-A)^{-\sigma/2}\Phi\|_{\HS(\U,\V)}<\infty$, this proves $\|\Psi_t(x)-\Psi_t(y)\|_{G_1}\le C\|x-y\|_{G_0}$. The same noise-coordinate summation and dominated-convergence argument as in Step 2 of the proof of Theorem \ref{tm:rate-k}, using the preceding covariance estimate for $v\in G_0$, proves that $\Psi_t:G_0\to G_1$ is everywhere G\^ateaux differentiable, with bounded linear derivative
\[
(\Psi_t'(x)v)(s)g_j =E(t-s)D^2F(x(s))\bigl(v(s),\Phi g_j\bigr),\qquad 0<s<t.
\]
Lemma \ref{lm:chain} therefore yields
\[
K_t=\Psi_t(X_t)\in\mathbb D^{1,p}(G_1),\qquad DK_t=\Psi_t'(X_t)DX_t.
\]
In coordinates, this gives
\[
D_u^iK_t(s)g_j =E(t-s)D^2F(X(s))\bigl(D_u^iX(s),\Phi g_j\bigr),\qquad u<s<t,
\]
with $D_u^iK_t(s)g_j=0$ for $u>s$.

Since $\tau,\eta<2H$, Minkowski's inequality in $L^{p/2}(\Omega)$, \eqref{eq:HH}, and \eqref{eq:K-bound}--\eqref{eq:K-bounds} yield
\begin{align}
\label{eq:K-domain} \|K_t\|_{L^p(\Omega;\HH\otimes\HS(\U,\V))}^2 +\|DK_t\|_ {L^p(\Omega;(\HH\otimes\U)^{\otimes2}\otimes\V)}^2 &\le C\sum_{c\in\{\tau,\eta\}} \int_0^t\int_0^t (t-r)^{-\frac{c}{2}} (t-v)^{-\frac{c}{2}} |r-v|^{2H-2}drdv \notag\\ &\le C\sum_{c\in\{\tau,\eta\}}t^{2H-c}\le C.
\end{align}
Then \eqref{eq:div-Lp} with $r=p$, together with \eqref{eq:K-domain}, gives $\|\deltaH(K_t)\|_{L^p(\Omega;\V)}\le C$.

For the correction term in \eqref{eq:I}, \eqref{eq:F2} and \eqref{eq:DXpoint} give
\begin{align}
\label{eq:limit-trace-bound} &\quad\left\|\alpha_H\int_0^t\int_0^s \sum_{j\ge1}E(t-s)D^2F(X(s)) \bigl(D_u^jX(s),\Phi g_j\bigr) (s-u)^{2H-2}duds \right\|_{L^p(\Omega;\V)} \notag\\ &\le C\int_0^t(t-s)^{-\frac{\eta}{2}} \int_0^s (s-u)^{\frac{\beta-1}{2}+2H-2}duds =\frac{C}{\nu} \int_0^t(t-s)^{-\frac{\eta}{2}}s^\nu ds \le C,
\end{align}
where $\sigma=1-\beta$ and $\nu>0$ follow from \eqref{eq:new-index} and \eqref{eq:nu}.

The preceding estimates show that all inhomogeneous terms in \eqref{eq:limit-U} are well-defined and uniformly bounded in $L^p(\Omega;\V)$. Since $DF$ is bounded, a standard fixed-point argument yields a unique solution $U$ in the class of strongly measurable processes uniformly bounded in time in $L^p(\Omega;\V)$ satisfying
\[
\sup_{t\in[0,T]}\|U(t)\|_{L^p(\Omega;\V)}\le C.
\]

We next establish convergence of $\tilde{U}^{\infty,n}$ to $U$ in $L^p(\Omega;\V)$ uniformly in time. Set
\begin{align*}
\mathcal I^n(t) :={}& \int_0^tE_n(t-s)P_nDF(X(s))dW^{H,n}_\Phi(s)\\ ={}&\deltaH(K_t^n) +\alpha_H\int_0^t\int_0^s \sum_{j\ge1}E(t-s)P_nD^2F(X(s)) \bigl(D_u^jX(s),P_n\Phi g_j\bigr) (s-u)^{2H-2}duds,
\end{align*}
where $K_t^n(s):=\mathbf1_{[0,t)}(s) E(t-s)P_nDF(X(s))P_n\Phi.$ The last equality follows as in the divergence decomposition above, using the relation between the Young and divergence integrals (cf. \cite[Proposition 2.1]{N06fBm}).

To verify the required absolute integrability of \cite[Proposition 2.1]{N06fBm}, set $\mathcal Q:=\HS(\U,\HS(\U,\V))$. By \eqref{eq:F2}, \eqref{eq:DXpoint}, and $\sigma=1-\beta$,
\[
\|D_uK_t^n(s)\|_{L^2(\Omega;\mathcal Q)} \le C(t-s)^{-\frac{\eta}{2}}(s-u)^{\frac{\beta-1}{2}}, \qquad 0\le u<s<t,
\]
while $D_uK_t^n(s)=0$ for $u>s$. Cauchy--Schwarz in $\HH$, \eqref{eq:double-int}, and $\eta<2H$ give
\begin{align*}
&\E\int_{[0,t]^4} \|D_uK_t^n(s)\|_{\mathcal Q} \|D_vK_t^n(r)\|_{\mathcal Q} |u-v|^{2H-2}|s-r|^{2H-2}\,du\,dv\,ds\,dr\\ &\qquad\le C\int_0^t\int_0^t (t-s)^{-\frac{\eta}{2}}(t-r)^{-\frac{\eta}{2}} |s-r|^{2H-2}\,ds\,dr<\infty.
\end{align*}
For $K_t^n$ itself, since $p\ge2$, Cauchy--Schwarz in $\Omega$ and \eqref{eq:K-bound} with the projections inserted give
\begin{align*}
&\alpha_H\E\int_0^t\int_0^t \|K_t^n(s)\|_{\HS(\U,\V)}\|K_t^n(r)\|_{\HS(\U,\V)} |s-r|^{2H-2}\,dr\,ds\\ &\quad\le\alpha_H\int_0^t\int_0^t \|K_t^n(s)\|_{L^2(\Omega;\HS(\U,\V))} \|K_t^n(r)\|_{L^2(\Omega;\HS(\U,\V))} |s-r|^{2H-2}\,dr\,ds\\ &\quad\le C\int_0^t\int_0^t (t-s)^{-\frac{\tau}{2}}(t-r)^{-\frac{\tau}{2}} |s-r|^{2H-2}\,dr\,ds \le Ct^{2H-\tau}<\infty,
\end{align*}
since $\tau<2H$. Together with the preceding fourfold integral estimate, this verifies, coordinatewise, both second-moment terms in the first integrability condition of \cite[Proposition 2.1]{N06fBm}. For the absolute-integrability condition \cite[(2.12)]{N06fBm} concerning the correction term in $\mathcal I^n(t)$, note that
\[
D_u^jK_t^n(s)g_j =E(t-s)P_nD^2F(X(s)) \bigl(D_u^jX(s),P_n\Phi g_j\bigr),\qquad u<s<t.
\]
The argument for \eqref{eq:limit-trace-bound}, with the projections inserted and the sum of the norms inside the integrals, yields
\begin{align*}
&\E\int_0^t\int_0^s\sum_{j\ge1} \|D_u^jK_t^n(s)g_j\|_{\V} (s-u)^{2H-2}\,du\,ds\\ &\quad\le C\int_0^t(t-s)^{-\frac{\eta}{2}} \int_0^s(s-u)^{\frac{\beta-1}{2}+2H-2}\,du\,ds =\frac{C}{\nu}\int_0^t(t-s)^{-\frac{\eta}{2}}s^\nu\,ds<\infty,
\end{align*}
because $\nu>0$ and $\eta/2<1$. The nonnegative random integral inside this expectation is therefore finite almost surely. This verifies the cited condition (2.12) and the preceding absolute integrability bounds.

We decompose the kernel difference as
\begin{align*}
K_t(s)-K_t^n(s) =\mathbf1_{[0,t)}(s) E(t-s)(I-P_n)DF(X(s))\Phi +\mathbf1_{[0,t)}(s) E(t-s)P_nDF(X(s))(I-P_n)\Phi.
\end{align*}
Estimating the divergence and correction terms as in the proof of Lemma \ref{lm:proj}, we obtain
\begin{align}
\label{eq:limit-noise-tail} \sup_{t\in[0,T]} \|\mathcal I^n(t)-\mathcal I(t)\|_{L^p(\Omega;\V)} &\le C\lambda_{n+1}^{-\vartheta} \sup_{t\in[0,T]} \sum_{c\in\{\tau,\eta\}} \left( \int_0^t\int_0^t (t-r)^{-\frac{c}{2}-\vartheta} (t-v)^{-\frac{c}{2}-\vartheta} |r-v|^{2H-2}drdv \right)^{\frac12} \notag\\ &\quad+ C\lambda_{n+1}^{-\vartheta} \sup_{t\in[0,T]} \int_0^t (t-s)^{-\frac{\eta}{2}-\vartheta}s^\nu ds +C\varepsilon_n \notag\\ &\le C\bigl(\lambda_{n+1}^{-\vartheta} +\varepsilon_n\bigr)\to0,
\end{align}
where $\vartheta$ is as in Lemma \ref{lm:proj}, and we used \eqref{eq:spec-tail}, \eqref{eq:noise-tail}, and \eqref{eq:K-bounds}.

The convergence of the deterministic forcings follows by applying the estimates \eqref{eq:S1}--\eqref{eq:S2} with $X$ in place of $X^m$ and $\frac{T}{2}$ in place of $m(s-\kappa_m(s))$, where the $L^{2p}$ spatial tail estimate in \eqref{eq:an} also holds for $X$. The projection error in the linear term is estimated as in \eqref{eq:S0}. Gronwall's inequality then yields
\begin{equation}
\label{eq:Un-pointwise} \sup_{t\in[0,T]} \|\tilde{U}^{\infty,n}(t)-U(t)\|_{L^p(\Omega;\V)}\to0.
\end{equation}

{\it Step 2. Uniform time increment estimates.} Throughout Steps 2 and 3, assume that $X_0\in L^\infty(\Omega;\V^\rho)$, and let $\ell\in[2,\infty)$ be arbitrary.

To pass from convergence uniformly in time in $L^\ell(\Omega;\V)$ to convergence in $L^\ell(\Omega;C([0,T];\V))$ by the argument for \eqref{eq:mesh-upgrade}, we need to establish uniform H\"older moment bounds for $\tilde U^m$, $\tilde U^{m,n}$, $\tilde{U}^{\infty,n}$, and $U$. By Kolmogorov's criterion, it suffices to prove, for some $\theta>0$ independent of $\ell$, that
\begin{align}
\label{eq:Y-increment} \sup_{m,n\ge1} \sup_{Y\in\{\tilde U^m,\tilde U^{m,n},\tilde{U}^{\infty,n},U\}} \|Y(b)-Y(a)\|_{L^\ell(\Omega;\V)} \le C(b-a)^\theta, \qquad 0\le a<b\le T.
\end{align}
For \eqref{eq:tildeUm} and \eqref{eq:tildeUmn}, we estimate the increments of $B_i^m$ and $I_i^{m,n}$, $i=1,2,3$. We also need to estimate the increments of $R_1^m$ and $R_2^m$ to control $U^m-\tilde U^m$ through \eqref{eq:Um} and \eqref{eq:tildeUm}. Specifically, fixing
\begin{equation}
\label{eq:path-exponent} 0<\theta< \min\left\{ 1-\frac{\delta}{2}, H-\frac{\tau}{2}, H-\frac{\eta}{2} \right\},
\end{equation}
we next verify that
\begin{equation}
\label{eq:forcing-increments}
\begin{aligned}
&\sum_{i=1}^3 \left( \|B_i^m(b)-B_i^m(a)\|_{L^\ell(\Omega;\V)} +\|I_i^{m,n}(b)-I_i^{m,n}(a)\|_{L^\ell(\Omega;\V)} \right) \le C(b-a)^\theta,\\ &\|R_1^m(b)-R_1^m(a)\|_{L^\ell(\Omega;\V)} \le Ch(b-a)^\theta, \qquad \|R_2^m(b)-R_2^m(a)\|_{L^\ell(\Omega;\V)} \le Ch^\nu(b-a)^\theta,
\end{aligned}
\end{equation}
where $C$ may depend on $\ell$ but is independent of $m,n,a,b$. By \eqref{eq:semi1} and \eqref{eq:semi2},
\begin{equation}
\label{eq:increment-semigroup} \|(E(b-a)-I)E(a-s)v\| \le C(b-a)^\theta(a-s)^{-c-\theta}\|v\|_{-2c}, \qquad s<a.
\end{equation}
For $Y(t)=\int_0^tE(t-s)f(s)ds$, where $f$ satisfies $\sup_{s\in[0,T]}\|f(s)\|_{L^\ell(\Omega;\V^{-2c})}\le M_\ell$ for some $0\le c<1-\theta$, \eqref{eq:increment-semigroup} gives
\begin{align}
\label{eq:det-convolution-increment} \|Y(b)-Y(a)\|_{L^\ell(\Omega;\V)} &\le \|(E(b-a)-I)Y(a)\|_{L^\ell(\Omega;\V)} +\left\|\int_a^bE(b-s)f(s)ds \right\|_{L^\ell(\Omega;\V)} \notag\\ &\le CM_\ell\left( (b-a)^\theta\int_0^a(a-s)^{-c-\theta}ds +\int_a^b(b-s)^{-c}ds \right) \notag\\ &\le CM_\ell(b-a)^\theta.
\end{align}

For $R_1^m$, $R_2^m$, $B_i^{m}$, and $I_i^m$, $i=1,2,$ they are in the form of $Y$ above with corresponding pairs $(c,M_\ell)$ given in the following chart:
\[
\begin{array}{c|cccc}
Y &B_1^m,\ I_1^{m,n} &B_2^m,\ I_2^{m,n} &R_1^m &R_2^m\\ \hline c &\frac{\delta}{2} &0 &\frac{\eta}{2} &\frac{\eta}{2}\\ M_\ell &C &C &Ch &Ch^\nu
\end{array}
\]
according to the arguments in the proof of Lemma \ref{lm:RB}, Theorem \ref{tm:rate} at moment order $2\ell$, and \eqref{eq:new-index} for $R_1^m$. Thus \eqref{eq:det-convolution-increment} gives \eqref{eq:forcing-increments} for these terms.

For $B_3^m=\mathcal J_m(I,I;\cdot)$ and $I_3^{m,n}=\mathcal J_m(P_n,P_n;\cdot)$, let $Q\in\{I,P_n\}$, and define
\[
\mathcal U_{a,b}(r) := m\int_{\max\{a,r\}}^{(\kappa_m(r)+h)\wedge b} E(b-s)QDF(X^m(\kappa_m(s)))E(s-r)Q\Phi ds
\]
for $\kappa_m(a)<r<b$ and set $\mathcal U_{a,b}(r)=0$ otherwise, and
\begin{align*}
\mathcal C_{a,b} :={}&m\alpha_H\int_a^b\int_0^{\kappa_m(s)} \int_{\kappa_m(s)}^s \sum_{j\ge1}E(b-s)QD^2F(X^m(\kappa_m(s)))\\ &\qquad\qquad\times \bigl(D_u^jX^m(\kappa_m(s)),E(s-r)Q\Phi g_j\bigr) (r-u)^{2H-2}drduds.
\end{align*}
Then the semigroup property and decomposition \eqref{eq:Jm} give
\begin{align}
\label{eq:JmQ} \mathcal J_m(Q,Q;b)-\mathcal J_m(Q,Q;a) = (E(b-a)-I)\mathcal J_m(Q,Q;a) +\deltaH(\mathcal U_{a,b})+\mathcal C_{a,b}.
\end{align}

The first term on the right-hand side, by \eqref{eq:semi2}, \eqref{eq:Jm}, and \eqref{eq:div-Lp} with $r=\ell$, satisfies
\begin{align}
\label{eq:B3I3} &\quad \|(E(b-a)-I)B_3^m(a)\|_{L^\ell(\Omega;\V)} +\|(E(b-a)-I)I_3^{m,n}(a)\|_{L^\ell(\Omega;\V)}\notag\\ &\le C(b-a)^\theta\left( \|(-A)^\theta B_3^m(a)\|_{L^\ell(\Omega;\V)} +\|(-A)^\theta I_3^{m,n}(a)\|_{L^\ell(\Omega;\V)} \right) \le C(b-a)^\theta.
\end{align}
The last inequality follows by applying the estimates for \eqref{eq:Jm} with the additional factor $(-A)^\theta$, since \eqref{eq:path-exponent} gives $\frac{\tau}{2}+\theta<H$ and $\frac{\eta}{2}+\theta<H$.

For $c\in\{\frac{\tau}{2},\frac{\eta}{2}\}$, when $b-a\ge h$, \eqref{eq:HH}, \eqref{eq:double-int}, \eqref{eq:DXm-sup}, and the fact $b-\kappa_m(a)\le2(b-a)$ give
\begin{align*}
&\quad\|\mathcal U_{a,b}\|_ {L^\ell(\Omega;\HH\otimes\HS(\U,\V))} +\|D\mathcal U_{a,b}\|_ {L^\ell(\Omega;(\HH\otimes\U)^{\otimes2}\otimes\V)}\\ &\le C \sum_{c\in\{\frac{\tau}{2},\frac{\eta}{2}\}} \left( \int_{\kappa_m(a)}^b\int_{\kappa_m(a)}^b (b-r)^{-c}(b-v)^{-c} |r-v|^{2H-2}drdv \right)^{\frac12}\\ &\le C \sum_{c\in\{\frac{\tau}{2},\frac{\eta}{2}\}} (b-\kappa_m(a))^{H-c} \le C\left( (b-a)^{H-\frac{\tau}{2}} +(b-a)^{H-\frac{\eta}{2}} \right),
\end{align*}
where we used the following estimate based on Lemma \ref{lm:int}:
\[
m\int_{\max\{a,r\}}^{(\kappa_m(r)+h)\wedge b} (b-s)^{-c}ds \le Cmh(b-r)^{-c}, \qquad \kappa_m(a)<r<b.
\]
When $b-a<h$, \eqref{eq:HH}, \eqref{eq:F}--\eqref{eq:F2}, and \eqref{eq:DXm-sup} give, for $a<s<b$,
\begin{align*}
&\quad \left\| \mathbf1_{[\kappa_m(s),s]}(\cdot) E(b-s)QDF(X^m(\kappa_m(s)))E(s-\cdot)Q\Phi \right\|_{\mathbb D^{1,\ell}(\HH\otimes\HS(\U,\V))}\\ &\le C(s-\kappa_m(s))^H \left( (b-s)^{-\frac{\tau}{2}} +(b-s)^{-\frac{\eta}{2}} \right) \le Ch^H \left( (b-s)^{-\frac{\tau}{2}} +(b-s)^{-\frac{\eta}{2}} \right).
\end{align*}
Minkowski's inequality, together with the facts $b-a<h$ and $H<1$, therefore yields
\begin{align*}
&\quad\|\mathcal U_{a,b}\|_ {L^\ell(\Omega;\HH\otimes\HS(\U,\V))} +\|D\mathcal U_{a,b}\|_ {L^\ell(\Omega;(\HH\otimes\U)^{\otimes2}\otimes\V)}\\ &\le Cmh^H \sum_{c\in\{\frac{\tau}{2},\frac{\eta}{2}\}} \int_a^b(b-s)^{-c}ds \le Ch^{H-1} \sum_{c\in\{\frac{\tau}{2},\frac{\eta}{2}\}} (b-a)^{1-c} \le C\left( (b-a)^{H-\frac{\tau}{2}} +(b-a)^{H-\frac{\eta}{2}} \right).
\end{align*}

By \eqref{eq:F2} and \eqref{eq:DXm-norm}, as in the proof of Lemma \ref{lm:proj},
\[
\|\mathcal C_{a,b}\|_{L^\ell(\Omega;\V)} \le Cmh\int_a^b(b-s)^{-\frac{\eta}{2}}ds \le C(b-a)^{1-\frac{\eta}{2}}.
\]
Hence, \eqref{eq:div-Lp} with $r=\ell$ and \eqref{eq:path-exponent} imply
\begin{align*}
\|\deltaH(\mathcal U_{a,b})+\mathcal C_{a,b} \|_{L^\ell(\Omega;\V)} \le C\left( (b-a)^{H-\frac{\tau}{2}} +(b-a)^{H-\frac{\eta}{2}} +(b-a)^{1-\frac{\eta}{2}} \right) \le C(b-a)^\theta.
\end{align*}
Together with \eqref{eq:JmQ} with $Q=I$ and $Q=P_n$, and \eqref{eq:B3I3}, it proves \eqref{eq:forcing-increments} for $B_3^m$ and $I_3^{m,n}$.

Applying \eqref{eq:K-bounds}, \eqref{eq:limit-trace-bound}, and \eqref{eq:increment-semigroup} to \eqref{eq:I} and its projected version gives
\begin{align}
\label{eq:limit-I-increment} \sup_{n\ge1} \|\mathcal I^n(b)-\mathcal I^n(a)\|_{L^\ell(\Omega;\V)} +\|\mathcal I(b)-\mathcal I(a)\|_{L^\ell(\Omega;\V)} \le C(b-a)^\theta.
\end{align}

By \eqref{eq:F}--\eqref{eq:F1}, \eqref{eq:new-index}, and Theorem \ref{tm:X}, for $Q\in\{I,P_n\}$,
\begin{align*}
\sup_{s\in[0,T]} \|QDF(X(s))AQX(s)\|_{L^\ell(\Omega;\V^{-\delta})} &\le C\left( 1+\sup_{s\in[0,T]} \|X(s)\|_{L^{2\ell}(\Omega;\V^\rho)}^2 \right)\le C,\\ \sup_{s\in[0,T]} \|QDF(X(s))QF(X(s))\|_{L^\ell(\Omega;\V)} &\le C\left( 1+\sup_{s\in[0,T]} \|X(s)\|_{L^\ell(\Omega;\V)} \right)\le C.
\end{align*}
Moreover, \eqref{eq:F}, \eqref{eq:tildeUm-bound}, Lemma \ref{lm:tildeUmn}, and Step 1 give uniformly in $m,n$ that
\[
\sup_{s\in[0,T]} \|QDF(X(s))Y(s)\|_{L^\ell(\Omega;\V)} \le C\sup_{s\in[0,T]} \|Y(s)\|_{L^\ell(\Omega;\V)} \le C, \qquad Y\in\{\tilde U^m,\tilde U^{m,n},\tilde{U}^{\infty,n},U\}.
\]
Then applying \eqref{eq:det-convolution-increment} with $c=\frac{\delta}{2}$ for $f(s)=QDF(X(s))AQX(s)$, and with $c=0$ for $f(s)=QDF(X(s))QF(X(s))$ and $f(s)=QDF(X(s))Y(s)$, we obtain
\begin{align*}
\|\tilde U^m(b)-\tilde U^m(a)\|_{L^\ell(\Omega;\V)} &\le C(b-a)^\theta +\sum_{i=1}^3 \|B_i^m(b)-B_i^m(a)\|_{L^\ell(\Omega;\V)} \le C(b-a)^\theta,\\ \|\tilde U^{m,n}(b)-\tilde U^{m,n}(a)\|_{L^\ell(\Omega;\V)} &\le C(b-a)^\theta +\sum_{i=1}^3 \|I_i^{m,n}(b)-I_i^{m,n}(a)\|_{L^\ell(\Omega;\V)} \le C(b-a)^\theta,\\ \|\tilde{U}^{\infty,n}(b)-\tilde{U}^{\infty,n}(a)\|_{L^\ell(\Omega;\V)} &\le C(b-a)^\theta +\frac{T}{2} \|\mathcal I^n(b)-\mathcal I^n(a)\|_{L^\ell(\Omega;\V)} \le C(b-a)^\theta,\\ \|U(b)-U(a)\|_{L^\ell(\Omega;\V)} &\le C(b-a)^\theta +\frac{T}{2} \|\mathcal I(b)-\mathcal I(a)\|_{L^\ell(\Omega;\V)} \le C(b-a)^\theta,
\end{align*}
where the last inequalities follow from \eqref{eq:forcing-increments} and \eqref{eq:limit-I-increment}. This proves \eqref{eq:Y-increment}.

{\it Step 3. Convergence in path space for bounded initial value.} Fix $\theta'\in(0,\theta)$. Since \eqref{eq:Y-increment}, \eqref{eq:forcing-increments}, and \eqref{eq:limit-I-increment} hold at every finite moment order, Kolmogorov's criterion, applied at an order greater than $\max\{\ell,(\theta-\theta')^{-1}\}$, gives continuous modifications satisfying
\begin{align}
\label{eq:path-moments} &\sup_{m,n\ge1}\left( \|\tilde U^m\|_{L^\ell(\Omega;C^{\theta'})} +\|\tilde U^{m,n}\|_{L^\ell(\Omega;C^{\theta'})} +\|\tilde{U}^{\infty,n}\|_{L^\ell(\Omega;C^{\theta'})} \right) +\|U\|_{L^\ell(\Omega;C^{\theta'})} \le C, \notag\\ &\|R_1^m\|_{L^\ell(\Omega;C^{\theta'})} +\|R_2^m\|_{L^\ell(\Omega;C^{\theta'})} \le C(h+h^\nu),
\end{align}
where $C^{\theta'}=C^{\theta'}([0,T];\V)$, and we used that all processes vanish at zero. The same argument gives
\[
\sup_{n\ge1} \|\mathcal I^n\|_{L^\ell(\Omega;C^{\theta'})} +\|\mathcal I\|_{L^\ell(\Omega;C^{\theta'})} \le C.
\]
We henceforth use these continuous modifications without changing the notation.

Repeating the argument for \eqref{eq:mesh-upgrade} with $\tilde U^{m,n}-\tilde U^m$, and setting $r_j=jT/N$, we obtain
\begin{align*}
\sup_{m\ge1} \|\tilde U^{m,n}-\tilde U^m\|_ {L^\ell(\Omega;C([0,T];\V))} \le(N+1)^{\frac1\ell} \sup_{m\ge1}\max_{0\le j\le N} \|\tilde U^{m,n}(r_j)-\tilde U^m(r_j)\|_ {L^\ell(\Omega;\V)} +C\left(\frac{T}{N}\right)^{\theta'}.
\end{align*}
First let $n\to\infty$ by Lemma \ref{lm:tildeUmn-tildeUm}, and then $N\to\infty$. Applying the same argument to $\tilde{U}^{\infty,n}-U$ and $\mathcal I^n-\mathcal I$, using \eqref{eq:Un-pointwise} and \eqref{eq:limit-noise-tail}, respectively, yields
\begin{equation}
\label{eq:path-spatial}
\begin{aligned}
&\lim_{n\to\infty}\sup_{m\ge1} \|\tilde U^{m,n}-\tilde U^m\|_ {L^\ell(\Omega;C([0,T];\V))}=0,\\ &\lim_{n\to\infty} \|\tilde{U}^{\infty,n}-U\|_{L^\ell(\Omega;C([0,T];\V))}=0, \qquad \lim_{n\to\infty} \|\mathcal I^n-\mathcal I\|_ {L^\ell(\Omega;C([0,T];\V))}=0.
\end{aligned}
\end{equation}

By \eqref{eq:F}, pathwise Gronwall's inequality, and \eqref{eq:path-moments},
\begin{equation}
\label{eq:path-remainders} \|U^m-\tilde U^m\|_{L^\ell(\Omega;C([0,T];\V))} \le C(h+h^\nu)\to0,\qquad m\to\infty.
\end{equation}

Combining \eqref{eq:path-spatial}, \eqref{eq:path-remainders}, and Proposition \ref{prop:finite-limit} in
\begin{align*}
\|U^m-U\|_{L^\ell(\Omega;C([0,T];\V))} \le{}& \|U^m-\tilde U^m\|_{L^\ell(\Omega;C([0,T];\V))} +\|\tilde U^m-\tilde U^{m,n}\|_ {L^\ell(\Omega;C([0,T];\V))}\\ &+\|\tilde U^{m,n}-\tilde{U}^{\infty,n}\|_ {L^\ell(\Omega;C([0,T];\V))} +\|\tilde{U}^{\infty,n}-U\|_{L^\ell(\Omega;C([0,T];\V))},
\end{align*}
and first letting $m\to\infty$, then $n\to\infty$, we obtain
\[
\lim_{m\to\infty}\|U^m-U\|_{L^\ell(\Omega;C([0,T];\V))}=0
\]
for bounded $X_0$ and every finite $\ell\in[2,\infty)$.

{\it Step 4. General initial value and moment convergence.} Return to the general initial value and the exponent $p\in[2,\infty)$ fixed in the standing assumptions.

As in the proof of Proposition \ref{prop:finite-limit}, set $A_R:=\{\|X_0\|_\rho\le R\}$ and $X_{0,R}:=X_0\mathbf1_{A_R}$, and denote the corresponding solutions by a subscript $R$. Then it holds almost surely on $A_R$ that $Y_R=Y$ for $Y\in\{X,X^m,\mathcal I^n,\tilde{U}^{\infty,n}\}$.

For $S\ge R$, it follows $\mathbf1_{A_R}(\tilde{U}_S^{\infty,n}-\tilde{U}_R^{\infty,n})=0$ and $\mathbf1_{A_R}(\mathcal I_S^n-\mathcal I_R^n)=0$. Passing to the limits in \eqref{eq:path-spatial} with $\ell=p$ gives, in $C([0,T];\V)$, that $\mathbf1_{A_R}(U_S-U_R)=0$ and $\mathbf1_{A_R}(\mathcal I_S-\mathcal I_R)=0$ almost surely. Since $\mathbb P(\cup_{R\in\mathbb N}A_R)=1$, the consistency above allows us to piece together $U_R$ and $\mathcal I_R$, respectively, into continuous processes for the general initial value.

For each fixed $t$, \eqref{eq:Un-pointwise} and \eqref{eq:limit-noise-tail} give, in $L^p(\Omega;\V)$,
\begin{align*}
\mathbf1_{A_R}U_R(t) &=\lim_{n\to\infty}\mathbf1_{A_R}\tilde{U}_R^{\infty,n}(t) =\lim_{n\to\infty}\mathbf1_{A_R}\tilde{U}^{\infty,n}(t) =\mathbf1_{A_R}U(t),\\ \mathbf1_{A_R}\mathcal I_R(t) &=\lim_{n\to\infty}\mathbf1_{A_R}\mathcal I_R^n(t) =\lim_{n\to\infty}\mathbf1_{A_R}\mathcal I^n(t) =\mathbf1_{A_R}\mathcal I(t).
\end{align*}
Hence these processes are continuous versions of $U$ and $\mathcal I$ from Step 1 and satisfy \eqref{eq:limit-U}--\eqref{eq:I}. We retain the same notation, so that $U=U_R$ on $A_R$.

Since $U^m=U_R^m$ on $A_R$, for every $\epsilon>0$,
\begin{align*}
\mathbb P\bigl( \|U^m-U\|_{C([0,T];\V)}>\epsilon \bigr) &\le \mathbb P(A_R^c) +\mathbb P\bigl( \|U_R^m-U_R\|_{C([0,T];\V)}>\epsilon \bigr)\\ &\le R^{-2p}\mathbb E\left[\|X_0\|_\rho^{2p}\right] +\epsilon^{-p} \|U_R^m-U_R\|_{L^p(\Omega;C([0,T];\V))}^p.
\end{align*}
Letting $m\to\infty$ by Step 3 and then $R\to\infty$ proves convergence in probability.

For any two continuous solutions $U$ and $\widehat U$ of \eqref{eq:limit-U}, \eqref{eq:F} gives
\[
\|U(t)-\widehat U(t)\| \le C\int_0^t\|U(s)-\widehat U(s)\|ds.
\]
Gronwall's inequality gives uniqueness of the solution to \eqref{eq:limit-U}.

Finally, Theorem \ref{tm:rate} yields $\sup_{m\ge1}\sup_{t\in[0,T]} \mathbb E\left[\|U^m(t)\|^p\right]\le C$. The convergence in probability above indicates that there exists a subsequence $\{U^{m_j}\}_{j\ge1}\subset\{U^m\}_{m\ge1}$ such that $U^{m_j}\to U$ almost surely in $C([0,T];\V)$. Then Fatou's lemma gives
\[
\mathbb E\left[\|U(t)\|^p\right] \le\liminf_{j\to\infty}\mathbb E\left[\|U^{m_j}(t)\|^p\right] \le C.
\]
Then for $1\le\ell<p$, $K>\epsilon>0$, and fixed $t$,
\begin{align*}
&\quad\mathbb E\left[\|U^m(t)-U(t)\|^\ell\right]\\ &=\mathbb E\left[\|U^m(t)-U(t)\|^\ell\left({\bf 1}_{\{\|U^m(t)-U(t)\|\le\epsilon\}}+{\bf 1}_{\{\epsilon<\|U^m(t)-U(t)\|\le K\}}+{\bf 1}_{\{\|U^m(t)-U(t)\|>K\}}\right)\right]\\ &\le \epsilon^\ell +K^\ell\mathbb P\bigl( \|U^m(t)-U(t)\|>\epsilon \bigr) +K^{\ell-p}\mathbb E\left[\|U^m(t)-U(t)\|^p\right]\\ &\le \epsilon^\ell +K^\ell\mathbb P\bigl( \|U^m(t)-U(t)\|>\epsilon \bigr) +CK^{\ell-p}.
\end{align*}
Letting $m\to\infty$, $K\to\infty$, and $\epsilon\to0$ completes the proof.
\end{proof}

\section*{Declaration of generative AI and AI-assisted
technologies in the manuscript preparation process}

The authors used ChatGPT (OpenAI) to assist with literature
searches, language editing, and the construction of examples
satisfying the assumptions on the coefficients. The authors
reviewed and verified the resulting content, including the
references and examples, and take full responsibility for
the final manuscript.


\begin{thebibliography}{}

\bibitem{CGW24}
M. Cai, S. Gan, and X. Wang, Weak approximations of stochastic partial differential equations with fractional noise, J. Comput. Math. 42 (2024), no.~3, 735--754. \href{https://doi.org/10.4208/jcm.2203-m2021-0194}{doi:10.4208/jcm.2203-m2021-0194}.

\bibitem{CSY24}
H. Chen, Z. Shen, and Q. Yu, Exponential Euler method for stiff SDEs driven by fractional Brownian motion, arXiv:2407.03546. \href{https://doi.org/10.48550/arXiv.2407.03546}{doi:10.48550/arXiv.2407.03546}.

\bibitem{FH20}
P. K. Friz and M. Hairer, {\it A Course on Rough Paths: With an Introduction to Regularity Structures}, 2nd ed., Universitext, Springer, Cham, 2020.

\bibitem{HH20}
J. Hong and C. Huang, Super-convergence analysis on exponential integrator for stochastic heat equation driven by additive fractional Brownian motion, arXiv:2007.02223. \href{https://doi.org/10.48550/arXiv.2007.02223}{doi:10.48550/arXiv.2007.02223}.

\bibitem{HJW25}
J. Hong, D. Jin, and X. Wang, Asymptotic error distribution of numerical methods for parabolic SPDEs with multiplicative noise, arXiv:2511.05251. \href{https://doi.org/10.48550/arXiv.2511.05251}{doi:10.48550/arXiv.2511.05251}.

\bibitem{HJWY24}
J. Hong, D. Jin, X. Wang, and G. Yang, Asymptotic error distribution of accelerated exponential Euler method for parabolic SPDEs, Ann. Appl. Probab., to appear.

\bibitem{HLN16}
Y. Hu, Y. Liu, and D. Nualart, Rate of convergence and asymptotic error distribution of Euler approximation schemes for fractional diffusions, Ann. Appl. Probab. 26 (2016), no.~2, 1147--1207. \href{https://doi.org/10.1214/15-AAP1114}{doi:10.1214/15-AAP1114}.

\bibitem{JP98}
J. Jacod and P. Protter, Asymptotic error distributions for the Euler method for stochastic differential equations, Ann. Probab. 26 (1998), no.~1, 267--307. \href{https://doi.org/10.1214/aop/1022855419}{doi:10.1214/aop/1022855419}.

\bibitem{JK09}
A. Jentzen and P. E. Kloeden, Overcoming the order barrier in the numerical approximation of stochastic partial differential equations with additive space-time noise, Proc. R. Soc. Lond. Ser. A Math. Phys. Eng. Sci. 465 (2009), no.~2102, 649--667. \href{https://doi.org/10.1098/rspa.2008.0325}{doi:10.1098/rspa.2008.0325}.

\bibitem{M10}
J. Maas, Malliavin calculus and decoupling inequalities in Banach spaces, J. Math. Anal. Appl. 363 (2010), no.~2, 383--398. doi:10.1016/j.jmaa.2009.08.041.

\bibitem{MPR17}
T. Mastrolia, D. Possama\"i, and A. R\'eveillac, On the Malliavin differentiability of BSDEs, Ann. Inst. Henri Poincar\'e Probab. Stat. 53 (2017), no.~1, 464--492. \href{https://doi.org/10.1214/15-AIHP723}{doi:10.1214/15-AIHP723}.

\bibitem{MN19}
R. Musina and A. I. Nazarov, A note on truncations in fractional Sobolev spaces, Bull. Math. Sci. 9 (2019), no.~1, 1950001. \href{https://doi.org/10.1142/S1664360719500012}{doi:10.1142/S1664360719500012}.

\bibitem{NN07}
A. Neuenkirch and I. Nourdin, Exact rate of convergence of some approximation schemes associated to SDEs driven by a fractional Brownian motion, J. Theoret. Probab. 20 (2007), no.~4, 871--899. \href{https://doi.org/10.1007/s10959-007-0083-0}{doi:10.1007/s10959-007-0083-0}.

\bibitem{N06fBm}
D. Nualart, Stochastic calculus with respect to fractional Brownian motion, Ann. Fac. Sci. Toulouse Math. (6) 15 (2006), no. 1, 63--78. \href{https://doi.org/10.5802/afst.1113}{doi:10.5802/afst.1113}.

\bibitem{N06}
D. Nualart, {\it The Malliavin Calculus and Related Topics}, 2nd ed., Probab. Appl. (N. Y.) Springer-Verlag, Berlin, 2006. \href{https://doi.org/10.1007/3-540-28329-3}{doi:10.1007/3-540-28329-3}.

\bibitem{PV14}
M. Pronk and M. Veraar, Tools for Malliavin calculus in UMD Banach spaces, Potential Anal. 40 (2014), no.~4, 307--344. \href{https://doi.org/10.1007/s11118-013-9350-0}{doi:10.1007/s11118-013-9350-0}.

\bibitem{vN10}
J. van Neerven, $\gamma$-Radonifying operators---a survey, in {\it The AMSI--ANU Workshop on Spectral Theory and Harmonic Analysis}, A. Hassell, A. McIntosh, and R. Taggart (eds.), Proc. Centre Math. Appl. Austral. Nat. Univ., vol. 44, Australian National University, Canberra, 2010, 1--62.

\bibitem{WQ15}
X. Wang and R. Qi, A note on an accelerated exponential Euler method for parabolic SPDEs with additive noise, Appl. Math. Lett. 46 (2015), 31--37. \href{https://doi.org/10.1016/j.aml.2015.02.001}{doi:10.1016/j.aml.2015.02.001}.

\bibitem{WQJ17}
X. Wang, R. Qi, and F. Jiang, Sharp mean-square regularity results for SPDEs with fractional noise and optimal convergence rates for the numerical approximations, BIT Numer. Math. 57 (2017), no. 2, 557--585. \href{https://doi.org/10.1007/s10543-016-0639-4}{doi:10.1007/s10543-016-0639-4}.

\bibitem{ZHL23}
H. Zhou, Y. Hu, and Y. Liu, Backward Euler method for stochastic differential equations with non-Lipschitz coefficients driven by fractional Brownian motion, BIT Numer. Math. 63 (2023), Paper No.~40, 37 pp. \href{https://doi.org/10.1007/s10543-023-00981-z}{doi:10.1007/s10543-023-00981-z}.

\end{thebibliography}
\end{document}